\documentclass{article}
\usepackage{amsmath}
\usepackage{amsthm}
\usepackage{amscd}
\usepackage{graphicx}
\newtheorem{theorem}{Theorem}[section]
\newtheorem{lemma}{Lemma}[section]
\newtheorem{definition}{Definition}[section]
\newtheorem{proposition}{Proposition}[section]
\newtheorem{corollary}{Corollary}[section]
\newtheorem{remark}{Remark}[section]
\newcommand{\D}{{\bf D}}

\newcommand{\ep}{\varepsilon}

\newcommand{\lan}{\lambda}
\newcommand{\Dim}{\mbox{\rm Dim}}
\newcommand{\N}{{\bf N}}
\newcommand{\R}{{\bf R}}
\newcommand{\Z}{{\bf Z}}
\newcommand{\p}{\partial}
\newcommand{\J}{{\bf J}}

\newcommand{\diminf}{\underline{\delta}_\lambda(x_0)}

\newcommand{\uL}{\underline{L}}
\newcommand{\dimsup}{\overline{\delta}_\lambda(x_0)}

\newcommand{\Diminf}{\underline{\Delta}_\lambda(x_0)}

\newcommand{\ellinf}{{\underline{\ell}}}
\newcommand{\ellsup}{{\overline{\ell}}}
\newcommand{\overtau}{{\overline{\tau}_\lan(x_0)}}

\newcommand{\cark}{\chi_{{}_{\widetilde A_k}}}
\newcommand{\diam}{\text{\rm diam\,}}

\newcommand{\Pibf}{{\bf\Pi}}

\begin{document}

\title{Quantitative recurrence properties of expanding maps}


\author{J.L. Fern\'andez\thanks{Research supported by Grant BFM2003-04780
from Ministerio de Ciencia y Tecnolog\'{\i}a, Spain}
\\ \small Universidad Aut\'onoma de Madrid \\ \small Madrid , Spain \\
\small {\tt joseluis.fernandez@uam.es}
\and M.V. Meli\'an$^*$ \\ \small Universidad Aut\'onoma de Madrid \\ \small Madrid , Spain \\
\small {\tt mavi.melian@uam.es} \and D. Pestana \thanks{Research
supported by Grants BFM2003-04780 and BFM2003-06335-C03-02
Ministerio de Ciencia y Tecnolog\'{\i}a, Spain}
 \\ \small Universidad Carlos III de Madrid \\ \small Madrid, Spain \\
\small {\tt domingo.pestana@uc3m.es}}
\date{}
\maketitle

\abstract{Under a map $T$, a point $x$ recurs at rate given by a sequence $\{r_n\}$ near a point $x_0$ if $d(T^n(x),x_0)< r_n$ infinitely often. Let us fix $x_0$, and consider
the set of those $x$'s. In this paper, we study the size of this set for expanding maps and obtain its measure and sharp lower bounds on its dimension involving  the entropy of $T$, the local dimension near $x_0$ and the upper limit of 
$\frac 1n \log\frac 1{r_n}$. We apply our results in several concrete examples including subshifts of finite type, Gauss transformation and inner functions.
}

\section{Introduction}

The pre-images under a mixing transformation $T$  distribute
themselves somehow regularly along the base space. In this paper
we aim to quantify this regularity by studying both the measure
and dimension of some recurrence sets. More precisely, we study the
behaviour of pre-images under expanding transformations, i.e.
transformations which locally increase distances.

\medskip

Throughout this paper $(X,d)$ will be a locally complete separable
metric space endowed with a finite measure $\lan$ over the
$\sigma$-algebra $\cal A$ of Borel sets. We further assume
throughout that the support of $\lan$ is equal to $X$ and that
$\lan$ is a non-atomic measure.
We recall that a measurable transformation $T:X\longrightarrow X$
preserves the measure $\lan$ if $\lan(T^{-1}(A))=\lan(A)$ for
every $A\in{\cal A}$.
The classical recurrence theorem of Poincar\'e (see, for example,
\cite{Fu}, p.61) says that

\medskip

\noindent {\bf Theorem A  (H. Poincar\'e).} {\it If $T:X \longrightarrow X$ preserves the measure $\lan$,
then $\lan$-almost every point of $X$
is recurrent, in the sense that
$$
\liminf_{n\to \infty} d(T^n(x),x)=0 \,.
$$}
Here and hereafter $T^n$ denotes the $n$-th fold composition $T^n=T\circ T \circ \cdots \circ T$.
M. Boshernitzan obtained in \cite{Bo} the following quantitative version of
Theorem A.

\

\noindent {\bf Theorem B (M. Boshernitzan).} {\it If  the
Hausdorff $\alpha$-measure $H_\alpha$ on $X$ is $\sigma$-finite
for some $\alpha>0$ and  $T:X \longrightarrow X$ preserves the
measure $\lan$, then for $\lan$-almost all $x\in X$,
$$
\liminf_{n\to\infty} n^{1/\alpha} \, d(T^n(x),x) < \infty \,.
$$
Besides, if $H_\alpha (X)=0$,  then for $\lan$-almost all $x\in X$,
\begin{equation} \label{Bos}
\liminf_{n\to\infty} n^{1/\alpha} \, d(T^n(x),x) = 0
\end{equation}
and when the measure $\lan$ agrees with $H_\alpha$ for some $\alpha>0$, then for $\lan$-almost all $x\in X$,
$$
\liminf_{n\to\infty} n^{1/\alpha} \, d(T^n(x),x) \le 1 \,.
$$
}

L. Barreira and B. Saussol  \cite{BS} have obtained a
generalization of (\ref{Bos}) when $X\subseteq {\bf R}^N$ in terms
of the lower pointwise dimension of $\lan$ at the point $x\in X$
instead of the Hausdorff measure of $X$ and other authors have
also obtained new quantitative recurrence results relating various
recurrence indicators with entropy and dimension, see e.g.
\cite{ACS}, \cite{BGI}, \cite{G}, \cite{G2} and \cite{STV}.

\medskip

It is natural to ask if the orbit $\{T^n(x)\}$ of the point $x$
comes back not only to every neighborhood of $x$ itself as
Poincar\'e's Theorem asserts, but whether it also visits every
neighborhood of a previously chosen point $x_0 \in X$. Under the
additional hypothesis of ergodicity it is easy to check that for
any $x_0\in X$, we have that
\begin{equation} \label{poincareergodico}
\liminf_{n\to \infty}  d(T^n(x),x_0) = 0 \,, \qquad \hbox{for $\lan$-almost all $x\in X$.}
\end{equation}
Recall that the
transformation $T$ is {\it ergodic} if the only $T$-invariant sets
(up to sets of $\lan$-measure zero) are trivial, {\it i.e.} they
have zero $\lan$-measure or their complements have zero
$\lan$-measure.

\

In order to obtain a quantitative version of
(\ref{poincareergodico}) along the lines of Theorem B we need
stronger mixing properties on $T$. In \cite{FMP} we studied {\it
uniformly mixing} transformations. For these transformations we
obtained, for example, that, given a decreasing sequence  $\{r_n\}$
of positive numbers tending to zero as $n\to\infty$, if
$$\displaystyle\sum_{n=1}^\infty \lan(B(x_0,r_n)) = \infty,$$
then
$$
\lim_{n\to\infty} \dfrac {\# \{i\le n: \ d(T^i (x),x_0) \le
r_i\,\}} {\sum_{j=1}^n \lan(B(x_0,r_j))} = 1 \,, \qquad \hbox{for
$\lan$-almost every $x\in X$}\,,
$$
and therefore
$$
\liminf_{n\to \infty} \frac {d(T^n(x),x_0)}{r_n} \le 1 \,, \qquad
\hbox{for $\lan$-almost every $x\in X$}\,.
$$
Here and hereafter the notation $\# A$ means the number of elements of the set $A$.

\

\noindent {\bf Expanding maps.}

\medskip

In this paper we consider the recurrence properties of the orbits under expanding maps, a context which encompasses many interesting examples:  subshift of finite
type, in particular Bernoulli shifts, Gauss transformation and
continued fractions expansions, some inner functions and expanding
endomorphisms of compact manifolds.

An expanding map  $T:X\longrightarrow X$ does not, in general,
preserves the given measure  $\lan$ for a given expanding system, but among all the measures invariant
under $T$ there exists a unique probability measure $\mu$ which is
locally absolutely continuous with respect to $\lan$ and has good
mixing properties (see Theorem E). We will refer to $\mu$ as the
absolutely continuous invariant probability measure (ACIPM). For
the complete definition of expanding maps we refer to Section
\ref{expandingmaps}. However we will describe here their main
properties in an informal way. An expanding system $(X,d,\lan,T)$
has  an associated  Markov partition ${\cal P}_0$ of $X$ in such a way that
$T$ is injective  in each block $P$ of ${\cal P}_0$
and $T(P)$ is a union of blocks of ${\cal P}_0$. Also there
exists a positive measurable function ${\bf J}$ on $X$, the {\it
Jacobian} of $T$, such that
$$
\lan(T(A)) = \int_A {\bf J} \, d\lambda\,, \qquad  \hbox{if $A$ is contained in some block of ${\cal P}_0$}
$$
and ${\bf J}$  has the following distortion property for all $x,y$ in the same block of ${\cal P}_0$:
$$
\left|  \frac {\J(x)}{\J(y)}  -1 \right| \le C_1\,
d(T(x),T(y))^\alpha \,.
$$
Here $C_1>0$ and $0<\alpha\le 1$ are absolute constants. This property allows us to compare the ratio $\lan (A)/\lan(A')$ with  $\lan (T(A))/\lan(T(A'))$for $A,A'$ contained in the same block of ${\cal P}_0$.

Finally, the reason for the name  expanding is the property that if $x,y$ belong
to the same block of the partition ${\cal P}_n=\vee_{j=0}^{n}
T^{-j}({\cal P}_0)$, then
$$
d(T^n(x),T^n(y))\ge C_2 \beta^n d(x,y)
$$
with absolute constants $C_2>0$ and $\beta>1$.

In this paper we are interested in studying for expanding systems the size of the set of
points of $X$ where $\liminf_{n\to\infty} d(T^n(x),x_0)/r_n =0$, where $\{r_n\}$ is a given sequence of positive numbers and $x_0$
is a previously chosen point in $X$. To do this we study the size of the set
$$
{\cal W}(x_0,\{r_n\}) = \{x\in X: \ d(T^n (x),x_0)< r_n \
\mbox{for infinitely many} \ n\}\,.
$$
Our first objective is to study the relationship between the measure of this set and how fast goes to zero  the sequence  of radii.
We use  the following definitions of local dimension.

\begin{definition} \label{dimensiones}
The lower and upper ${\cal P}_0$-dimension of the measure $\mu$ at the
point $x\in X$ are defined, respectively, by
$$
\underline{\delta}_\mu(x) = \liminf_{n\to\infty} \frac{\log \mu(P(n,x))}{\log\hbox{\rm
diam}(P(n,x))}\,, \qquad
\overline{\delta}_\mu(x) = \limsup_{n\to\infty} \frac{\log \mu(P(n,x))}{\log\hbox{\rm
diam}(P(n,x))}\,.
$$
Here and hereafter $P(n,x)$ denotes the block of the partition ${\cal
P}_n=\vee_{j=0}^{n} T^{-j}({\cal P}_0)$ which contains the point
$x\in X$.
\end{definition}

We have  the following result.

\begin{theorem}
\label{intromedidaradios} Let $(X,d,\lan,T)$ be an expanding system with
finite entropy $H_\mu({\cal P}_0)$ with respect to the partition
${\cal P}_0$.
Let $\{r_n\}$ be a non increasing sequence of positive numbers.
Then for $\lan$-almost all point $x_0\in X$ we have that
$$
\hbox{if} \qquad\sum_{n=1}^\infty r_n^{\delta} = \infty \, \quad
\hbox{for some $\delta>\overline{\delta}_\lambda(x_0),\quad $ then
\quad ${\cal W}(x_0,\{r_n\})$ has full $\lan$-measure}\,,
$$
and we can conclude that, for $\lan$-almost all point
$x_0\in X$,
$$
\liminf_{n\to\infty} \frac {d(T^n(x),x_0)}{r_n} = 0 \,, \qquad
\hbox{for $\lan$-almost all $x\in X$}\,.
$$
\end{theorem}

In particular, we have, for $\alpha>\dimsup$, that for $\lan$-almost all point $x_0\in X$
$$
\liminf_{n\to\infty} n^{1/\alpha} d(T^n(x),x_0) = 0 \,, \qquad
\hbox{for $\lan$-almost all $x\in X$}\,.
$$

Further results about the points $x_0$ which satisfy the conclusions in Theorem \ref{intromedidaradios} are included in Section \ref{medida}. There is also a quantitative
version when the system has the Bernoulli property (see Theorem
\ref{elbueno}).

\smallskip

If the sequence $\{r_n\}$ tends to zero in such a way that $\sum_n
\lan(B(x_0,r_n))<\infty$ then it is easy to check that $\lan({\cal
W}(x_0,\{r_n\}))=0$ and therefore
$$
\liminf_{n\to\infty} \frac {d(T^n(x),x_0)}{r_n} \ge 1 \,, \qquad
\hbox{for $\lan$-almost all $x\in X$}\,.
$$
As a consequence we get that if $\lan(B(x_0,r))\le C\,r^\Delta$ for all $r$, then for $\alpha<\Delta$
$$
\liminf_{n\to\infty} n^{1/\alpha} d(T^n(x),x_0) = \infty \,, \qquad
\hbox{for $\lan$-almost all $x\in X$}\,.
$$
\indent Theorem \ref{intromedidaradios} is sharp in the sense that if the series $\sum_{n=1}^\infty r_n^{\delta}$ diverges for $\delta<\overline{\delta}_\lambda(x_0)$, then it can happen that the set ${\cal W}(x_0,\{r_n\})$ has zero $\lan$-measure. Consider, for instance, the sequence $r_n=1/n^{1/\delta}$ with $\delta<1$, and an expanding system with $X\subset\R$ and $\lan$ the Lebesgue measure. Then 
$\dimsup=1$ and $\sum_n r_n^\delta=\infty$. But $\sum_n 
\lan(B(x_0,r_n))<\infty$ and therefore $\lan({\cal W}(x_0,\{r_n\}))=0$.
On the other hand,
for expanding systems with extra mixing properties, Theorem \ref{intromedidaradios} still holds for
$\delta=\overline{\delta}_\lambda(x_0)$ (see, Theorem 3 in \cite{FMP}).

\smallskip

Even in the case that the set ${\cal W}(x_0,\{r_n\})$ has zero $\lan$-measure
we have proved that this set is large since we have
obtained a positive lower bound for its dimension.

\smallskip

In this paper we use two different notions of dimension:
the $\lan$-grid dimension  ($\Dim_{\Pibf,\lan}$),  considering
coverings with blocks of the partitions ${\cal P}_n$, and the
$\lan$-Hausdorff dimension  ($\Dim_\lan$), when we consider
coverings with balls of small diameter (see Section
\ref{grids} for the definitions). We remark that $\Dim_\lan$  is equal to $1/N$ times the usual Hausdorff dimension when $\lan$ is the Lebesgue measure in $X=\R^N$. To obtain lower bounds for the
dimension, we have constructed a Cantor-like set contained in
${\cal W}(x_0,\{r_n\})$. The elements of the different families of
the Cantor set are certain blocks of some partitions ${\cal P}_n$.
In our construction these blocks have controlled $\mu$-measure and
they are in a certain sense well distributed. The main difficulty while estimating the dimension of this Cantor set  is  that does not
have a fixed pattern and the ratio between  the measure of  a
`parent' and his `son' can be very big depending on the sequence
of radii. Our approach is contained in Theorem \ref{azulesyrojos}.

The main tools in the construction of the Cantor set are: (1) good estimates for
the measure $\mu$ of some blocks of ${\cal P}_n$, obtained as a consequence of
Shannon-McMillan-Breimann Theorem  (see Theorem D); (2)
good estimates of the ratio between $\lan(P(n+1,x))$ and
$\lan(P(n,x))$ due to the distortion property of $\bf J$. An
extra difficulty is that the measures $\lan$ and $\mu$ are only
comparable in each block of the partition ${\cal P}_0$.

\

In order to obtain lower bounds for $\Dim_\lan$ we relate it with
$\Dim_{\Pibf,\lan}$ and to do so we have required an  extra
condition of regularity over the `grid' $\Pibf=\{{\cal P}_n\}$
(see Section \ref{grids}). The required condition is trivially
fulfilled in the one dimensional case.  When the
measure $\lan$ of a ball is comparable to a power of its diameter we have
also obtained an estimate of $\Dim_{\lan}$ without assuming
the regularity condition on the grid.

\

\begin{theorem}  \label{introdimensionradios}
Let $(X,d,\lan,T)$ be an expanding system with finite
entropy $H_\mu({\cal P}_0)$ with respect to the partition ${\cal
P}_0$,  and let us  us consider the grid $\Pibf=\{{\cal P}_n\}$. Let $\{r_n\}$ be a non increasing sequence of positive numbers, and let $U$ be
an open set  in $X$ with $\mu(U)>0$. Then, for almost all $x_0\in X$,
$$
\frac 1{\Dim_{\Pibf,\lan} ({\cal W}(x_0,\{r_n\})\cap U)} -1  \le  \frac {\dimsup \, \ellsup}{h_\mu} \,,
$$
where $\ellsup= \limsup_{n\to\infty} \frac 1n\log\frac{1}{r_n}$
and $h_\mu$ is the entropy of $T$ with respect to $\mu$.

\medskip

\noindent Moreover, for almost all $x_0\in X$, the Hausdorff dimensions of the set
${\cal W}(U,x_0,\{r_n\})$ verify:

\smallskip

\begin{enumerate}

\item[\rm 1.] If the grid $\Pibf$ is $\lan$-regular then
$$
\frac 1{\Dim_\lan  ({\cal W}(U,x_0,\{r_n\}))} -1 \le \frac {\dimsup \, \ellsup}{h_\mu} \,.
$$
\item[\rm 2.] If $\lan$ is a doubling measure verifying that $\lan(B(x,r))\le C\,r^s$ for all ball $B(x,r)$, then
$$
\Dim_\lan  ({\cal W}(U,x_0,\{r_n\})) \ge 1 - \frac {\dimsup\ellsup}{s\log\beta}\,.
$$
\end{enumerate}
As a consequence, we obtain the same estimates for the Hausdorff dimensions of the set
$$
 \left\{x\in U: \ \liminf_{n\to\infty}  \frac{d(T^n(x),x_0)}{r_n}=0\right\}.
$$
\end{theorem}

Observe that, for instance, if $X\subset\R$, we obtain that, for any $\alpha>0$,
$$
\Dim_\lan\{x\in U: \ \liminf_{n\to\infty}  e^{n\alpha} d(T^n(x),x_0)=0\}\ge  \frac {h_\mu}{h_\mu +
\dimsup \, \alpha} \, .
$$

Theorem \ref{introdimensionradios} is sharp since for some expanding systems we get equality, see Theorem \ref{innerp}.
As in Theorem \ref{intromedidaradios} we have chosen to state  the above result for  almost all  $x_0\in X$ and we refer to Section \ref{lowerdimension} for more precise results concerning to the set of points $x_0$ where this kind of results holds.

\

Results related to these two theorems above can be found in
\cite{Ja}, \cite{Be}, \cite{BaSc}, \cite{Su}, \cite{MP},
\cite{DMPV}, \cite{CK} and \cite{KM}. See also, \cite{Da}, \cite{FM},
\cite{BiJo}, \cite{Ur}  \cite{Do} and \cite{FM2}.

\

\noindent {\bf Coding.}

\medskip

It is a well known fact that an expanding map induces a coding on the points of $X$ (see  Section \ref{codes}). Via this coding the above results are in certain sense a consequence of analogous results involving symbolic dynamic. More precisely, each point $x$ of the set
$$
X_0: = \bigcap_{n=0}^\infty \bigcup_{P\in {\cal P}_n} P
$$
can be codified as $x=[\;i_0 \;i_1\:\ldots\;]$ where
$
P(0,T^n(x))=P_{i_n}\in{\cal P}_0$   for all $ n=0,1,2,\ldots$
Notice that if $x=[\;i_0 \; i_1 \; i_2 \ldots\;]$ then $T(x)=[\;i_1\:i_2\;i_3\;\ldots\;]$, i.e. $T$ acts as the left shift on the set of all codes.

Given an increasing sequence $\{t_k\}$ of positive integers and a point
$x_0\in X_0$ we  study the size of set
$$
\widetilde{\cal W}(x_0,\{t_n\}) = \{x\in X_0: \ T^k(x) \in
P(t_k,x_0) \ \hbox{for infinitely many $k$}\} .
$$
If $x\in X_0$ and $T^k(x)\in P(t_k,x_0)$ then $P(j, T^k(x))=P(j,x_0)$ for $j=0,1,\ldots, t_k$ and it follows that $\widetilde{\cal W}(x_0,\{t_n\})$
can be also described as the set of points
$x=[\;m_0\;m_1\ldots\;]\in X_0$ such that
$$
m_k=i_0,\, m_{k+1}=i_1,\;\dots\;,\, m_{k+t_k}=i_{t_k}
$$
for infinitely many $k$, where $x_0=[\;i_0\;i_1\ldots\;]$. For this set, we have the following analogue of Theorem \ref{intromedidaradios}:

\begin{theorem}  \label{intromedidacode}
Let $(X,d,\lan,T)$ be an expanding system. Let
$x_0$ be a  point of $X_0$ such that
$\underline{\delta}_\lambda(x_0)>0$ and let  $\{t_n\}$ be a non
decreasing sequence of positive integers numbers.
$$
\hbox{If} \qquad \sum_{n=1}^\infty \lan(P(t_n,x_0)) = \infty, \qquad \hbox{then}
\qquad \lan(\widetilde{{\cal W}}(x_0,\{t_n\}))=\lan (X).
$$
Moreover, if the partition ${\cal P}_0$ is finite or if the system
has the Bernoulli property, $($i.e. if\, $T(P)=X$ $($mod $0)$ for
all $P\in {\cal P}_0)$, then we have the following quantitative
version:
\begin{equation} \label{cuantitativointro}
\lim_{n\to\infty} \frac {\#\{i\le n: \ T^i(x)\in
P(t_i,x_0)\}}{\sum_{j=1}^n \mu(P(t_j,x_0))} =1\,,  \qquad
\hbox{for $\lan$-almost every $x\in X$ , }
\end{equation}
where $\mu$ is the ACIPM associated to the system.
\end{theorem}

Property (\ref{cuantitativointro}) is related to  the  decay of
the correlation coefficients of the indicator functions of
$\{P(n,x_0)\}$, see \cite{S} and \cite{Ph}. For
expanding systems with the Bernoulli property L.S. Young \cite{Y}
has proved that this decay is exponential.

\smallskip

As with ${\cal W}(x_0,\{r_n\})$, it is easy to see using the Borel-Cantelli lemma  that if
$\sum_n\lan(P(t_n,x_0)) < \infty$, then $\lan(\widetilde{{\cal W}}(x_0,\{t_n\}))=0$.

\begin{theorem} \label{introdimensioncode} Let $(X,d,\lan,T)$ be an expanding system with finite entropy $H_\mu({\cal P}_0)$ with respect to the
partition ${\cal P}_0$ where $\mu$ is ACIPM associated to the
system, and let us consider the grid $\Pibf=\{{\cal P}_n\}$.  Let $\{t_n\}$ be a non decreasing sequence of positive
integers and let $U$ be an open set in $X$ with $\mu(U)>0$. Then,
for almost all point $x_0\in X_0$,
$$
\frac 1{\Dim_{\Pibf,\lan} ({\cal \widetilde W}(x_0,\{t_n\})\cap U)}-1 \le \frac
{\overline{L}(x_0)}{h_\mu}  \,,
$$
where $\overline{L}(x_0)= \limsup_{n\to\infty} \frac 1n
\log\frac{1}{\lan(P( t_n,x_0))}$ and $h_\mu$ is the entropy of $T$
with respect to $\mu$. Moreover, if the grid $\Pibf$ is
$\lan$-regular, then
$$
\frac 1{\Dim_{\lan} ({\cal \widetilde W}(x_0,\{t_n\})\cap U)} \le \frac {\overline{L}(x_0)}{h_\mu} \,.
$$
\end{theorem}

As in the previous theorems,  for the sake of simplicity, we have stated this last result for almost all point $x_0$ in $X$, but we refer to Section \ref{lowerdimension} for a more precise statement concerning the points $x_0$ which satisfy the conclusions.  Also, we should mention that, up a $\lan$-zero measure set, we have that $\overline{L}(x_0)/h_{\mu}=\limsup_{n\to\infty} t_n/n$ (see Theorem \ref{dimensioncode2}).

\

\

\noindent {\bf Applications.}

\medskip

The generality of the definition of expanding systems allows us to
apply our results in a broad kind of situations. In the final
section, we  have obtained results for Markov transformations,
subshifts of finite type (in particular, Bernoulli shifts), the
Gauss transformation, some inner functions and expanding
endomorphisms. In the case of Bernoulli shifts we also give a
precise upper bound of the dimension by using a large deviation
inequality. As an example, for the Gauss map $\phi$, which acts on
the continued fractions expansions as the left shift,  we have the
following results:

\begin{theorem}
\begin{itemize}
\item[]
\item[\rm (1)] If $\alpha>1$ then, for almost all $x_0\in [0,1]$, and more precisely,
if $x_0$  is an irrational number with continued fraction expansion $[\,i_0\,,\;i_1\,, \; \dots ]$ such that
$\log i_n = o(n)$ as $n\to\infty$, we have that
$$
\liminf_{n\to\infty} n^{1/\alpha}  |\phi^n(x)-x_0| = 0 \,, \qquad
\text{for almost all $x\in [0,1]$.}
$$

\item[\rm (2)] If $\alpha <1$, then for all $x_0\in [0,1]$ we have that
$$
\liminf_{n\to\infty} n^{1/\alpha} |\phi^n(x)-x_0| = \infty \,,
\qquad \text{for almost all $x\in [0,1]$.}
$$

\item[\rm (3)] If $x_0$ verifies the same hypothesis than in part $(1)$, then
$$
\Dim \left\{x\in [0,1]: \ \liminf_{n\to\infty} n^{1/\alpha}
|\phi^n(x)-x_0| =0 \right\} = 1\,, \qquad \text{for any
$\alpha>0$}.
$$
and
$$
\Dim \left\{x\in [0,1]: \ \liminf_{n\to\infty} e^{n\kappa}
|\phi^n(x)-x_0| =0 \right\} \ge \frac {\pi^2}{\pi^2+6\kappa\log2}
\,, \qquad \text{for any $\kappa>0$}.
$$
\end{itemize}
\end{theorem}

\begin{theorem}
Let $x_0\in [0,1]$ be an irrational number with continued fraction expansion $x_0=[\,i_0,i_1,\dots \,]$ and let $t_n$ be a non decreasing sequence of natural numbers. Let $\widetilde W$ be the set of
points $x=[\,m_0,m_1,\dots\,]\in [0,1]$ such that
$$
m_n=i_0\,, \; m_{n+1}=i_1\,,\; \dots\;, \; m_{n+t_n} = i_{t_n} \;,
\qquad \text{for infinitely many $n$}.
$$

\begin{itemize}
\item[\rm (1)] \ $\lan(\widetilde W)=1$, if
$$
 \sum_n \frac 1{(i_0+1)^2\cdots (i_{t_n}+1)^2} =\infty \,.
$$

\item[\rm (2)] \ $\lan(\widetilde W)=0$, if
$$
\sum_n \frac 1{i_0^2\cdots i_{t_n}^2} <\infty \,.
$$

\item[\rm (3)] In any case, if  $\log i_n = o(n)$ as $n\to\infty $, then
$$
\Dim (\widetilde W) \ge \frac {\pi^2}{\pi^2+6\log 2\,\limsup_{n\to\infty}
\frac 1n \log  (i_0+1)^2\cdots (i_{t_n}+1)^2} \,.
$$
\end{itemize}
\end{theorem}

\

The techniques developed in this paper for expanding maps and therefore for one-sided Bernoulli shifts, can be extended to bi-sided Bernoulli shifts. This has allowed us \cite{Anosov} to get results on recurrence for Anosov flows.

\

The outline of the paper is as follows: In Section
\ref{grids} we give our two definitions of dimension and prove some general results for computing the  dimensions of a kind of
Cantor-like sets with the particular feature that the ratio between the size of a 'son' and his 'parent'  decays very fast. Section \ref{SMB} contains
some consequences of Shannon-McMillan-Breiman Theorem.
 In section \ref{expandingmaps} we give the
complete definition of an expanding system and  in Section
\ref{codes} we recall how to associate a code to the points of
$X$. In  Section \ref{PropEM} we prove some general properties of
expanding maps.  The precise statements and proofs of Theorems
\ref{intromedidaradios} and \ref{intromedidacode}, and some
consequences of them, are contained in Section \ref{medida}. The dimension results are included in Section
\ref{dimensionresults}. More general versions of Theorems
\ref{introdimensionradios} and \ref{introdimensioncode} are
included in Section \ref{lowerdimension}.  In
Section \ref{upperbounds} we include  an upper bound of the
dimension. Finally, Section \ref{aplicaciones} contains several
applications of the above results.

\

\noindent {\it Acknowledgements:} We want to thank  to  J. Gonzalo,   R. de la Llave, V. Mu\~noz and R. P\'erez Marco 
 for helpful conversations about this work.  We are particularly indebted to A. Nicolau for his encouragement and stimulating discusions.

\

\noindent {\it  A few words about notation.}  There are many estimates in this paper involving absolute constants. These are usually denoted by capital letters like $C$. Occasionally, we shall indicate a constant $C$ depending on some parameter $\alpha$ by $C(\alpha)$. The symbol $\# D$ denotes the number of elements of the set $D$. By  $A\asymp B$ we mean that there exist absolute constants $C_1,C_2>0$ such that $ C_1B\leq A\leq C_2 B$.

\section{Grids and dimensions.} \label{grids}

\

Along this section $(X,d,{\cal A},\lan)$ will be a finite measure
space with a compatible metric. Compatible means that $\cal A$ is the the
$\sigma$-algebra of the Borel sets of $d$. We recall that we are assuming that the
measure $\lan$ is non-atomic and its support is $X$.

\begin{definition}
Given a set $F\subset X$ and $0<\alpha\leq 1$, we define the
$\alpha$-dimensional $\lan$-Hausdorff measure of $F$ as
$$
{\cal H}_{\lan}^{\alpha}(F)=\lim_{\ep\to 0}{\cal
H}^{\alpha}_{\lan,\, \ep}(F)
$$
with
$$
{\cal H}^{\alpha}_{\lan,\, \ep}(F)=\inf\sum_i(\lan(B_i))^{\alpha}
$$
where the infimum is taken over all the coverings $\{B_i\}$ of $F$
with balls such that $\diam(B_i)\leq\ep$ for all $i$.
\end{definition}

It is not difficult to check that ${\cal H}^{\alpha}$ is a regular Borel
measure, see e.g. \cite{Matti}. Observe that if $X\subset \R^N$ and $\lan$ is Lebesgue
measure, then ${\cal H}^{\alpha}$ is comparable with the usual
$N\alpha$-dimensional Hausdorff measure.

\begin{definition}
The $\lan$-Hausdorff dimension of $F$ is defined as
$$
{\mbox{\rm Dim}}_{\lan}(F)=\inf\{\alpha \, : \, {\cal
H}_{\lan}^{\alpha}(F)=0\}=\sup\{\alpha \, : \, {\cal
H}_{\lan}^{\alpha}(F)>0\}\,.
$$
\end{definition}

If $X\subset \R^N$ and $\lan$ is Lebesgue measure, then the
$\lan$-Hausdorff dimension coincides with $1/N$ times the usual
Hausdorff dimension.

\begin{definition}
A grid is a collection $\Pibf=\{{\cal P}_n\}$ of partitions of $X$
each of them constituted by disjoint open sets, and such
that for all $P_n\in {\cal P}_n$ there exists a unique
$P_{n-1}\in{\cal P}_{n-1}$ such that $P_n\subset P_{n-1}$, and
$\sup_{P\in{\cal P}_n} \diam(P) \to 0$ as $n\to\infty$..
\end{definition}

\begin{definition}
Given a grid $\Pibf=\{{\cal P}_n\}$ of $X$ and  $0<\alpha\leq 1$,
the $\alpha$-dimensional $\lan$-grid measure of any subset
$F\subset X$ is defined as
$$
{\cal H}_{\Pibf,\lan}^{\alpha}(F)=\lim_{n\to \infty}{\cal
H}^{\alpha}_{\Pibf,\lan, n}(F)
$$
with
$$
{\cal H}^{\alpha}_{\Pibf,\lan, n}(F)=\inf\sum_i(\lan(P_i))^{\alpha}
$$
where the infimum is taken over all the  coverings $\{P_i\}$ of
$F$ with sets  $P_i\in\cup_{k\geq n}{\cal P}_k$.

The $\lan$-grid Hausdorff dimension of $F$ is defined as
$$
\Dim_{\Pibf,\lan}(F)=\inf\{\alpha \, : \, {\cal
H}_{\Pibf,\lan}^{\alpha}(F)=0\}=\sup\{\alpha \, : \, {\cal
H}_{\Pibf,\lan}^{\alpha}(F)>0\}\,.
$$
\end{definition}

As before we have that ${\cal H}_{\Pibf,\lan}^{\alpha}$ is a Borel measure.

\begin{remark} \rm
If $X\subseteq \R$ and $\lan$ is Lebesgue measure we have that ${\cal H}_{\lan}^{\alpha}(F)\le {\cal
H}_{\Pibf,\lan}^{\alpha}(F)$ and therefore, for any $F\subset \R$,
$$
\Dim_\lan (F) \le \Dim_{\Pibf,\lan}(F) \,.
$$
Also, if $X=[0,1]$, ${\cal P}_n$ denotes the family of dyadic intervals with length $1/2^{n+1}$ and $\lan$ is Lebesgue measure, then we have, for any $F\subset [0,1]$, that
$$
\Dim_\lan (F) = \Dim_{\Pibf,\lan}(F) \,.
$$
\end{remark}

In order to compute the $\lan$-grid Hausdorff dimension we will use  the
following result which parallels Frostman lemma.

\begin{lemma}\label{DimBill}
Let $\Pibf=\{{\cal P}_n\}$ be a grid of $X$. For each $n\in\N$, let
${\cal Q}_n$ be a subcollection of ${\cal P}_n$ and let $F$ be a set with
$$
F \subseteq \bigcap_{n}\bigcup_{Q\in{\cal Q}_n} Q \,.
$$
If there exist a measure $\nu$ such that $\nu(F)>0$, a real
number $0<\gamma\leq 1$  and a positive constant $C$ such that, for all $x\in F$,
$$
\nu(Q(k,x)) \le C\, (\lambda(Q(k,x)))^\gamma \,,
$$
where $Q(k,x)$ denotes the block of ${\cal Q}_k$ which contains $x$, then,
$$
\Dim_{\Pibf,\lan}(F)\geq \gamma \,.
$$
\end{lemma}

\begin{proof}
It follows from the fact that ${\cal H}_{\Pibf,\lan}^{\alpha}(F)\ge \nu(F)/C>0$.
\end{proof}

The following result allows to obtain a lower bound for  the $\lan$-grid Hausdorff dimension of Cantor-like sets.

\begin{theorem}\label{azulesyrojos}
Let $\Pibf=\{{\cal P}_n\}$ be a grid of $X$ and let $\{d_j\}$ and $\{\widetilde{d}_j\}$ be two increasing sequences of natural numbers tending to infinity verifying that
$d_{j-1}<\widetilde{d}_j < d_j $ for each $j$.

\smallskip

Consider two collections $\{{\cal J}_j\}$ and   $\{\widetilde{{\cal J}}_j\}$ of  subsets of $X$ such that:
\begin{itemize}
\item[$(i)$] $\widetilde{\cal J}_0={\cal J}_0=\{J_0\}$ and for each $j$, $
{\cal J}_j \subseteq {\cal P}_{d_j}$ and $\widetilde{\cal J}_j \subseteq {\cal P}_{\widetilde{d}_j}$.
\item[$(ii)$] For each $J_j\in {\cal J}_j$ $(j\ge 1)$ there exists a unique $\widetilde{J}_j \in
\widetilde{\cal J}_j$ such that  $\hbox{\rm closure} (J_j) \subset \widetilde{J}_j$. Reciprocally, for each  $\widetilde{J}_j \in \widetilde{\cal J}_j$ there exists a unique ${J}_j \in \widetilde{\cal J}_j$ such that  $\hbox{\rm closure} (J_j) \subset \widetilde{J}_j$.
\item[$(iii)$] For each $\widetilde{J}_j \in \widetilde{\cal J}_j$ $(j\ge 1)$ there exists a unique
$J_{j-1}\in {\cal J}_{j-1}$ such that $\widetilde{J}_j \subset J_{j-1}$.
\end{itemize}

Let ${\cal C}$ be the Cantor-like set defined by
$$
{\cal C}= \bigcap_{j=0}^\infty \bigcup_{J_j\in{\cal J}_j} J_j=  \bigcap_{j=0}^\infty \bigcup_{\widetilde{J}_j\in\widetilde{\cal J}_j}\widetilde{ J}_j \;.
$$

Assume that   the pattern of ${\cal C}$ has the following additional properties:
\begin{itemize}
\item[$(1)$] There exist two sequences $\{\alpha_j\}$ and $\{\beta_j\}$ of positive numbers such that
$$
\alpha_j \le \displaystyle\frac{\lan(\widetilde{J}_j)}{\lan(J_{j-1})} \le \beta_j \,.
$$
\item[$(2)$] There exists a sequence $\{\gamma_j\}$ of positive numbers such that
$$
\frac{\lan(J_j)}{\lan(\widetilde{J}_{j})} \ge \gamma_j \;.
$$
\item[$(3)$] There exists a sequence $\{\delta_j\}$ with $0<\delta_j\leq 1$ such that
$$
\lan (\widetilde{\cal J}_j \cap J_{j-1}) \ge \delta_j\, \lan(J_{j-1})\,.
$$
\item[$(4)$] There exists an absolute constant $\Lambda$ such that for all $j$ large enough
$$
\frac{1}{\delta_{j+1}}\frac {\beta_j}{\delta_j} \frac {\beta_{j-1}}{\delta_{j-1}} \cdots \frac {\beta_1}{\delta_1} \le
[( \alpha_j \gamma_j)\, ( \alpha_{j-1} \gamma_{j-1})\,\cdots\,( \alpha_1 \gamma_1)]^\Lambda\,.
$$
\end{itemize}

Then
$$
\Dim_{\Pibf,\lan} ({\cal C}) \ge \Lambda\,.
$$
\end{theorem}

\begin{remark} \rm
Observe that in the special case when the two families $\widetilde{\cal J}_j$ and ${\cal J}_j$ coincide and $\alpha_j=\alpha$, $\beta_j=\beta$, $\delta_j=\delta$, then the above result is the usual Hungerford's Lemma $(\Lambda=\log(\beta/\delta)/\log\alpha)$, see e.g. \cite{Pom2}.
\end{remark}

\begin{proof}
We construct a probability measure $\nu$ supported on  ${\cal C}$
in the following way: We define $\nu(J_0)=1$ and for each set
$J_j\in {\cal J}_j$ we write
$$
\nu(J_j)=\nu(\widetilde{J}_j) = \frac {\lan(\widetilde{J}_j)}{\lan(\widetilde{\cal J}_j
\cap J_{j-1})} \, \nu(J_{j-1})
$$
where $J_{j-1}$ and $\widetilde{J}_j$ denote the unique sets in
${\cal J}_{j-1}$ and $\widetilde{\cal J}_{j}$ respectively, such
that $J_j \subset \widetilde{J}_j \subset J_{j-1}$. As usual, for
any Borel  set $B$, the $\nu$-measure of $B$ is defined by
$$
\nu(B) = \nu(B\cap {\cal C}) = \inf  \sum_{U\in {\cal U}}  \nu(U) \,
$$
where the infimum is taken over all the coverings $\cal U$ of $B\cap {\cal C}$
with sets in $\bigcup {\cal J}_j$.

We will show that there exists a positive constant $C$ such that for all $x\in{\cal C}$ and $m$ large enough,
\begin{equation} \label{frostmangrid}
\nu(P(m,x))\leq C\,(\lan(P(m,x))^\Lambda
\end{equation}
and therefore, from Lemma \ref{DimBill}, we get  the result.

To prove (\ref{frostmangrid}) let us suppose first that $P(m,x)=J_j$ for some $J_j\in{\cal
J}_j$. From properties (1)-(3) we have that
\begin{align}
\nu(J_j) & \le \frac {\beta_j}{\delta_j} \, \nu(J_{j-1}) \le \frac {\beta_j}{\delta_j} \frac {\beta_{j-1}}{\delta_{j-1}} \cdots \frac {\beta_1}{\delta_1} \,, \notag \\
\lan(J_j) & \ge \alpha_j \gamma_j \, \lan(J_{j-1}) \ge ( \alpha_j \gamma_j)\, ( \alpha_{j-1} \gamma_{j-1})\,\cdots\,( \alpha_1 \gamma_1) \lan(J_0) \,.\notag
\end{align}
and  follows from property (4) that
\begin{equation} \label{nuj}
\nu(\widetilde{J}_j)=\nu({J}_j) \le C\, \delta_{j+1}\lan(J_j)^\Lambda \,.
\end{equation}
This condition is stronger than (\ref{frostmangrid}) for $\widetilde{J}_j$ and $J_j$  and we  will use it to get (\ref{frostmangrid}) in general.

Now,  let us suppose that $P(m,x)\not= J_j$ for all $j$ and  for
all $J_j\in{\cal J}_j$. Since $x\in{\cal C}$ there exist
$J_j\in{\cal J}_j$ and $J_{j+1}\in{\cal J}_{j+1}$ such that
$$
J_{j+1}\subset P(m,x)\subset J_j \,.
$$
If $P(m,x)\subset \widetilde {J}_{j+1}$, then from the
definition of $\nu$ and (\ref{nuj}) for $\widetilde{J}_{j+1}$ we get
$$
\nu(P(m,x))=\nu(\widetilde{J}_{j+1})=\nu(J_{j+1})\leq C
\,(\lan(J_{j+1}))^{\Lambda}\leq C\, (\lan(P(m,x)))^{\Lambda}\,.
$$

Otherwise $P(m,x)$ contains sets of the family ${\widetilde{\cal
J}}_{j+1}$ and we have that
\begin{align} \label{nocantor}
\nu(P(m,x)) & =
\sum_{\substack{{\tilde{J}}_{j+1} \in {\widetilde{\cal J}}_{j+1} \\
{\widetilde{J}}_{j+1} \subseteq P(m,x)}}  \nu (J_{j+1}) \,=
\sum_{\substack{{\tilde{J}}_{j+1} \in {\widetilde{\cal J}}_{j+1} \\
{\widetilde{J}}_{j+1} \subseteq P(m,x)}} \frac{\lan({\widetilde{J}_{j+1}})}{\lan({\widetilde{\cal J}}_{j+1}
\cap J_j)}\,\nu (J_{j}) \notag \\
& = \frac{\nu (J_{j}) }{\lan({\widetilde{\cal J}}_{j+1}\cap J_j)}\,\sum_{\substack{{\tilde{J}}_{j+1}
\in {\widetilde{\cal J}}_{j+1} \\
{\widetilde{J}}_{j+1} \subseteq
P(m,x)}}\lan({\widetilde{J}_{j+1}}) \leq \frac{\nu (J_{j})
}{\lan({\widetilde{\cal J}}_{j+1}\cap J_j)}\ \lan(P(m,x)) \,.
\end{align}
And using  property (3) and (\ref{nuj}) we obtain
that
$$
\nu(P(m,x))\leq \frac{C}{(\lan(J_j))^{1-\Lambda}}\,\lan(P(m,x))
$$
But $\lan(J_j))\geq \lan(P(m,x))$ and  so we get
$$
\nu(P(m,x))\leq C\,(\lan(P(m,x)))^{\Lambda}\,.
$$
\end{proof}

\begin{remark} \rm  Notice that if we define
$$
\nu(J_j) = \frac {\lan({J}_j)}{\lan({\cal J}_j
\cap J_{j-1})} \, \nu(J_{j-1})
$$
then instead of (\ref{nocantor}) we get

\begin{align}
\nu(P(m,x)) & \le
\frac{\nu (J_{j}) }{\lan({{\cal J}}_{j+1}\cap J_j)}\,\sum_{\substack{{\tilde{J}}_{j+1}
\in {\widetilde{\cal J}}_{j+1} \\
{\widetilde{J}}_{j+1} \subseteq
P(m,x)}}\lan({{J}_{j+1}}) \notag \\
& \le
 \frac{\nu (J_{j}) \,\omega_{j+1}}{\lan({{\cal J}}_{j+1}\cap J_j)}\,\sum_{\substack{{\tilde{J}}_{j+1}
\in {\widetilde{\cal J}}_{j+1} \\
{\widetilde{J}}_{j+1} \subseteq
P(m,x)}}\lan(\widetilde{J}_{j+1})
 \leq \frac{1}{\delta_{j+1}} \frac{\nu (J_{j})}{\lan(J_j)}\frac{\omega_{j+1}}{\gamma_{j+1}}
 \ \lan(P(m,x)) \,. \notag
\end{align}
where
$$
\gamma_j\le \frac{\lan(J_j)}{\lan(\widetilde{J}_{j})} \le \omega_j \; .
$$
Hence if $\nu(J_j)\leq C \delta_{j+1} (\lan(J_j))^{\Lambda}$ we get that
$$
\nu(P(m,x))\leq \frac{C}{(\lan(J_j))^{1-\Lambda}}\,\,\frac{\omega_{j+1}}{\gamma_{j+1}}\, \lan(P(m,x))
$$
 and we will need
$$
\frac{\omega_{j+1}}{\gamma_{j+1}}\,\leq\frac{1}{ ( \lan(P(m,x)))^\ep}
$$
in order to get that the dimension is greater than $\Lambda-\ep$. We recall that in this case the upper bound for $\lan(P(m,x))$ is $\lan(J_j)$ and $\lan(J_j)\leq (\omega_j\beta_j)(\omega_{j-1}\beta_{j-1})\cdots (\omega_1\beta_1)$.
\end{remark}

\begin{corollary} \label{azulesyrojos2}
Under the same hypotheses that in Theorem {\rm \ref{azulesyrojos}} we have that if  $\delta_j=\delta>0$ and
$$
\alpha_j = e^{-N_j a}\,, \qquad \beta_j= e^{-N_j b} \,, \qquad \gamma_j = e^{-N_j c}\,,
$$
then
$$
\Dim_{\Pibf,\lan} ({\cal C}) \ge \frac b{a+c} - \frac {\log(1/\delta)}{a+c}\lim_{j\to\infty} \frac {j}{N_1 +\cdots +N_j}\,.
$$
\end{corollary}

The next definition states some kind of regularity on the
distribution of the blocks of the partitions. This property will
allow us to relate the Hausdorff dimension with the grid Hausdorff
dimension.

\begin{definition}\label{gridregular}
Let $\Pibf=\{{\cal P}_n\}$ be a grid of $X$. We will say that $\Pibf$ is $\lan$-regular if
there exists a positive constant $C$ such that for
all ball $B$
$$
\lan(\cup\{P: P\in{\cal P}_n \, , \, P\cap B\not=\emptyset\})\leq
C\, \lan(B)
$$
for all $n$ such that $\sup_{P\in{\cal P}_n}\lan(P)\leq\lan(B)$.
\end{definition}

\begin{remark}  \label{regularreal}
\rm It is clear from the definition  that any grid of
 $X\subset \R$  is $\lan$-regular (we can take $C=3$)
 \end{remark}

\noindent{\bf An example:}
Let $X$ be the square $[0,1]\times[0,1]$ in $\R^2$ and let us denote by  $\lan$ the Lebesgue measure. Consider the grid $\Pibf=\{{\cal P}_n\}$ defined as follows:  the elements of ${\cal P}_0$ are the four open rectangles  obtained by  dividing  the square $[0,1]\times[0,1]$ through the lines  $x=a$ and $y=b$, with $\frac12<b<a<1$; the elements of ${\cal P}_n$ are getting by dividing each rectangle of ${\cal P}_{n-1}$ in four rectangles using the same proportions.
We will see that this is not a regular grid.

Let us consider the ball $B_k$ with diameter $(1-b)^k$ and contained in the square $[0,(1-b)^k]\times[(1-b)^k,1]$.
It is easy to see that $\sup_{P\in{\cal P}_n}\lan(P)=(ab)^n$, and therefore $(ab)^n\leq \lan(B_k)$ implies
$$
n\geq \frac{\log c+2k\log(1/(1-b))}{\log{1/(ab)}}
$$

Therefore, if $\Pibf$ is regular, then  for $n=n(k)= 2k \frac{\log(1/(1-b))}{\log{1/(ab)}}+C$ the quotient
$$
C_k:=\frac{\lan(\cup\{P: P\in{\cal P}_n \, , \, P\cap B_k\not=\emptyset\})}{\lan(B_k)}
$$
has to be bounded. But, it is easy to see that the elements of
${\cal P}_n$ whose closure intersects to  $[0,1]\times\{0\}$ are
rectangles of width $a^n$, and hence,  since $b<a$,
$$
C_k\ge \frac {(1-b)^ka^{n(k)}}{(1-b)^{2k}} \to \infty \qquad \hbox{as $k\to\infty$}.
$$
Therefore, this grid is not regular. On the other hand it is clear that any grid in $X$
whose elements are all squares is regular.

\

The following result gives a lower  bound for  the Hausdorff dimension of Cantor like sets
which are constructed
using a regular subgrid with some control into the quotient between the size of parents and sons.

\begin{proposition} \label{haches}
Let $\Pibf=\{{\cal P}_n\}$ be a grid of $X$ and let $\{{\cal Q}_n\}$
be a $\lan$-regular subgrid of ${\Pibf}$. Let us suppose that there exist strictly
non increasing sequences $\{a_n\}$, $\{b_n\}$ of positive numbers
such that $\lim_{n\to\infty} b_n=0$ and for all $Q\in {\cal Q}_n$
$$
a_n\leq\lan(Q)\leq b_n\,.
$$
Then, for any subset $F \subseteq \bigcap_n\bigcup_{Q\in{\cal Q}_n} Q$,
\begin{equation} \label{Hausmess}
{\cal H}_{\Pibf,\lan}^{\alpha}(F)\leq C\, {\cal
H}_{\lan}^{1-(1-\alpha)\eta}(F)
\end{equation}
for all $\alpha$ and $\eta$ such that
\begin{equation}\label{limgrid}
\limsup_{n\to\infty}\frac{\log(1/a_n)}{\log
(1/b_{n-1})}<\eta < \frac 1{1-\alpha}\,,
\end{equation}
where $C$ is an absolute positive constant. In particular,
\begin{equation} \label{Hausdims}
\frac {1-\Dim_\lan (F)}{1- \Dim_{\Pibf,\lan}(F)} \le
\limsup_{n\to\infty}\frac{\log(1/a_n)}{\log (1/b_{n-1})}\,.
\end{equation}
\end{proposition}

\begin{proof} Let us consider a ball $B$
such that $B\cap F\not=\emptyset$ and let $n=n(B)$ be the smallest
integer such that $b_n\leq\lan(B)$. Then
$$
b_n\leq\lan(B) <b_{n-1}.
$$
We denote by ${\cal Q}(B)$ the collection of elements in ${\cal
Q}_n$  whose  intersection with $B$ is not empty. Then the
collection ${\cal Q}(B)$ is a covering of $B\cap F$, that is
\begin{equation}\label{covering}
B\cap F\subset \bigcup\{Q\, :\, Q\in{\cal Q}(B)\},
\end{equation}
and moreover by the Definition \ref{gridregular} and the election of $n=n(B)$,
$$
\sum_{Q\in{\cal Q}(B)}(\lan(Q))^{\alpha}=\sum_{Q\in{\cal
Q}(B)}\frac1{(\lan(Q))^{1-\alpha}}\,\lan(Q)\leq C\,
\frac{1}{a_n^{1-\alpha}}\, \lan(B) \,.
$$
We may assume that $n$ is large because $\diam(B)$ is small and so
the above inequality and (\ref{limgrid}) imply that
\begin{equation} \label{medidacov}
\sum_{Q\in{\cal Q}(B)}(\lan(Q))^{\alpha}\leq C' (\lan(B))^{1-(1-\alpha)\eta} \,.
\end{equation}
The inequality (\ref{Hausmess}) follows now from (\ref{covering})
and (\ref{medidacov}). To prove (\ref{Hausdims}) let us observe
that we can assume that
$$
\limsup_{n\to\infty}\frac{\log(1/a_n)}{\log (1/b_{n-1})} < \frac 1{1- \rm
Dim_{\Pibf,\lan}(F)}
$$
since in other case (\ref{Hausdims}) is trivial. Let us choose now
$\alpha$ and $\eta$ such that
\begin{equation} \label{etayalpha}
\limsup_{n\to\infty}\frac{\log(b_{n-1}/a_n)}{\log (1/b_n)} <\eta <
\frac 1{1-\alpha}< \frac 1{1-\Dim_{\Pibf,\lan}(F)}.
\end{equation}
Then ${\cal H}_{\Pibf,\lan}^{\alpha}(F)>0$ and by (\ref{Hausmess})
we have also that ${\cal H}_{\lan}^{1-(1-\alpha)\eta}(F)>0$.  Since $\alpha$ y $\eta$ are
arbitrary numbers verifying (\ref{etayalpha}), the ineguality (\ref{Hausdims})
follows.
\end{proof}

\section{Some consequences of Shannon-McMillan-Breiman Theorem} \label{SMB}

Along this section $(X,{\cal A},\mu)$ will be a finite measure space and $T:X\longrightarrow X$ will be a measurable transformation.
A {\it partition} of $X$ is a family ${\cal P}$ of measurable sets with positive
measure satisfying
\begin{itemize}
\item[1.]  If $P_1,P_2 \in {\cal P}$ then $\mu(P_1 \cap P_2)=0$.
\item[2.]  $\mu\left(X\setminus \cup_{P\in{\cal P}} P\right) =0$.
\end{itemize}
It follows from these properties that ${\cal P}$ must be finite or numerable.
The {\it entropy} of a partition ${\cal P}$ is defined as
$$
H_\mu({\cal P}) = \sum_{P\in{\cal P}}  \mu(P) \, \log \frac 1{\mu(P)}\,.
$$
If $T:X\longrightarrow X$ preserves the measure $\mu$, then the {\it entropy of $T$
with respect to the partition ${\cal P}$} is
$$
h_\mu(T,{\cal P}) = \lim_{n\to\infty} \frac 1n \, H_\mu \big( \vee_{j=0}^{n-1} T^{-j} {\cal P} \big).
$$
This limit exists since the sequence in the right hand side is
decreasing. Hence $h_\mu(T,{\cal P})\le H_\mu({\cal P})$.

Finally, the {\it entropy} $h_\mu(T)$ of the endomorphism $T$ is the supremum of
$h_\mu(T,{\cal P})$ over all the partitions ${\cal P}$ of $X$ with
entropy $h_\mu(T,{\cal P})<\infty$.

If the partition ${\cal P}$ is generating, i.e. if
$\vee_{j=0}^\infty T^{-j}({\cal P})$  generates $\cal A$, then, by
the  Kolmogorov-Sinai  Theorem ([M, p. 218-220]), we get

\

\noindent {\bf Theorem C } {\it Let  $(X,{\cal A}, \mu)$ be a
probability space and $T:X\longrightarrow X$ be a measure
preserving transformation.  If  $\cal P$ is a generating partition
of $X$ and  the entropy $H_{\mu}({\cal P})$ is finite, then
$h_{\mu}(T)=h_{\mu}(T,{\cal P})$}.

\

Let $P(n,x)$ denotes the element of the partition $\vee_{j=0}^n
T^{-j}({\cal P})$  which contains the point $x\in X$. It follows
from the definition of partition, that for almost every $x\in X$,
$P(n,x)$ is defined for all $n$. Entropy is a measure of how fast
$\mu(P(n,x))$ goes to zero. The following fundamental result,
which is due to Shannon, McMillan and Breiman,  formalizes this
assertion:

\

\noindent{\bf Theorem D}([M, p. 209]) {\it Let  $(X,{\cal A},\mu)$
be a probability  space and let $T:X\longrightarrow X$ be a
measure preserving ergodic transformation.  Let ${\cal P}$ be a
partition with finite entropy $H_{\mu}({\cal P})$. Then,
$$
\lim_{n\to\infty} \frac 1n\log{\frac{1}{\mu(P(n,x))}}=h_{\mu}(T,{\cal P}) \,,
$$
for $\mu$-almost every $x\in X$.}

\

We will need later the following consequence of Theorem D.

\begin{lemma} \label{egorof}
Let $(X,{\cal A},\mu)$ be a probability space and let  $T:X\longrightarrow X$  be a measure
preserving ergodic transformation and ${\cal P}$ be a partition with finite entropy $H_{\mu}({\cal P})$. Then,
given $\ep>0$ there exists a decreasing sequence of sets $\{E_N^{\ep}\}_{N\in\N}$ such that
\begin{equation} \label{e_n}
\mu(E_N^{\ep})\to 0 \quad  \mbox{  as   }   \quad N\to\infty
\end{equation}
and for all $x\in X\setminus E_N^{\ep}$
\begin{equation} \label{SM}
e^{-j(h_\mu+\ep)}<\mu(P(j,x))<e^{-j(h_\mu-\ep)}\,, \qquad \mbox{
for all } j\geq N\,.
\end{equation}
with $h_\mu=h_{\mu}(T,{\cal P})$ the entropy of $T$ with respect to
$\mu$ and the partition $\cal P$.
\end{lemma}

\begin{proof}
Given $\ep>0$ we define  for all  $j\in\N$  the sets
$$
F_{j}^{\ep}=\left \{x\in X \, : \,\left  |
\frac{1}{j}\log\frac{1}{\mu(P(j,x))}-h_\mu\right |<\ep\right \}\,.
$$
By Theorem D we know that for almost  every  $x\in X$
$$
\lim_{n\to\infty}\frac{1}{n}\log\frac{1}{\mu(P(n,x))} =h_\mu\,.
$$
Therefore there is a set $S$ with $\mu(S)=0$ such that for all
$x\in X\setminus S$ there exists $n(x)\in\N$ such that
$$
 \left | \frac{1}{j}\log\frac{1}{\mu(P(j,x))}-h_\mu\right |<\ep \qquad \mbox{ for all } j\geq
 n(x)\,.
$$
Hence,
$$
X\setminus S\subset \bigcup_{N\in\N}\bigcap_{j\geq N}F_{j}^{\ep}
$$
or equivalently
\begin{equation} \label{med0}
\bigcap_{N\in\N}\bigcup_{j\geq N}(X\setminus F_{j}^{\ep})\subset
S\,.
\end{equation}
We define
$$
E_N^{\ep}=\bigcup_{j\geq N}(X\setminus F_j^{\ep})\,.
$$
By definition $E_{N+1}^{\ep}\subset E_N^{\ep}$ for all $N\in\N$
and by (\ref{med0}) $\mu(\bigcap_N E_N^{\ep})=0$, therefore
$$
\mu(E_N^{\ep})\to 0 \qquad \mbox{ when }\qquad  N\to\infty\,.
$$
Moreover, if $x\in X\setminus E_N^{\ep}$, then $x\in F_j^{\ep}$
for all $j\geq N$, and therefore
$$
e^{-j(h_\mu+\ep)}<\mu(P(j,x))<e^{-j(h_\mu-\ep)} \qquad \mbox{ for
all }\qquad j\geq N\,.
$$
\end{proof}

\begin{proposition} \label{SMsinmalos}
Let $(X,{\cal A},\mu)$ be a probability space, let
$T:X\longrightarrow X$  be a measure preserving mixing
transformation, and ${\cal P}$ be a partition with finite entropy
$H_{\mu}({\cal P})$. Let us denote
$$
X_0 = \bigcap_{n=0}^\infty \bigcup_{P\in \vee_{j=0}^n T^{-j}({\cal
P})} P \,.
$$
Let $P_1,P_2$ be two fixed elements of ${\cal P}$. For $\ep>0$ let
$\{E_M^\ep\}$ be the decreasing sequence of sets given by Lemma
{\rm \ref{egorof}} . If ${\cal S}_{N,M}$ denotes the collection of
the sets $P(N,x)$ verifying
$$
x\in X_0\setminus  E_M^\ep \,, \qquad P(N,x) \subset P_1\,, \qquad
T^N(P(N,x))=P(0,T^N(x))=P_2 \,,
$$
then, for all $M$ and $N$ large enough depending on $P_1$ and $P_2$,
\begin{equation} \label{Sn}
\mu({\cal S}_{N,M}) := \mu \big(\bigcup_{S\in{{\cal S}_{N,M}}} S
\big) \ge \frac 12 \,  \mu (P_1) \, \mu  (P_2) \,.
\end{equation}
\end{proposition}

\begin{proof} We have that
\begin{align}
\mu (P_1) & = \mu ({\cal S}_{N,M}) + \sum_{P\in{\cal P}\setminus
\{P_2\}}
\sum_{\substack{ P(N,x)\subset P_1 \ s.t. \ x\in X_0\setminus E_M^{\ep} \\ T^N(P(N,x)) = P}} \mu(P(N,x)) + \mu (P_1 \cap E_M^\ep) \notag \\
& \le \mu({\cal S}_{N,M}) + \sum_{P\in{\cal P}\setminus \{P_2\}} \mu (P_1 \cap T^{-N} (P)) +  \mu (P_1 \cap E_M^\ep) \notag \\
& = \mu({\cal S}_{N,M}) + \mu (P_1) - \mu (P_1 \cap T^{-N} (P_2))
+  \mu (P_1 \cap E_M^\ep) \notag \,.
\end{align}
Notice that $\lim_{M\to\infty} \mu(E_M^\ep) = 0$  by Lemma
\ref{egorof}  and
$$
\lim_{N\to\infty} \mu (P_1 \cap T^{-N} (P_2)) = \mu(P_1) \,
\mu(P_2)
$$
because $T$ is mixing. Hence, for $M$ and $N$ large enough,
$$
\mu ({\cal S}_{N,M}) \ge \frac 12 \,  \mu (P_1) \, \mu  (P_2)\,.
$$
\end{proof}

\section{Expanding maps.} \label{expandingmaps}

We will say that $(X,d,{\cal A},\lambda,T)$ is an {\it expanding
system} if $(X,\cal A,\lambda)$ is a finite measure space, $\lan$
is a non-atomic measure and the support of $\lan$ is equal to $X$,
$(X,d)$ is a locally complete separable metric space, $\cal A$ is
its Borel $\sigma$-algebra and $T:X \longrightarrow X$ is an {\it
expanding map}, i.e. a measurable transformation satisfying the
following properties:

\begin{itemize}
\item[(A)] There exists  a collection of open sets
${\cal P}_0 = \{P_i\}$ of $X$ such that $\displaystyle \sup_{P\in
{\cal P}_0} \diam (P) <\infty$, and
\begin{itemize}
\item[(1)] $\lambda (P_i)>0$,

\item[(2)] $P_i\cap P_j = \emptyset$ if $i\ne j$,

\item[(3)] $\lambda (X \setminus \cup_i P_i) =0$,

\item[(4)] The restriction  $T\big |_{P_i}$ of $T$ to the set
$P_i$  is injective,

\item[(5)] For each $P_i$, if $P_j\cap T(P_i)\ne \emptyset$, then $P_j \subseteq T(P_i)$.
\item[(6)] For each $P_i$, if  $P_j \subseteq T(P_i)$,  then the map $T\big |_{P_i}^{-1}: T(P_i)\cap P_j\longrightarrow P_i $ is open.

\item[(7)] There is a natural number $n_0>0$ such that
$\lambda (T^{-n_0}(P_i)\cap P_j) > 0$, for all $P_i, \, P_j \in
{\cal P}_0$\,.
\end{itemize}

\item[(B)] There exists a measurable map $\J:X \longrightarrow
[0,\infty)$, $\J>0$ in $\cup_{P\in{\cal P}_0} P$,  such that for
all $P_i \in {\cal P}_0$ and for all Borel subset $A$ of $P_i$ we
have that
$$
\lambda (T(A)) = \int_A \J\, d\lambda \,.
$$
and moreover there exist absolute constants $0<\alpha\le 1$ and
$C_1>0$, such that for all $x,y \in P_i$
$$
\left|  \frac {\J(x)}{\J(y)}  -1 \right| \le C_1\,
d(T(x),T(y))^\alpha \,.
$$
\item[(C)] Let us define inductively  the following collections $\{{\cal P}_i\}$ of open  sets:
$$
{\cal P}_1=\bigcup_{P_i\in{\cal P}_0}\{ (T\big|_{P_i})^{-1}(P_j)\,
:\, P_j\in{\cal P}_0 \,, \ P_j \subset T(P_i) \},
$$
and, in general,
$$
{\cal P}_{n}=\bigcup_{P_i\in{\cal
P}_0}\{(T\big|_{P_i})^{-1}(P_j)\, :\, P_j\in{\cal P}_{n-1} \,, \
P_j \subset T(P_i)\}.
$$
Then, there exist absolute constants $\beta >1$ and $C_2>0$ such that for all $x,y$ in the  same element  of ${\cal P}_n$ we  have that
$$
d(T^n(x),T^n(y)) \ge C_2 \beta^n d(x,y)\,.
$$
\end{itemize}

\

\begin{remark}  \label{nota1}
\rm

\begin{itemize}
\item[] \item[1.] It is easy to see that each family ${\cal P}_n$
verifies the properties (A.1), (A.2) and (A.3). Also notice that,
for each $n$, $T({\cal P}_n)$ is equal to  ${\cal P}_{n-1}$ (mod
$0$) in the sense that the image of each element of ${\cal P}_n$
is an element of ${\cal P}_{n-1}$ (mod $0$). \item[2.] From the
properties (A.1), (A.2) and (A.3) it follows that ${\cal P}_0$ is
finite or numerable.
\item[3.] As a consequence of property (C) we have, since
$\sup_{P\in {\cal P}_0}\diam(P)<\infty$, that
$$
\lim_{n\to\infty} \Big( \sup_{P\in{\cal P}_n} \hbox{diam}(P) \Big) =0 \,.
$$
Therefore (see for example [M, p.13]) we deduce that the partition ${\cal P}_0$ is {\it generating}.
This means that ${\cal A}=\bigvee_{n=0}^\infty {\cal P}_n$ (mod 0).
\end{itemize}
\end{remark}

We will define also, for each $n\in \N$,  the function
$$
\J_n(x) = \J(x) \cdot \J(T(x)) \cdots \J(T^{n-1}(x)) \,, \qquad
\hbox{for }  x\in \bigcup_{P\in {\cal P}_{n-1}}P.
$$
Then it follows easily that
\begin{equation} \label{njacobiano}
\int_A f(T^n(x)) \, {\bf J}_n(x)\, d\lambda(x) = \int_{T^n(A)}
f(x)\, d\lambda (x) \,, \qquad \hbox{for all $f\in L^1(\mu)$}\,,
\end{equation}
and, in particular,
$$
\lambda (T^n(A)) = \int_A \J_n \, d\lambda
$$
for each measurable set $A$ contained in some element of ${\cal P}_{n-1}$.

\

Notice that in the definition of an expanding map the measure $\lambda$ it is not required to have special dynamical properties. However it is a remarkable fact that it is possible to find an invariant measure which is essentially comparable to $\lambda$ and has very interesting dynamical properties. More concretely it is known the following result

\

\noindent {\bf Theorem E} ([M, p.172]). {\it  Let  $(X,d,{\cal
A},\lambda,T)$ be an expanding system. Then, there exist a unique
probability measure $\mu$ on $\cal A$ which is absolutely
continuous with respect to $\lambda$ and such that

\begin{itemize}
\item[\rm (i)]  $T$ preserves the measure $\mu$.

\item[\rm (ii)] $d\mu/d\lambda$ is H\"older continuous.

\item[\rm (iii)]  For each $P_i \in {\cal P}_0$ there exist a
positive constant $K_i$ such that
$$
\frac 1{K_i} \le \frac {d\mu}{d\lambda} (x) \le K_i \,, \qquad
\hbox{for all $x\in P_i$}\,.
$$
\item[\rm (iv)] $T$ is exact with respect to $\mu$. \item[\rm (v)]
$\mu(B)= \dfrac 1{\lan (X)} \lim_{n\to\infty} \lambda (T^{-n}(B))$
for every $B\in {\cal A}$.
\end{itemize}
}

In what follows we will refer to $\mu$ as the ACIPM measure
associated to the expanding system.

\begin{remark} \label{fullsets}
\rm Notice that by part (iii) and property (A.3) of expanding maps the measures $\lan$ and $\mu$ have the same zero measure sets and therefore the same full measure sets.
\end{remark}

\begin{remark} \label{comparableabsoluta}
\rm We recall that the condition (A.4) in the definition of
expanding maps says that $T\big| P_i$ must be injective for all
$P_i \in {\cal P}_0$. If we strengthen this condition by
requiring also that
$$
\inf_{P\in {\cal P}_0} \lan (T(P)) >0  \qquad \text{and} \qquad
\sup_{P\in {\cal P}_0} \diam (T(P)) <\infty \,,
$$
or, in particular, if $T:P \longrightarrow X$ is bijective (mod
$0$) for all $P\in {\cal P}_0$ and $X$ is bounded, then, a slight
modification of the proof of Theorem E in \cite{M} (using Remark
\ref{lema3} instead of [M, Lemma 1.5]),  allows to obtain the
property (iii) of $\mu$ with an absolute constant $K$. Therefore
with this additional assumption one have that
$$
\frac 1K \, \lambda (A) \le \mu(A) \le K \, \lambda(A) \,, \qquad
\hbox{for all $A\in {\cal A}$}\,.
$$
Of course, this condition also holds if the partition ${\cal P}_0$
is finite.
\end{remark}

\

We recall that since $\sup_{P\in{\cal P}_0} \diam(P)<\infty$ then the partition ${\cal P}_0$ is  generating, Therefore, for expanding systems $h_\mu(T)=h_\mu(T,{\cal P}_0)$.

\begin{remark} \label{entropia}
\rm By definition of entropy, if $H_\mu({\cal P}_0)$ is finite, then $h_\mu(T)\le H_\mu({\cal P}_0)<\infty$. Also, since $T$ is exact with respect to $\mu$ we have that an expanding
system is Kolmogorov (\cite{M}, p. 158). Also, $X$ is a Lebesgue space (\cite{M}, p. 81). As a
consequence it follows that $h_\mu(T)>0$, (\cite{M}., p. 225).
\end{remark}

For expanding maps there exists an alternative way of computing
the entropy of $T$:

\

\noindent{\bf Theorem F}([M, p. 227]) {\it  Let  $(X,d,{\cal
A},\lambda,T)$ be an expanding system and let $\mu$  be the ACIPM
measure associated to the system. If the entropy $H_{\mu}({\cal
P}_0)$ of the partition ${\cal P}_0$ is finite, then $\log \J$ is
integrable and
$$
h_{\mu}(T)=\int_X\log \J \, d\mu \,.
$$
}

\

\subsection{A code for expanding maps} \label{codes}

We will denote by $P(n,x)$ the element of the collection ${\cal
P}_n$ which contains the point $x\in X$.   Observe that for each
$n$,  $P(n,x)$  is well defined for $x$ belonging to
$$
\Upsilon_n:= \cup \{P: P\in {\cal P}_n \}
$$
and $\Upsilon_n$ has full $\lambda$-measure for property (A.3) for ${\cal P}_n$, see Remark \ref{nota1}.1. Therefore if $x$ belongs to the set
\begin{equation}  \label{defdeX0}
X_0:=\bigcap_{n=0}^\infty \Upsilon_n = \bigcap_{n=0}^\infty \bigcup_{P\in {\cal P}_n} P
\end{equation}
then $P(n,x)$ is
well defined for all $n$. Moreover, if $x\in \Upsilon_n$
then from the definition of ${\cal P}_n$ we have that $T(x)\in \Upsilon_{n-1}$. Hence,
if $x\in X_0$ we have that $T^\ell (x) \in X_0$ for all $\ell\in\N$, and so
$P(n,T^\ell(x))$ is well defined for all $n, \ell \in\N$.
This set has full $\lambda$-measure since $X\setminus X_0
\subseteq \cup_{n\ge 0} X\setminus \Upsilon_n$ and
this set has zero $\lambda$-measure by (A.3) for ${\cal
P}_n$. Hence, for almost every $x\in X$, $P(n,T^\ell(x))$ is
defined for all $n, \ell\in\N$.
An easy consequence of the definition of $P(n,x)$ that we will use in the sequel  is that
\begin{equation} \label{TdeIn}
T(P(n,x)) = P(n-1, T(x)) \,,  \qquad n\ge 1 \,.
\end{equation}

If $x\in X_0$  then, since $P(n+1,x)\subset P(n,x)$ and
$\diam(P(n,x))\to 0$ when $n\to\infty$, we have that
$$
\bigcap_n P(n,x)=\{ x\}
$$
and so the sequence $\{P(n,x)\}_n$ determines to the point $x$. Moreover, from (\ref{TdeIn}) we have that $T^n(P(n,x))=P(0,T^n(x))$, and it is not difficult to see that
\begin{align}
P(k,x) & =T\big |_{P(0,x)}^{-1}T\big |_{P(0,T(x))}^{-1}\ldots T\big |_{P(0,T^{k-1}(x))}^{-1}\Big ( P(0,T^k(x))\Big )\notag \\
& =T\big |_{P(0,x)}^{-1}T\big |_{P(0,T(x))}^{-1}\ldots T\big |_{P(0,T^{k-2}(x))}^{-1}\Big (T^{-1}\big (P(0,T^k(x))\big )\cap P(0,T^{k-1}(x)) \Big)\notag \\
& =T\big |_{P(0,x)}^{-1}T\big |_{P(0,T(x))}^{-1}\ldots T\big |_{P(0,T^{k-3}(x))}^{-1}\Big (T^{-2}\big (P(0,T^k(x))\big )\cap T^{-1}\big ( P(0,T^{k-1}(x))\big )\cap P(0,T^{k-2}(x))\Big )\notag \\
&=\ldots =\bigcap_{n=0}^k T^{-n}(P(0,T^n(x))) . \notag
\end{align}
Hence
$$
\bigcap_{n=0}^\infty T^{-n}(P(0,T^n(x)))=\bigcap_{n=0}^\infty
P(n,x)=\{ x\}
$$
and the sequence $\{P(0,T^n(x))\}_n$ also determines the point $x$.

We will also define the set $X_0^+$ as the union of $X_0$ with the set of points $x\in X$ verifying that
there exists a sequence $\{P_n\}$, with $P_n \in {\cal P}_n$ and $P_{n+1} \subset P_n$, such that
$$
\bigcap_{n=0}^\infty \text{closure}(P_n)=\{ x\} \,.
$$
We remark that for points $x\in X_0^+\setminus X_0$ the sequence $\{P_n\}$ is not uniquelly determinated by $x$. From now on, for each $x\in X_0^+\setminus X_0$ we make an election of $\{P_n\}$ and we denote $P_n$ by $P(n,x)$. Also by $P(0,T^n(x))$ we mean $T^n(P(n,x))$.
We are extending in this way the definition of $P(n,x)$ and $P(0,T^n(x))$ given for points in $X_0$ in such a way that for points in $ X_0^+\setminus X_0$ we also have that $T^n(P(n,x))=P(0,T^n(x))$.

\begin{definition}\label{code}
If  $x\in X_0^+$, then  we will code $x$ as the sequence $\{ i_0,
i_1, \ldots\}$ and we will write $x=[\;i_0 \;i_1\:\ldots\;]$ if
and only if
$$
P(0,T^n(x))=P_{i_n}\in{\cal P}_0,\qquad \text{  for all } \quad  n=0,1,2,\ldots
$$
\end{definition}
\begin{remark} If $x=[\;i_0 \; i_1 \; i_2 \ldots\;]$ then $T(x)=[\;i_1\:i_2\;i_3\;\ldots\;]$. Therefore $T$ acts as the left shift on the space of all codes.
\end{remark}

\

\subsection{Some properties of expanding maps}\label{PropEM}

Let $(X,d,{\cal A},\lambda,T)$ be an expanding system. Following
\cite{M} we have

\begin{proposition}  \label{cotajacobiano}
There exists an absolute constant $C>0$ such that for all $x_0 \in
X_0^+$ and for all natural number $n$ we have that if $x,y \in
P(n,x_0)$ then
\begin{equation} \label{lachachi}
\frac {\J_{s}(x)}{\J_{s}(y)} \le C \,, \qquad \hbox{for
$s=1,\dots,n$}\,.
\end{equation}
Moreover, if $\sup_{P\in {\cal P}_0} \diam(T(P))<\infty$, then {\rm (\ref{lachachi})} holds for $s=n+1$.
\end{proposition}

\begin{proof} We will prove the lemma for the case $s=n+1$.
If $x,y\in P(n.x_0)$ we have, from properties (B) and (C), that
\begin{align}
\frac{\J_{n+1}(x)}{\J_{n+1}(y)} & = \prod_{k=0}^n \frac
{\J(T^k(x))}{\J(T^k(y))}
\le \prod_{k=0}^n (1+C \, d(T^{k+1}(x),T^{k+1}(y))^\alpha ) \notag \\
& = (1+C \, d(T^{n+1}(x),T^{n+1}(y))^\alpha) \prod_{k=0}^{n-1} (1+C\,  d(T^{k+1}(x),T^{k+1}(y))^\alpha) \notag \\
& \le (1+C \, [\hbox{diam\,}T(P(0,T^n(x_0)))]^\alpha) \prod_{k=0}^{n-1} (1+C\, \beta^{-\alpha(n-(k+1))} d(T^{n}(x),T^{n}(y))^\alpha) \notag
\end{align}
since $x,y\in P(n,x_0)$ implies that $T^{k+1}(x),T^{k+1}(y)\in P(n-(k+1),T^{k+1}(x_0))$ for $k=0,\dots,n-1$. Therefore,
\begin{align}
\frac{\J_{n+1}(x)}{\J_{n+1}(y)} &
\le (1+C \, [\hbox{diam\,}T(P(0,T^n(x_0)))]^\alpha) \prod_{j=0}^{n-1} (1+C\, \beta^{-\alpha j} [\hbox{diam\,} P(0,T^n(x_0))]^\alpha) \notag \\
&\le (1+ C\, \overline{D}^\alpha) \exp \Big[ C\, D^\alpha \sum_{j=0}^{\infty}
\beta^{-\alpha j} \Big] \le C \,, \notag
\end{align}
where $D=\sup_{P\in {\cal P}_0} \diam(P)$ and $\overline{D}=\sup_{P\in {\cal P}_0} \diam(T(P))$.
\end{proof}

An easy consequence of the above bound, that we will often use, is the following one:

\begin{proposition} \label{elchachi}
If $P$ is an element of ${\cal P}_n$, i.e. if $P=P(n,x)$ for some
$x\in X$, and $P'$ is a measurable subset of $P$ then
$$
\frac 1C \frac {\lambda(T^{j}(P'))}{\lambda (T^{j}(P))}\le \frac
{\lambda(P')}{\lambda(P)} \le C \frac {\lambda(T^{j}(P'))}{\lambda
(T^{j}(P))}\, , \qquad \hbox{ for } j=1,\dots, n\,,
$$
with $C$ an absolute constant. Moreover, if $\sup_{P\in {\cal P}_0}\diam(T(P))<\infty$, then the above
inequality is true for $j=n+1$.
\end{proposition}

\begin{proof}
Using (\ref{njacobiano}) we get
$$
 \frac {\inf_{y\in P}{\bf  J_{j}}(y)}{\sup_{x\in P}{\bf J_{j}}(x)} \, \frac {\lambda(P')}{\lambda(P)} \le
 \frac {\lambda(T^{j}(P'))}{\lambda (T^{j}(P))} = \frac {\int_{P'} {\bf J_{j}} \, d\lambda}{\int_{P} {\bf J_{j}} \,
 d\lambda} \le \frac {\sup_{x\in P} {\bf J_{j}}(x)}{\inf_{y\in P}{\bf  J_{j}}(y)} \, \frac {\lambda(P')}{\lambda(P)}
$$
and the result follows from Proposition \ref{cotajacobiano}.
\end{proof}

\begin{lemma} \label{elquefunciona}
If $A\in {\cal A}$ and $Q\in {\cal P}_m$ for some m, then
$$
\lambda (T^{-\ell}(A)\cap Q) = \int_A \sum_{y\in T^{-\ell}(x)\cap
Q} \frac 1{\J_\ell(y)} \, d\lambda (x) \,, \qquad \hbox{for }
\ell=1,2,\dots \,.
$$
\end{lemma}

\begin{proof} We may assume that $T^{-\ell}(A)\cap Q\neq \emptyset$ and $A\subseteq P\in {\cal P}_m$. The general result follows from the fact that ${\cal P}_m$ is a partition of $X$.
Then we have that $T^{-\ell}(A)\cap Q$ is a union
of some elements $B_1,B_2,\ldots$ such that $B_i\subset P_i\in {\cal P}_{\ell+m}$ and
$T^\ell \big|_{B_j}: B_j \longrightarrow A$ is bijective for all
$j$. Let us denote by $S_j$ its inverse map,
$S_j:A\longrightarrow B_j$. Then
$$
\lambda (T^{-\ell}(A)\cap Q) = \sum_j \lambda (B_j) = \sum_j
\lambda (S_j(A)) \,.
$$
But using (\ref{njacobiano}) we deduce that
$$
\lambda (S_j(A)) = \int_{S_j(A)} d\lambda = \int_{S_j(A)} \frac
{\J_\ell (x)}{\J_\ell ((S_j (T^{\ell}(x)))} \, d\lambda(x) =
\int_{T^{\ell}(S_j(A))} \frac 1{\J_\ell(S_j(x))} \, d\lambda
(x)\,.
$$
Therefore, since $T^{\ell}(S_j(A))=A$, we have
$$
\lambda (T^{-\ell}(A)\cap Q) = \sum_j \int_A \frac 1{\J_\ell
(S_j(x))} \, d\lambda (x) = \int_A \sum_j \frac 1{\J_\ell
(S_j(x))} \, d\lambda (x)\,.
$$
If we denote, for each $j$, $y=S_j(x)$ we have that $y\in
T^{-\ell}(x)\cap B_j$ and that $y$ is unique. This observation
completes the proof.
\end{proof}

\begin{lemma} \label{nucleoacotado}
If $x\in P_0\in {\cal P}_0$ and $z\in Q\in {\cal P}_m$, then
$$
\sum_{y\in T^{-\ell}(x)\cap Q} \frac 1{\J_\ell(y)} \le
\begin{cases}
C\, \lan(P(\ell,z))& \text{ if  $\ell< m$},
\\
&
\\
C\, \lan(Q)& \text{ if  $\ell \geq m.$}
\end{cases}
$$
with $C>0$ a constant depending on $P_0$.
\end{lemma}

\begin{proof}
Using (\ref{TdeIn}), (\ref{njacobiano}) and Proposition \ref{cotajacobiano} we deduce that
$$
\lambda (P_0) = \lambda (P(0,x)) =\lambda (P(0,T^{\ell}(y))=
\lambda (T^\ell (P(\ell,y))) = \int_{P(\ell,y)} \J_\ell \,
d\lambda \asymp  \J_\ell (y) \lan(P(\ell,y))  \,.
$$
Therefore
\begin{equation} \label{unoentrejotaele}
\frac 1{\J_\ell (y)} \asymp \frac {\lan(P(\ell,y))}{\lambda (P_0)}
\,.
\end{equation}
If $\ell\geq m$ we have that $P(\ell,y)\subset Q$ for all $y\in T^{-\ell}(x)\cap Q$ and so
$$
\sum_{y\in T^{-\ell}(x)\cap Q} \frac 1{\J_\ell(y)} \asymp \frac
1{\lambda (P_0)} \sum_{y\in T^{-\ell}(x)\cap Q} \lambda
(P(\ell,y)) \le \frac {C}{\lambda (P_0)} \, \lambda (Q)
$$
On the other hand, if $\ell < m$ the map $T^{\ell}$ is injective
in $Q$ and therefore there is at most one point $y\in
T^{-\ell}(x)\cap Q$. Since $P(\ell,y)\supset P(m,y)=Q$ we also
have that $P(\ell,z)=P(\ell,y)$ for any $z\in Q$. Therefore, in
this case, the result follows from (\ref{unoentrejotaele}).
\end{proof}

\begin{remark} \label{lema3}  \rm
Under the same hypotheses for $x$ and $Q$, and if
$$
C_0:=\inf_{P\in{\cal P}_0} \lan(T(P)) >0 \qquad \text{and} \qquad
D_0:=\sup_{P\in {\cal P}_0} \diam (T(P)) <\infty ,
$$
then
$$
\sum_{y\in T^{-\ell}(x)\cap Q} \frac 1{\J_\ell(y)} \le
\begin{cases}
C\, \lan(P(\ell-1,z))& \text{ if  $\ell< m+1$},
\\
&
\\
C\, \lan(Q)& \text{ if  $\ell \geq m+1.$}
\end{cases}
$$
with $C$ a constant depending on $C_0$ and $D_0$.

\begin{proof}
Notice that from Proposition \ref{cotajacobiano} we get that
$$
\lan(T(P(0,T^{\ell-1}(y))))=\lambda (T^\ell (P(\ell-1,y))) = \int_{P(\ell-1,y)} \J_\ell \,
d\lambda \asymp  \J_\ell (y) \lan(P(\ell-1,y))  \,.
$$
Therefore,
$$
\frac 1{\J_\ell (y)} \asymp \frac {\lan(P(\ell-1,y))}{\lan(T(P(0,T^{\ell-1}(y))))} \le \frac 1{C_0}\,\lan(P(\ell-1,y))
\,.
$$
The rest of the proof is similar to the proof of Lemma
\ref{nucleoacotado}.
\end{proof}
\end{remark}

\begin{proposition}\label{yaesta}
Let $\mu$ be the ACIPM measure associated to the expanding system.
Let $A\in {\cal A}$ and $Q\in {\cal P}_m$ with $A, Q\subset P_0
\in {\cal P}_0$. Then, we have that
$$
 \mu(T^{-\ell}(A)\cap Q) \leq \begin{cases}
C\, \mu(A)\mu(P(\ell,z))& \text{ if  $\ell < m$},
\\
&
\\
C\, \mu(A)\mu(Q)& \text{ if  $\ell \ge m.$}
\end{cases}
$$
where $z$ is any point of $Q$ and $C>0$ is a constant depending on $P_0$.
\end{proposition}

\begin{proof}
By Theorem E we know that
\begin{equation}\label{compara}
\mu(V)\asymp\lambda(V) \hbox{ for all measurable set } V\subset P_0
\end{equation}
and  the result is a consequence of lemmas \ref{elquefunciona} and \ref{nucleoacotado}.
\end{proof}

\begin{remark} \rm
If $\inf_{P\in{\cal P}_0} \lan(T(P)) >0$ and $\sup_{P\in {\cal
P}_0} \diam (T(P)) <\infty$, then $\lan$ and $\mu$ are comparable
in the whole $X$ and it is not necessary in the statement of the
Proposition \ref{yaesta} that $A,Q\subset P_0$.
\end{remark}

Recall now the definitions of lower and upper ${\cal
P}$-dimensions, see Definition \ref{dimensiones}. Since the
sequence $\{P(n,x_0)\}_{n\in\N}$  is defined for all $x_0\in
X_0^+$, we have that $\dimsup$ and $\diminf$ are also defined for
$x_0\in X_0^+$,

\

\begin{lemma}\label{compdiam1}
Let $x_0\in X_0^+$ such that $\underline{\delta}_\lambda(x_0)>0$.
Given $0<\ep<\diminf$ there exists $N\in\N$ such that for all
$n\geq N$
$$
\lambda(P(n,x_0))\leq \beta^{-n(\underline{\delta}_\lambda(x_0)-\ep)}
$$
with $\beta>1$ the constant in the property {\rm (C)} of
expanding maps.
\end{lemma}

\begin{proof} By definition of $\underline{\delta}_\lambda(x_0)$ we have that  for $n$ large enough
$$
\lambda(P(n,x_0))\leq  (\hbox{diam}(P(n,x_0)))^{\underline{\delta}_\lambda(x_0)-\ep/2}
$$
Now, if $x_0\in X_0$, from the property (C) of expanding maps we get that
$$
C_2\beta^{n}\hbox{diam}(P(n,x_0))\leq
\hbox{diam}(P(0,T^{n}(x_0)))\le D=\sup_{P\in {\cal P}_0}
\diam(P)<\infty.
$$
The result follows for $x_0\in X_0$ from these two inequalities.
If $x_0\in X_0^+\setminus X_0$, then $P(n,x_0)=P(n,x)$ with $x\in
P(n,x_0)\cap X_0$ and from this fact we conclude that also
$C_2\beta^{n}\hbox{diam}(P(n,x_0))\le \sup_{P\in {\cal P}_0}
\diam(P)$ for these points.
\end{proof}

Another quantity that we will need is the following:

\begin{definition} \label{tau}
The rate of decay of the measure $\lan$ at $x\in X_0$ with respect to the partition ${\cal P}_0$ is defined as
$$
\overline{\tau}_\lan(x)=\limsup_{n\to\infty}\frac 1n \log\frac
{\lan(P(n,x))}{\lan(P(n+1,x))}\,.
$$
\end{definition}

Notice that we can extend the definition of $\overtau$  to
 all $x_0\in X_0^+$.

\begin{lemma}
\label{justificacion} If the entropy $H_\mu({\cal P}_0)$ of the
partition ${\cal P}_0$ with respect to the ACIPM measure
associated to the expanding system is finite, then the set of
points $x_0 \in X_0^+$ verifying
\begin{equation} \label{tau0}
\overtau=\lim_{n\to\infty}\frac 1n \log\frac
{\lan(P(n,x_0))}{\lan(P(n+1,x_0))}=0
\end{equation}
has full $\lan$-measure. Besides {\rm (\ref{tau0})} holds if
$\,\sup_{P\in {\cal P}_0} \diam (T(P))<\infty$ and $x_0 \in X_0^+$
verifies
$$
\inf_j \lan(P(0,T^j(x_0))) >0 \,.
$$
In particular, if the partition ${\cal P}_0$ is finite, then {\rm
(\ref{tau0})} holds for all $x_0 \in X_0^+¼$.
\end{lemma}

\

\begin{proof}.
From part (iii) of Theorem E we know that the measures $\mu$ and
$\lan$ are comparable in each element of the partition ${\cal
P}_0$  and as a consequence the zero measure sets are the same for
$\mu$ and $\lan$. Hence from Theorem D we have that for
$\lan$-almost all $x\in X$
$$
\lim_{n\to\infty}\frac
1n\log\frac{1}{\lan(P(n,x))}=\lim_{n\to\infty}\frac
1n\log\frac{1}{\mu(P(n,x))}=h_{\mu},
$$
and therefore (\ref{tau0}) holds.

Now let us prove that if  $\inf_j \lan(P(0,T^j(x_0)))>0$ and
$\sup_{P\in {\cal P}_0}\diam(T(P))<\infty$, then {\rm
(\ref{tau0})} holds for all $x_0 \in X_0^+$. By Proposition
\ref{elchachi} and (\ref{TdeIn}) we have that for all $x_0\in X_0$
$$
\frac{\lan(P(n,x_0))}{\lan(P(n+1,x_0))}\le C
\frac{\lan(T(P(0,T^n(x_0)))}{\lan(P(0,T^{n+1}(x_0)))}\le C
\frac{\lan(X)}{\inf_j \lan(P(0,T^j(x_0)))}  < C'
$$
with $C'>1$ a constant. This implies (\ref{tau0})  for all $x_0\in
X_0$. If $x_0 \in X_0^+\setminus X_0$, then $P(n+1,x_0)=P(n+1,x)$
for $x\in P(n+1,x_0)\cap X_0$ and, since $P(n+1,x_0)\subset P(n,x_0)$, we also have that $x\in P(n,x_0)\cap X_0$ and therefore also $P(n,x)=P(n,x_0)$. The result for these points follows now from the last chain of inequalities.
\end{proof}

\

\begin{lemma}\label{compdiam2} Let $x_0\in X_0^+$ be a point such that
$\overline{\delta}_\lambda(x_0)<\infty$ and  $\overtau<\infty$.
Given $\ep>0$ there exists $N\in\N$ such that for all $n\geq N$
$$
\lambda(P(n,x_0))\geq (\hbox{\rm diam}(P(n-1,x_0)))^{\overline{\delta}_\lambda(x_0)+\overtau/\log \beta+\ep}
$$
where $\beta>1$ is the constant in the property {\rm (C)} of
expanding maps.

\end{lemma}

\begin{proof}
By  Definition \ref{tau} we have that for $n$ large
enough
$$
\frac{1}{n-1} \log \frac {\lan(P(n-1,x_0))}{\lan(P(n,x_0))} <
\overtau +\frac 13 \ep\log\beta \,.
$$
Hence, for $n$ large enough,
\begin{equation}\label{lcompmed2}
\lambda(P(n,x_0))\geq \beta^{-(n-1)\ep/3} e^{-(n-1)\overtau} \lambda(P(n-1,x_0))\, .
\end{equation}

But from the property (C) of expanding maps, if $x_0\in X_0$ we have that,
$$
C_2\beta^{n-1}\hbox{diam}(P(n-1,x_0))\leq
\hbox{diam}(P(0,T^{n-1}x_0))\le D
$$
with $D=\sup_{P\in {\cal P}_0} \diam(P)$.  If $x_0\in
X_0^+\setminus X_0$, we obtain the same conclusion since
$P(n-1,x_0)=P(n-1,x)$ for $x\in P(n-1,x_0)\cap X_0$. Therefore, in
any case, we get that
\begin{equation} \label{diamn-1}
\beta^{-(n-1)\ep/3} e^{-(n-1)\overtau} \geq C(\hbox{diam}(P(n-1,x_0)))^{\ep/3+\overtau/\log\beta}\
\end{equation}
with $C>0$.

Finally from the definition of $\overline{\delta}_\lambda(x_0)$ we
have that for $n$ large
\begin{equation}\label{deltamed}
\lambda (P(n-1,x_0)) \ge \hbox{diam}
(P(n-1,x_0))^{\dimsup+\ep/3} \,.
\end{equation}
The result follows from (\ref{lcompmed2}), (\ref{diamn-1}), and (\ref{deltamed}).
\end{proof}

\

Using  lemma  \ref{egorof} we can define an important subset of
$X_0$ which also has  full $\lan$-measure.  We will refer to this
set in the rest of the paper. The following lemma summarizes its
properties.

\begin{lemma} \label{defdeX1}
Let $(X,{\cal A},\lan,T)$ be an expanding system such that the
entropy $H_\mu({\cal P}_0)$ of the partition ${\cal P}_0$, with
respect to the unique $T$-invariant probability measure which is
absolutely continuous with respect to $\lan$, is finite. Let $X_1$
denote the subset of $X_0$
$$
X_1 = X_0 \setminus [\cup_{m=1}^\infty \cap_{N} E_{N}^{1/m}]\,.
$$
with $\{E_N^{1/m}\}$ the sets given by Lemma {\rm \ref{egorof} } for $\ep=1/m$. Then
$\lan (X_1)=\lan(X)$ and moreover,  if $x_0\in X_1$ then:
\begin{itemize}
\item[\rm (i)] $P(n,T^{\ell}(x_0))$ is well defined for all $n,\ell \in \N$.
\item[\rm (ii)] For all positive integer $m$ there exists $N\in\N$ such that for all $n\geq N$
\begin{equation}\label{pronto}
\frac{1}{M}\,e^{-n(h_{\mu}+1/m)} < \lan(P(n,x_0)) < M\,
e^{-n(h_{\mu}-1/m)}\,
\end{equation}
with $M>0$ depending on $P(0,x_0)$. \item [\rm(iii)]
$$
\overtau=\lim_{n\to\infty}\frac 1n \log\frac
{\lan(P(n,x_0))}{\lan(P(n+1,x_0))}=0\,.
$$
\item[\rm (iv)] $ \overline{\delta}_\lambda(x_0)\le
h_{\mu}/\log\beta<\infty$, with $\beta$ the constant given by
property {\rm (C)} of expanding maps.
 \end{itemize}
\end{lemma}

\begin{proof}
In the proof of Lemma \ref{egorof} we saw that $\mu(\cap_{N}
E_{N}^{1/m})=0$.  Then, by Theorem E, we obtain that
$\lan(\cap_{N} E_{N}^{1/m})=0$, and therefore
$\lan[\cup_{m=1}^\infty \cap_{N} E_{N}^{1/m}]=0$. Hence we have
that $\lan(X_1)=\lan(X_0)$, but when we  defined $X_0$ (see
(\ref{defdeX0})) we showed that $\lan(X_0)=\lan(X)$.

The property  (i)  is satisfied for all points in $X_0$ and
therefore also in $X_1$. If the point $x_0\in X_1$, then, for all
positive integer $m$, $x_0$ does not belong to $\cap_{N}
E_{N}^{1/m}$, and so from Lemma \ref{egorof} we have that for all
$n$ large enough
$$
e^{-n(h_{\mu}+1/m)} < \mu(P(n,x_0)) < e^{-n(h_{\mu}-1/m)}\,.
$$
From part (iii) of Theorem E, we conclude that (\ref{pronto})
holds.

Moreover, from (\ref{pronto}) we also get that
$$
\frac{1}{n} \log \frac {\lan(P(n,x_0))}{\lan(P(n+1,x_0))}
<\frac{2\log M+h_{\mu}+1/m}{n}+\frac 2m   \longrightarrow \frac 2m   \qquad \hbox{ as }\quad n\to\infty
$$
Therefore by taking $m\to\infty$, we get the property (iii).

By property (C) of expanding maps we also have that
$$
\hbox{diam\,}(P(n,x_0)) \le C\, \beta^{-n} \,,
$$
and therefore, using again (\ref{pronto}), we get that
$$
\frac {\log\lan (P(n,x_0))}{\log\hbox{diam\,}(P(n,x_0))} \le \frac
{n(h_\mu+1/m) +\log M}{n\log\beta -\log C} \longrightarrow \frac
{h_\mu+1/m}{\log\beta}
$$
as $n\to\infty$, and so by letting $m\to\infty$ we obtain that
$\dimsup \leq h_\mu/\log\beta<\infty$.
 \end{proof}

\section{Measure results} \label{medida}

\

We want to study the size of the set
$$
{\cal W}(x_0,\{r_n\}) = \{x\in X: \ d(T^n (x),x_0)< r_n \
\mbox{for infinitely many} \ n\}
$$
where $\{r_n\}$ is a given sequence of positive numbers and $x_0$
is an arbitrary point in $X$. Observe that if the sequence
$\{r_n\}$ is constant this set is $T$-invariant, but, in general, this
is not the case.

We are also interested in the size of another set that we will see that is closely related with
${\cal W}(x_0,\{r_n\})$. This set is
$$
\widetilde{\cal W}(x_0,\{t_n\}) = \{x\in X: \ T^k(x) \in
P(t_k,x_0) \ \hbox{for infinitely many $k$}\}
$$
with $\{t_k\}$ an increasing sequence of positive integers and
$x_0\in X_0^+$.


\

Notice also that  if $x\in\widetilde{\cal W}(x_0,\{t_n\})$  then
$P(m,x)$ is well defined for infinitely many $m$. and so it is
well defined for all $m$. Therefore $\widetilde{\cal
W}(x_0,\{t_n\}) \subset X_0$.

\

Let us denote
$$
A_k=T^{-k}(B(x_0,r_k)) \qquad \hbox{and} \qquad  \widetilde{A}_k = T^{-k}(P(t_k,x_0))
$$
With these notations, we have that
\begin{align}
{\cal W}(x_0,\{r_n\}) = \{x\in X: \ x\in A_n \ \mbox{for
infinitely many $n$}\} = \bigcap_{k=1}^\infty \bigcup_{n=k}^\infty
A_n \,, \notag \\
\widetilde{\cal W}(x_0,\{t_n\}) = \{x\in X: \ x\in \widetilde{A}_n
\ \mbox{for infinitely many $n$}\} = \bigcap_{k=1}^\infty
\bigcup_{n=k}^\infty \widetilde{A}_n \,. \notag
\end{align}

\

The following result on the size of these sets is a consequence of
the direct part of Borel-Cantelli lemma and Theorem E.

\medskip

\begin{proposition} \label{serieconverge}
Let $(X,d,{\cal A},\lambda,T)$ be an expanding system.

\begin{itemize}
\item[{\rm i)}] Let $x_0\in \cup_{P\in {\cal P}_0} P$ and let $\{r_n\}$ be a
sequence of positive numbers.
$$
\hbox{If} \qquad \sum_{n=1}^\infty \lan(B(x_0, r_n))<\infty \qquad
\hbox{then} \qquad \lan ({\cal W}(x_0,\{r_n\}))=0.
$$
\item[{\rm ii)}] Let $x_0\in X_0^+$ and let $\{t_n\}$ be a non decreasing sequence of
positive integers.
$$
\hbox{ If \qquad } \sum_{n=1}^\infty \lan(P(t_n,x_0))<\infty ,
\qquad \hbox{ then } \qquad \lan (\widetilde {\cal
W}(x_0\{t_n\}))=0.
$$
\end{itemize}
\end{proposition}

\begin{proof} i) First of all, let us observe that $\lim_{n\to\infty} r_n=0$.
Let $\mu$ be the ACIPM associated to the system. We have that $\mu(B(x_0,r_k))=\mu(A_k)$ for all $k\in\N$ and
therefore
$$
\sum_{n=1}^\infty \mu(A_n) = \sum_{n=1}^\infty \mu(B(x_0,r_n)) < \infty\,,
$$
since for $r_n$ small enough, $B(x_0,r_n) \subset P(0,x_0)$ and
$\lan$ and $\mu$ are comparable in that set by Theorem E. From
Borel-Cantelli lemma it follows that $\mu({\cal
W}(x_0,\{r_n\}))=0$ and using the Remark \ref{fullsets}, we
conclude that $\lan ({\cal W}(x_0,\{r_n\}))=0$. The same argument
works for part ii).
\end{proof}

\begin{corollary}   \label{hartos}
Let $(X,d,{\cal A},\lambda,T)$ be an expanding system. Let $x_0\in
\cup_{P\in {\cal P}_0} P$ and let $\{r_n\}$ be a sequence of
positive numbers. If $\ \sum_{n=1}^\infty \lan (B(x_0,
r_n))<\infty$ then
$$
\liminf_{n\to \infty} \frac {d(T^n(x),x_0)}{r_n} \ge 1 \,, \qquad \hbox{for $\lan$-almost every $x\in X$}\,.
$$
\end{corollary}

\begin{corollary}  \label{hartosinf}
Let $(X,d,{\cal A},\lambda,T)$ be an expanding system. Let $x_0\in
\cup_{P\in {\cal P}_0} P$ such that
$$
0<\Diminf := \liminf_{r \to 0} \frac{\log \lan(B(x_0,r))}{\log r}<\infty\,,
$$
and let
$\{r_n\}$ be a sequence of positive numbers such that $\
\sum_{n=1}^\infty r_n^{\Diminf-\ep}<\infty$ for some
$0<\ep<\Diminf$. Then
$$
\liminf_{n\to \infty} \frac {d(T^n(x),x_0)}{r_n} = \infty \,, \qquad \hbox{for $\lan$-almost every $x\in X$}\,.
$$
If there exists a constant $\Delta(x_0)$ such that $\lan(B(x_0,r))\le Cr^{\Delta(x_0)}$ for all $r$ small enough, then the conclusion holds when $\
\sum_{n=1}^\infty r_n^{\Delta(x_0)}<\infty$.
\end{corollary}

\begin{proof}
By definition of $\Diminf$, we have that for any $r$ small enough
$$
\lan(B(x_0,r)) \le r^{\Diminf-\ep} \,.
$$
Now, for any $m\in \N$, since $\lim_{n\to\infty} r_n=0$, we have that for $n$ big enough, (depending on $m$),
$$
\lan(B(x_0,m r_n)) \le (m r_n)^{\Diminf-\ep} \,.
$$
Therefore, for all $m\in\N$,
$$
\sum_n \lan(B(x_0,m r_n)) \le C_m \, m^{\Diminf-\ep}  \sum_n
r_n^{\Diminf-\ep} <\infty \,.
$$
From Corollary \ref{hartos} we get that, for all $m\in\N$,
$$
\liminf_{n\to \infty} \frac {d(T^n(x),x_0)}{r_n} \ge m \,, \qquad \hbox{for $\lan$-almost every $x\in X$}\,.
$$
The result follows now from the fact that
$$
\left\{ x\in X : \  \liminf_{n\to\infty} \frac {d(T^n(x),x_0)}{r_n} = \infty \right\}
= \bigcap_{m=1}^\infty \left\{ x\in X : \  \liminf_{n\to\infty} \frac {d(T^n(x),x_0)}{r_n} \ge m  \right\}\,.
$$
\end{proof}

\medskip

\begin{remark} \label{ensoporte}
\rm If the measures $\lan$ and $\mu$ are comparable in $X$, then part i) of
Proposition \ref{serieconverge} and its corollaries hold for all
$x_0 \in X$. In particular, this happens if
$$
\inf_{P\in {\cal P}_0} \lan (T(P)) >0  \qquad \text{and} \qquad
\sup_{P\in {\cal P}_0} \diam (T(P)) <\infty \,,
$$
see Remark \ref{comparableabsoluta}. Also, part i) of Proposition  \ref{serieconverge} and its
corollaries hold for those $x_0$ such that the set of elements $P\in {\cal P}_0$ such
that $x_0$ belongs to $\p P$ is finite. For example,
if $X\subseteq \R$ or if the partition ${\cal P}_0$ is finite, then all $x_0\in X$ satisfy the above condition.
\end{remark}

\medskip

\begin{theorem}  \label{elbuenocode}
Let $(X,d,{\cal A},\lambda,T)$ be an expanding system. Let
$x_0$ be a  point of $X_0^+$ such that
$\underline{\delta}_\lambda(x_0)>0$ and let  $\{t_n\}$ be a non
decreasing sequence of positive integers numbers.
$$
\hbox{If} \qquad \sum_{n=1}^\infty \lan(P(t_n,x_0)) = \infty, \qquad \hbox{then}
\qquad \lan(\widetilde{{\cal W}}(x_0,\{t_n\}))=\lan (X).
$$
Moreover, if the partition ${\cal P}_0$ is finite or if the system
has the Bernoulli property, $($i.e. if\, $T(P)=X$ $($mod $0)$ for
all $P\in {\cal P}_0)$, then we have the following quantitative
version:
$$
\lim_{n\to\infty} \frac {\#\{i\le n: \ T^i(x)\in
P(t_i,x_0)\}}{\sum_{j=1}^n \mu(P(t_j,x_0))} =1\,,  \qquad
\hbox{for $\lan$-almost every $x\in X$.}
$$
\end{theorem}

\medskip

In the proof of Theorem \ref{elbuenocode} we will use the following classical result.

\medskip

\noindent {\bf Lemma (Payley-Zygmund Inequality).} {\it Let
$(X,{\cal A}, \mu)$ be a probability space and let
$Z:X\longrightarrow \bf R$ be a positive random variable. Then, for
$0<\tau<1$,
$$
\mu[Z>\tau \,E(Z)] \ge (1-\tau)^2 \frac {E(Z)^2}{E(Z^2)} \,,
$$
where $E(\cdot)$ denotes expectation value.
}

\

\begin{proof}[Proof of Theorem $\ref{elbuenocode}$]  Let $\mu$ be the ACIPM associated to the system. For $j\ge k$, we have that
$$
\mu(\widetilde A_k \cap \widetilde A_j) = \mu
(T^{-k}[P(t_k,x_0))\cap T^{-(j-k)}(P(t_j,x_0))]) = \mu
(P(t_k,x_0))\cap T^{-(j-k)}(P(t_j,x_0)))\,,
$$
and by using Proposition \ref{yaesta} with $\ell=j-k$, $n=t_j$ and $m=t_k$,
and using again that $T$ preserves the measure $\mu$,
we conclude that
\begin{equation}\label{akaj}
 \mu(\widetilde A_k\cap \widetilde A_j) \leq \begin{cases}
C \mu(\widetilde A_j)\mu(P(j-k,x_0))& \text{ if  $j-k< t_k$},
\\
&
\\
C \mu(\widetilde A_j)\mu(\widetilde A_k)& \text{ if  $j-k \ge t_k.$}
\end{cases}
\end{equation}
with $C>0$ depending on $P(0,x_0)$. Let us denote by $Z_n$ and $Z$ the counting functions
$$
Z_n=\sum_{k=1}^n \cark \qquad \hbox{and} \qquad  Z=\sum_{k=1}^\infty \cark \,,
$$
where $\cark$ is the characteristic function of $\widetilde{A}_k$. Observe
that ${\widetilde{ \cal W}}(x_0,\{r_n\})=\{x\in X_0: \
Z(x)=\infty\}$.

If we compute the expectation value of $Z_n^2$ (with respect to $\mu$), we obtain
$$
E(Z_n^2) = E\Big[\sum_{k=1}^n \cark +  \sum^n_{\substack{k,j=1 \\
k\ne j}} \chi_{{}_{\widetilde A_k\cap \widetilde A_j}}\Big]  =
\sum_{k=1}^n \mu(\widetilde A_k) +  2 \sum^n_{\substack{k,j=1 \\
k<j}} \mu(\widetilde A_k\cap \widetilde A_j)
$$
and using (\ref{akaj}) we get
$$
E(Z_n^2) \le E(Z_n) + 2C \sum^n_{\substack{k,j=1 \\
k<j}} \mu(\widetilde A_j) \, \mu(P(j-k,x_0) + 2C
\sum^n_{\substack{k,j=1
\\ k<j}} \mu(\widetilde A_j) \, \mu(\widetilde A_k)\,.
$$
But $\mu(\widetilde A_j)\le \mu(\widetilde A_k)$ for all $j> k$
because $\{t_n\}$ is non decreasing. Therefore
\begin{equation} \label{Zn2}
E(Z_n^2)  \le E(Z_n) + 2C \sum_{k=1}^n \mu(\widetilde A_k)
\sum_{j=k+1}^n \mu(P(j-k,x_0)) +C\, E(Z_n)^2 \,.
\end{equation}

Since by Theorem E the measures $\lambda$ and $\mu$ are comparable
in $P(0,x_0)$ we get from Lemma  \ref{compdiam1} that
\begin{equation}\label{acotado}
\sum_{s=1}^{\infty}\mu(P(s,x_0)) \leq C\,
\sum_{s=1}^{\infty}\beta^{-s(\underline{\delta}_\lambda(x_0)-\ep)}\leq
C'
\end{equation}
with $C'$ a positive constant. From (\ref{Zn2}), and
(\ref{acotado}) we obtain that
\begin{equation}\label{masacotado}
E(Z_n^2)  \le  \Big(1+ 2CC' \Big) \, E(Z_n) + C\,E(Z_n)^2 \,.
\end{equation}
By applying Paley-Zygmund Lemma we obtain from (\ref{masacotado}) that
\begin{align} \label{casicasi}
\mu [\{x\in X: \ Z(x) > \tau \, E(Z_n)\}] & \ge \mu [\{x\in X: \ Z_n(x) > \tau \, E(Z_n)\}] \notag \\
&\ge (1-\tau)^2  \frac {E(Z_n)}{1+ 2CC'+ C\,E(Z_n)} \,.
\end{align}

Using again that $\lan$ and $\mu$ are comparable in $P(0,x_0)$, we get that
$$
E(Z_n) = \sum_{k=1}^n \mu(\widetilde{A}_k)  = \sum_{k=1}^n \mu(P(t_k,x_0)) \ge C
\sum_{k=1}^n \lan(P(t_k,x_0))
$$
and from the hypothesis of the theorem, we obtain that $E(Z_n) \to \infty$ as $n\to\infty$.
Hence, we have from (\ref{casicasi}) that
$$
\mu [\{x\in X: \ Z(x) = \infty\}] \ge \frac 1C (1-\tau)^2 \,,
\qquad \hbox{for $0<\tau<1$}\,.
$$
and we conclude that $\widetilde {\cal W}(x_0,\{t_n\})$
has positive $\mu$-measure. If we denote, for each $n\in\bf N$
$$
\widetilde {\cal W}_n(x_0,\{t_n\}) = \{x\in X: \ T^{k-n}(x) \in
P(t_k,x_0) \ \hbox{for infinitely many $k$ with $k\ge n$}\}\,,
$$
it is easy to see that
$$
\widetilde{\cal W}(x_0,\{t_n\}) = T^{-n}(\widetilde{\cal
W}_n(x_0),\{t_n\}) \qquad \hbox{for each $n\in\bf N$}
$$
and since $T$ is exact with respect to $\mu$ (see Theorem E) it
follows that $\widetilde{\cal W}(x_0,\{t_n\})$ has full
$\mu$-measure. Therefore from Remark \ref{fullsets} we conclude
that $\widetilde{\cal W}(x_0,\{t_n\})$ has full $\lan$-measure.

Finally, if the system has the Bernoulli property then the
correlation coefficients of the sets $\{P(n,x_0)\}_{n\in \N}$ have
exponential decay, see \cite{Y}. Concretely, she proves that
\begin{equation} \label{UM}
|\mu(T^{-\ell}(P(n,x_0))\cap
P(m,x_0))-\mu(P(n,x_0))\mu(P(n,x_0))|\le
C\,\mu(P(n,x_0))\,e^{-\alpha\ell}
\end{equation}
for some absolute positive constants $C$ and $\alpha$ and for all
$m,n,\ell\in\N$. The same argument used in the proof of Theorem 1 in \cite{FMP},
gives the quantitative version.

If the partition ${\cal P}_0$ is finite, then the dynamical system $(X,{\cal A},\mu,T)$ is
isomorphic via coding to a (one-sided) subshift of finite type. The stochastic matrix $M$ of
this subshift is defined in the following way: $p_{i,j}=\mu(T^{-1}(P_j)\cap P_i)/\mu(P_i)$ where $P_i,P_j \in {\cal P}_0$.
Property (A.7) implies that $M$ verifies that $M^{n_0}$ has all its entries positive, see for example \cite{KH}, p.158 or Lemma 12.2 in \cite{M}. This implies that the shift $\sigma$ is mixing, see for example Proposition 12.3 in \cite{M},  and moreover
(\ref{UM}) follows from the Perron-Frobenius theorem (see, for example \cite{KH} or \cite{M}; see also \cite{CK}).
\end{proof}

We state now the following corollary of this proof.

\begin{corollary}
\label{casitodopuntocode} Let $(X,d,{\cal A},\lambda,T)$ be an
expanding system with finite entropy $H_\mu({\cal P}_0)$ with
respect to the partition ${\cal P}_0$ where $\mu$ is the ACIPM
associated to the system. Let $\{t_n\}$ be a non decreasing
sequence of positive integers.

Then for $\lan$-almost all point $x_0\in X_0^+$, more concretely if $x_0\in X_1$ $($see definition in Lemma {\rm \ref{defdeX1})}, we have that
$$
\hbox{If} \qquad \sum_{n=1}^\infty \lan(P(t_n,x_0)) = \infty \qquad \hbox{then} \qquad
\lan(\widetilde{{\cal W}}(x_0,\{t_n\}))=\lan (X).
$$
\end{corollary}

\begin{proof} If $x_0\in X_1$ we have from Lemma \ref{defdeX1} that for all
$m\in\N$ there exists $N\in\N$ such that
$$
e^{-n(h_\mu+1/m)} \le \mu(P(n,x_0)) \le e^{-n(h_\mu-1/m)}
$$
for all $n\ge N$. Therefore we can substitute the inequality
(\ref{acotado}) by
$$
\sum_{s=1}^\infty \mu(P(s,x_0)) \le C\sum_{s=1}^\infty
e^{-s(h_\mu-1/m)} \le C' < \infty\,,
$$
since $h_\mu>0$ (see Remark \ref{entropia}) and we can take $m$
large enough so that $0<1/m<h_\mu$. Hence we do not need now the
hypothesis $\diminf>0$.
\end{proof}

\medskip

\begin{theorem}  \label{elbueno}
Let $(X,d,{\cal A},\lambda,T)$ be an expanding system. Let $x_0$
be a point of $X_0^+$ such that
$$
\overtau<\infty \qquad \text{and} \qquad
0<\underline{\delta}_\lambda(x_0) \le
\overline{\delta}_\lambda(x_0) <\infty \,,
$$
and let  $\{r_n\}$ be a non increasing sequence of positive
numbers.
$$
\hbox{If} \qquad\sum_{n=1}^\infty
r_n^{\overline{\delta}_\lambda(x_0)+\overtau/\log\beta+\ep} = \infty \, \quad
\hbox{for some $\ep>0,\quad $ then} \qquad \lan({\cal W}(x_0,\{r_n\}))=
\lan(X)\,.
$$
Moreover, if the partition ${\cal P}_0$ is finite or if the system
has the Bernoulli property, $($i.e. if\, $T(P)=X$$ ($mod $0)$ for
all $P\in {\cal P}_0)$, then we have the following quantitative
version:
$$
\liminf_{n\to\infty} \frac {\#\{i\le n: \ d(T^i(x),x_0)\le
r_i\}}{\sum_{j=1}^n
r_j^{\overline{\delta}_\lambda(x_0)+\overtau/\log\beta+\ep}} \geq
C \,, \qquad \hbox{for $\lan$-almost every $x\in X$},
$$
with $C$ a positive constant depending on $x_0$ and on the
comparability constants between $\lan$ and $\mu$ at $P(0,x_0)$.
\end{theorem}

\begin{remark}\rm
If the correlation coefficients of the balls $\{B(x_0,r_n)\}$ had
exponential decay, i.e. if they verify the relations
$$
|\mu(T^{-\ell}(B(x_0,r_n))\cap
B(x_0,r_m))-B(x_0,r_n)B(x_0,r_m)|\le
C\,\mu(B(x_0,r_n))\,e^{-\alpha\ell}
$$
for some absolute positive constants $C$ and $\alpha$ and for all
$n,\ell \in\N$, then using the same arguments that in Theorem 1
in \cite{FMP} we would have
$$
\liminf_{n\to\infty} \frac {\#\{i\le n: \ d(T^i(x),x_0)\le r_i\}}{\sum_{j=1}^n \mu(B(x_0,r_j)} =1 \,, \qquad
\hbox{for $\lan$-almost every $x\in X$},
$$
\end{remark}

\begin{remark} \label{tau=0}
\rm We recall that by Lemma {\rm \ref{justificacion}} we
know that $\overtau=0$ for $\lan$-almost all $x_0\in X$. We have
also that $\overtau=0$ if $\inf_j \lan(P(0,T^j(x_0)))>0$ and
$\sup_{P\in {\cal P}_0}\diam(T(P))<\infty$. In particular, if the
partition ${\cal P}_0$ is finite and $\sup_{P\in {\cal
P}_0}\diam(T(P))<\infty$, then $\overtau=0$ for all $x_0\in
X_0^+$.
\end{remark}

\begin{corollary} \label{hartos2}
Under the same hypothesis than Theorem {\rm {\ref{elbueno}}}, if
\begin{equation} \label{seriedivergeep}
\sum_{n=1}^\infty r_n^{\dimsup+\overtau/\log\beta+\ep}=\infty, \qquad \hbox{  for
some } \quad \ep>0\,,
\end{equation}
then,
$$
\liminf_{n\to\infty} \frac {d(T^n(x),x_0)}{r_n} = 0\,, \qquad
\text{for $\lan$-almost all $x\in X$.}
$$
\end{corollary}

\begin{proof}
From Theorem \ref{elbueno} it follows easily that if the radii $r_n$ verify (\ref{seriedivergeep}) then
$$
\liminf_{n\to\infty} \frac {d(T^n(x),x_0)}{r_n} \le 1\,, \qquad
\text{for $\lan$-almost all $x\in X$.}
$$
But notice that for any $m\in \N$ the radii $r_n/m$ also verify  (\ref{seriedivergeep}) and therefore we get that
for any $m\in \N$
$$
\liminf_{n\to\infty} \frac {d(T^n(x),x_0)}{r_n} \le \frac 1m \,, \qquad
\hbox{for $\lan$-almost all $x\in X$.}
$$
The result follows now from the fact that
$$
\left\{ x\in X : \  \liminf_{n\to\infty} \frac {d(T^n(x),x_0)}{r_n} = 0 \right\}
= \bigcap_{m=1}^\infty \left\{ x\in X : \  \liminf_{n\to\infty} \frac {d(T^n(x),x_0)}{r_n} \le \frac 1m  \right\}\,.
$$
\end{proof}

\

\begin{proof}[Proof of Theorem $\ref{elbueno}$]  Let $\mu$ be the ACIPM associated to the system.
Given $x_0\in X_0^+$ and the sequence
$r_k$ we define $t_k$ as the smallest integer so that
\begin{equation} \label{defdetk}
P(t_k,x_0)\subset  B(x_0, r_k).
\end{equation}
Hence, $\widetilde{\cal W}(x_0,\{t_n\})\subset {\cal W}(x_0,\{r_n\})$.

Moreover, since $\dimsup<\infty$ and $\overtau<\infty$, from Lemma \ref{compdiam2}, we get that
\begin{equation} \label{sonlas8}
\sum_{k=1}^n \lan(P(t_k,x_0))\geq \sum_{k=1}^n
(\hbox{diam}(P(t_k-1,x_0)))^{\overline{\delta}_\lambda(x_0)+\overtau/\log\beta+\ep}
\end{equation}
But  from the definition of $t_k$ we have that
$$
P(t_k-1,x_0)\not\subset B(x_0,r_k)
$$
and since $x_0\in P(t_k-1,x_0)$ we can conclude that
$$
 \hbox{diam}(P(t_k-1,x_0))\geq r_k \,.
$$
Therefore we get that
$$
\sum_{k=1}^\infty \lan(P(t_k,x_0)) \geq \sum_{k=1}^{\infty}
r_k^{\overline{\delta}_\lambda(x_0)+\overtau/\log\beta+\ep} =\infty
$$
and from Theorem \ref{elbuenocode} we conclude that $\lan({\cal W}(x_0,\{r_n\}))=\lan(X)$.

Now from (\ref{defdetk}), (\ref{sonlas8}) and the fact that $\lan$ and $\mu$ are comparable on $P(0,x_0)$ we have that
$$
\frac {\#\{i\le n: \ d(T^i(x),x_0)\le r_i\}}{\sum_{j=1}^n r_j^{\overline{\delta}_\lambda(x_0)+\overtau/\log\beta+\ep}} \ge C\, \frac {\#\{i\le n: \ T^i(x)\in P(t_i,x_0)\}}{\sum_{j=1}^n \mu(P(t_j,x_0))}\,.
$$
Hence, the quantitative version follows from Theorem \ref{elbuenocode}.
\end{proof}

We have also the following corollary of the proof of Theorem \ref{elbueno}.

\

\begin{corollary}
\label{casitodopunto} Let $(X,d,{\cal A},\lambda,T)$ be an
expanding system with finite entropy $H_\mu({\cal P}_0)$ with
respect to the partition ${\cal P}_0$ where $\mu$ is the ACIPM
associated to the system. Let $\{r_n\}$ be a non increasing
sequence of positive numbers. Then for $\lan$-almost all point
$x_0\in X$, more concretely if $x_0\in X_1$ $($see definition in
Lemma {\rm \ref{defdeX1})}, we have that
$$
\hbox{if} \qquad\sum_{n=1}^\infty
r_n^{\overline{\delta}_\lambda(x_0)+\ep} = \infty \, \quad
\hbox{for some $\ep>0,\quad $ then} \qquad \lan({\cal
W}(x_0,\{r_n\}))= \lan(X)\,.
$$
In particular, we conclude that, for $\lan$-almost all point
$x_0\in X$,
$$
\liminf_{n\to\infty} \frac {d(T^n(x),x_0)}{r_n} = 0 \,, \qquad
\hbox{for $\lan$-almost all $x\in X$}\,.
$$
\end{corollary}

\begin{proof}  As in the proof of Theorem \ref{elbueno} we define $t_k$ by (\ref{defdetk}) and as a consequence $\widetilde{\cal W}(x_0,\{t_n\})\subset {\cal W}(x_0,\{r_n\})$.
Since $x_0\in X_1$ we have, from Lemma \ref{defdeX1}, that $\overtau=0$ and
$\dimsup<\infty$ and therefore we get the inequality (\ref{sonlas8}) with $\overtau=0$.

From Lemma \ref{defdeX1} we also have that for all $m\in\N$ there exists $N\in\N$ such that
$$
e^{-n(h_\mu+1/m)} \le \mu(P(n,x_0)) \le e^{-n(h_\mu-1/m)}
$$
for all $n\ge N$. The same argument given in Corollary
\ref{casitodopuntocode} allows us to avoid the condition
$\diminf>0$. The fisrt part of the corollary follows from these
facts as in Theorem \ref{elbueno}. The last assertion follows from
Corollary \ref{hartos2}.
\end{proof}

\begin{corollary}
Under the same hypotheses
than Corollary {\rm \ref{casitodopunto}} we have that if
$\Diminf=\dimsup:=D(x_0):=D$ and
$$
\sum_{k=1}^\infty \lan(B(x_0,r_k))^{1+\ep} = \infty \, \quad
\hbox{for some $\ep>0,\quad $ then} \qquad \lan({\cal W}(x_0))=
\lan(X)\,.
$$
\end{corollary}

\begin{proof} From the definition of $\Diminf$ (see Corollary \ref{hartosinf}), the condition $\sum_{k=1}^\infty
\lan(B(x_0,r_k))^{1+\ep} = \infty$ implies that $\sum_{k=1}^\infty
r_k^{(1+\ep)(D-\ep')}=\infty$. But if $\ep'$ is small enough we
have that $(1+\ep)(D-\ep')=D+\ep''$ with $\ep''>0$. Since
$D=\dimsup$ we conclude that $\sum_k r_k^{\dimsup+\ep''}=\infty$.
\end{proof}

\section{Dimension estimates} \label{dimensionresults}

\subsection{Lower bounds for the dimension} \label{lowerdimension}

 Our lower estimate of the dimension is based into the construction of a Cantor like set. Our argument requires to compare  the measures $\lan$ and $\mu$  several times because we use some consequences of the Shannon-McMillan-Breiman Theorem for
the measure $\mu$  (see Section \ref{SMB}) and also some
consequences of the definition of expanding maps involving the
measure $\lan$. We have already mentioned that the measures $\lan$
and $\mu$ are comparable into the blocks of the partition ${\cal
P}_0$, but in order to control the comparability constants  in our
proof, we need the following definition:

\begin{definition}\label{aproximable}
We will say that a point $x_0\in X_0$ is {\underline
{approximable}} if there exist an increasing sequence ${\cal
I}(x_0)=\{p_i\}$ of natural numbers such that
 for all $A\in\cal A$ contained in $P(0,T^{p_i}(x_0))$
for some $i$, we have that
$$
\frac 1K \, \lan(A) \le \mu(A) \le K\,\lan(A) \,.
$$
with  $K>1$ a constant depending on $x_0$.
\end{definition}

\begin{remark} \label{approx}
\rm A mixing version of Poincare's  Recurrence Theorem (see \cite{FMP},
Theorem A') shows that for $\lan$-almost all point $x_0$ there
exists an increasing sequence $\{p_i\}$ such that
$P(0,T^{p_i}(x_0))=P(0,T^{p_1}(x_0))$ for all $i$. Therefore, the
set of approximable points have full $\lan$-measure.
\end{remark}

\begin{remark} \label{approxfinito}
\rm From part (iii) of Theorem E we have that if the partition
${\cal P}_0$ is finite then any point in $X_0^+$ is an
approximable point. More generally,  from Remark
\ref{comparableabsoluta} we have that if $\inf_{P\in{\cal P}_0}
\lan(T(P))>0$ and $\sup_{P\in {\cal P}_0} \diam (T(P)) <\infty$,
then any point in $X_0^+$ is an approximable point.
\end{remark}

The next theorem  contains a lower bound for the Hausdorff and the grid Hausdorff  dimensions of ${\cal W}(U,x_0,\{r_n\})$ with respect to the grid $\Pibf=\{{\cal P}_n\}$.
As we mentioned in Section \ref{grids},
in order to get results for the $\lan$-Hausdorff dimension we need an extra property of regularity. More precisely, we ask $\Pibf$ to be $\lan$-regular (see Definiton \ref{gridregular}). We recall that any grid on $\R$ is $\lan$-regular.

\medskip

\begin{theorem} \label{dimension} Let $(X,d,{\cal A},\lambda,T)$ be an
expanding system with finite entropy $H_\mu({\cal P}_0)$
with respect to the partition ${\cal P}_0$ where $\mu$ is the
ACIPM associated to the system. Let us consider
the grid $\Pibf=\{{\cal P}_n\}$. Let
$\{r_n\}$ be a non increasing sequence of positive numbers and let
$U$ be an open set in $X$ with $\mu(U)>0$. Then, for all approximable point $x_0\in X_0$, the grid Hausdorff dimensions of the set
$$
{\cal W}(U,x_0,\{r_n\})=\{ x\in U\cap X_0: \  d(T^n (x),x_0) < r_n \ \mbox{for infinitely many} \ n \}
$$
verify
\begin{equation} \label{dimgridporabajo}
\Dim_{\Pibf,\lan} ({\cal W}(U,x_0,\{r_n\})) = \Dim_{\Pibf,\mu} ({\cal W}(U,x_0,\{r_n\})) \ge \frac {h_{\mu}}{h_\mu +
\dimsup \, \ellsup} \,.
\end{equation}
where $\ellsup= \limsup_{n\to\infty} \frac 1n\log\frac{1}{r_n}$
and $h_\mu$ is the entropy of $T$ with respect to $\mu$.

\smallskip

Moreover, for all approximable point $x_0\in X_0$, the Hausdorff dimensions of the set
${\cal W}(U,x_0,\{r_n\})$ verify:

\smallskip

\begin{enumerate}

\item[\rm 1.] If the grid $\Pibf$ is $\lan$-regular then
\begin{equation} \label{dimporabajo}
\Dim_\lan  ({\cal W}(U,x_0,\{r_n\}))=\Dim_\mu  ({\cal
W}(U,x_0,\{r_n\})) \ge \frac {h_{\mu}}{h_\mu + \dimsup \, \ellsup} \, \Big(
1-\frac{\overline{\tau}_\lan(x_0)\dimsup\ell^2}{h_{\mu}^2\log\beta}\Big)
\,,
\end{equation}
\item[\rm 2.] If $\lan$ is a doubling measure verifying that $\lan(B(x,r))\le C\,r^s$ for all ball $B(x,r)$, then
\begin{equation} \label{dimdiametros}
\Dim_\lan  ({\cal W}(U,x_0,\{r_n\}))=\Dim_\mu  ({\cal W}(U,x_0,\{r_n\}))\ge 1 - \frac {\dimsup\ellsup}{s\log\beta}\,.
\end{equation}
\end{enumerate}
Here $\beta$ is the constant appearing in the property {\rm (C)} of expanding
maps.
\end{theorem}

\

\begin{remark} \rm
Recall from Remark \ref{entropia} that $0<h_\mu<\infty$. Recall
also that from Remark \ref{approx} we know that the set of
approximable points has full $\lan$-measure, and from Lemma
{\rm\ref{defdeX1}} we know that all point $x_0$  in $X_1$
satisfies $\dimsup<\infty$ and $\overtau=0$. Therefore since $X_1$
has full $\lan$-measure we have that Theorem {\rm \ref{dimension}}
holds  with $\overtau=0$ for $\lan$-almost all $x_0\in X$.
\end{remark}

\begin{remark} \label{opeque–a}
\rm From Remark \ref{approxfinito} we have that if $\inf_{P\in{\cal P}_0} \lan(T(P))>0$
and $\sup_{P\in {\cal P}_0} \diam (T(P))<\infty$ then any point in
$X_0^+$ is an approximable point. Moreover, if $\sup_{P\in{\cal
P}_0}\diam(T(P))<\infty$ then,
 by Proposition \ref{elchachi},
$$
\frac {\lan(P(n,x_0))}{\lan(P(n+1,x_0))} \asymp \frac
{\lan(T(P(0,T^n(x_0))))}{\lan(P(0,T^{n+1}(x_0)))} \leq \frac
{\lan(X)}{\lan(P(0,T^{n+1}(x_0)))}
$$
and then $\overtau=0$ for all $x_0$ such that
$$
\log \frac 1{\lan(P(0,T^n(x_0)))} =o(n) \,, \qquad \text{ as  }
\quad n\to\infty\,.
$$
\end{remark}

\

First, let us observe that the $\lan$-Hausdorff dimension and $\mu$-Hausdorff
dimensions coincide for subsets of $\cup_{P\in {\cal P}_0} P$ and, in particular, for
subsets of $X_0$.

\begin{lemma} \label{dimensioncoincide}
If $A\in {\cal A}$ is a subset of $\cup_{P\in {\cal P}_0} P$, then
$$
\text{\rm Dim}_{\Pibf,\lan} (A) = \text{\rm Dim}_{\Pibf,\mu} (A)
\qquad \hbox{and} \qquad  \Dim_\lan (A) = \Dim_\mu (A) \,.
$$
\end{lemma}

\begin{proof} We will prove only the equality of grid-dimensions,
since the other proof is similar.
By properties (A.2) and (A.3) of expanding maps we have for the
$\alpha$-dimensional $\lan$-grid and $\mu$-grid Hausdorff measures
that
$$
{\cal H}_{\Pibf,\lan}^\alpha (A) = \sum_{i}  {\cal
H}_{\Pibf,\lan}^\alpha (A\cap P_i)\,, \qquad {\cal
H}_{\Pibf,\mu}^\alpha (A) = \sum_{i}  {\cal H}_{\Pibf,\mu}^\alpha
(A\cap P_i)\,.
$$
where ${\cal P}_0=\{P_i\}$. As a consequence of part (iii) of
Theorem E we get that
$$
{\cal H}_{\Pibf,\lan}^\alpha (A\cap P_i) \asymp  {\cal
H}_{\Pibf,\mu}^\alpha (A\cap P_i) \,,
$$
with constants depending on $i$. Therefore
$$
{\cal H}_{\Pibf,\lan}^\alpha (A) =0 \qquad  \Longleftrightarrow
\qquad {\cal H}_{\Pibf,\mu}^\alpha (A) =0 \,.
$$
\end{proof}

\begin{proof} [Proof of Theorem $\ref{dimension}$]
We may assume that $\dimsup$, $\overtau$ and $\ellsup$ are all finite,
since otherwise
the estimations (\ref{dimgridporabajo}) and (\ref{dimporabajo}) are trivial.
Since $\mu(U)>0$, the set $U$ contains a point $x\in X_0$. As $U$
is open, we have that there exists $r>0$ such that $B(x,r)\subset
U$, where by  $B(x,r)$ we denote the ball $\{y\in X : d(y,x)<
r\}$. Therefore, since $\sup_{P\in{\cal P}_n} \hbox{diam}(P)$ goes
to zero as $n\to\infty$,  there exists $N_0\in \N$ such that
$$
P(N_0,x)\subset B(x,r)\subset U.
$$

Let us write $J_0=P(N_0,x)$ and let $A_0$ denote the element of
the partition ${\cal P}_0$ such that $T^{N_0}(J_0)=A_0$.  To get
the desired result, it is enough to show (\ref{dimgridporabajo}) and (\ref{dimporabajo}) for
the set
$$
\widetilde W=\{ x\in J_0\cap X_0: d(T^n(x),x_0) < r_n  \ \mbox{for infinitely many}
\ n\}.
$$

Notice that we can assume that $\lim_{n\to\infty} r_n=0$. Otherwise, there exists $C>0$ such that
$r_n >C$ for all $n$, and since $\mu$ is ergodic, from Theorem A' in \cite{FMP}, we deduce that for $\mu$-almost every point in $J_0$,
\begin{equation} \label{inequalities}
d(T^n(x),x_0) < C  \ \mbox{for infinitely many} \ n \,.
\end{equation}
Using Remark \ref{fullsets} we conclude that (\ref{inequalities})
holds also for   $\lan$-almost every point in $J_0$ and therefore
$\Dim_{\Pibf,\lan} ({\cal W}(U,x_0,\{r_n\}))=\Dim_{\lan} ({\cal
W}(U,x_0,\{r_n\}))=1$. However,  in this case more is true, see Corollary
\ref{recurrentes}.

To obtain (\ref{dimgridporabajo}) we will construct, for each small $\ep>0$, a Cantor-like set $\cal C_\ep$$ \subset \widetilde W$ and we will prove using Corollary
\ref{azulesyrojos2} that
\begin{equation} \label{DiminfCantor}
\Dim_{\Pibf,\lan} ({\cal C}_\ep) = \Dim_{\Pibf,\mu} ({\cal
C}_\ep)\ge \frac {h_\mu
-2\ep}{h_\mu+2\ep+(1+\ep)(\dimsup+\ep) \, (\ellsup+\ep)+\ep}
\, .
\end{equation}

\

We construct now the Cantor-like set
$$
{\cal C}_{\ep}=\bigcap_{n=0}^{\infty}\bigcup_{J\in{\cal J}_n} J
$$
as follows: we start  with ${\cal J}_0=\{ J_0\}$ and we denote by
$G^{-1}_{J_0}$ the composition of the $N_0$ branchs of $T^{-1}$ such that
$J_0=G^{-1}_{J_0} (A_0)$.

Let ${\cal I}(x_0)=\{p_i\}$ denote the sequence associated to the approximable point $x_0$ given by Definition \ref{aproximable}.
Let $s_k=\diam({P}(p_k,x_0))$, and for each $s_k$ let $n(k)$ denotes the greatest natural number
such that $s_k \le r_{n(k)}$. We denote by $\cal D$  the set of these indexes $n(k)$.
Since $s_k\to 0$ as $k\to\infty$ and $r_n \to 0$ as $n\to\infty$ by hypothesis,
we have that $n(k)\to\infty$ as $k\to\infty$. We will write ${\cal D}=\{d_i\}$ with $d_i<d_{i+1}$ for all $i$.

Notice that  if $d\in{\cal D}$, then there exists $k(d)\in {\cal I}(x_0)$ such that
\begin{equation}\label{permitidos2}
r_{d+1} < \diam(P(k(d),x_0)) \le r_{d}
\end{equation}
and
\begin{equation}\label{permitidos1}
P(k(d),x_0)\subset B(x_0,\diam(P(k(d),x_0))\subset B(x_0, r_d)\,,
\end{equation}

Moreover from the property (C) of expanding maps we have that
$$
C_2\beta^{k(d)}\diam (P(k(d),x_0))\leq\diam (P(0,T^{k(d)} (x_0)))\le
\sup_{P\in {\cal P}_0} \diam(P) <\infty
$$
and using (\ref{permitidos2}) we get
$$
\frac{k(d)}{d+1}\leq\frac{1}{\log\beta}\left [\frac{C}{d+1}+\frac{1}{d+1}\log\frac{1}{r_{d+1}}\right ].
$$
Hence, for all $d$ large enough
\begin{equation} \label{kvsN}
\frac{k(d)}{d}\leq (1+\ep)\frac{\ellsup}{\log\beta}\,.
\end{equation}

To construct the family ${\cal J}_1$  we first choose a natural
number $N_1$  so that $d_1:=N_0+N_1\in\cal D$, and  large enough so that  (\ref{Sn}) holds with
$P_1=A_0$, $P_2=P(0,x_0)$, $N=N_1$ and $M:=M_1=[\ep N_0]$,  (\ref{kvsN}) holds for $d_1$, and also
\begin{equation} \label{d1}
d_1\leq (1+\ep)N_1\,, \qquad r_{d_1} > e^{-d_1(\ellsup+\ep)}.
\end{equation}
Let ${\cal S}_{N_1,M_1}$ denote the collection of  elements  in ${\cal
P}_{N_1}$ given by Proposition \ref{SMsinmalos}. We define
$\widetilde{\cal J}_1$ as the family of sets $G^{-1}_{J_0}(S)$
with $S\in{\cal S}_{N_1,M_1}$, see Figure 1.

\

\begin{center}
\includegraphics[scale=0.45]{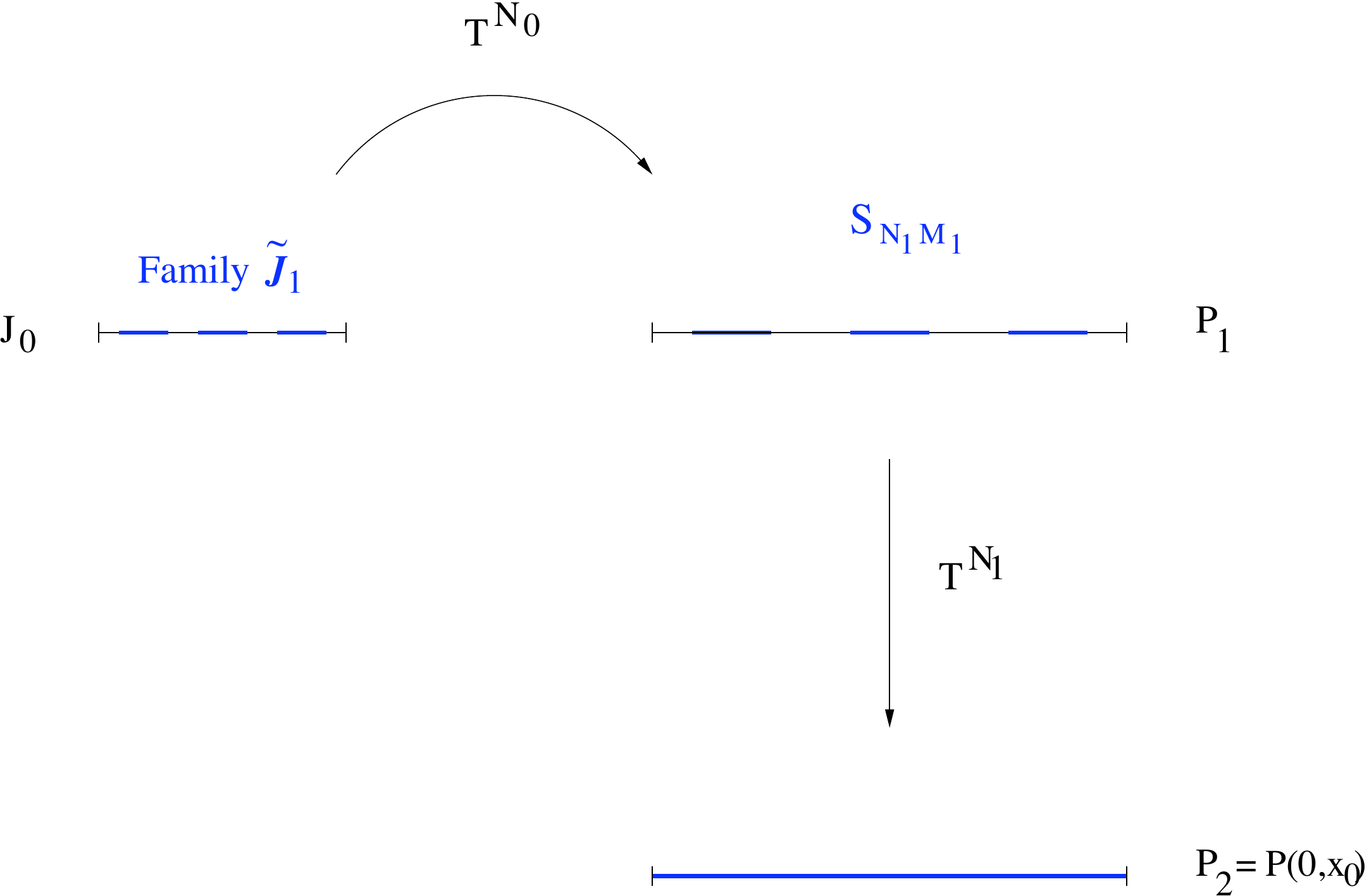}
\end{center}
\centerline{\it Figure $1$.}

\

Notice that by construction if $\widetilde J\in\widetilde{ \cal J}_1$, then
\begin{equation}\label{d1tilde}
T^{d_1}(\widetilde J)=P(0,x_0),\qquad \text{ with }\qquad d_1=N_0+N_1
\end{equation}
and since $T^{N_0}(\widetilde J)=S$ for some $S\in{\cal S}_{N_1,M_1}$
$$
T^{N_0}(\widetilde J)\cap (X_0\setminus E_{M_1}^\ep)\neq\emptyset
\qquad \text{ with } \quad M_1=[\ep N_0].
$$

 We remark that we will define later the
family ${\cal J}_1$ by taking an appropriate  subset of each one
of the elements of the family $\widetilde{\cal J}_1$.
From Proposition  \ref{elchachi}  we obtain that if $\widetilde J=G_{J_0}^{-1}(S)$ with $S\in{\cal S}_{N_1,M_1}$, then
\begin{equation}\label{cociente1}
\frac 1C \, \frac{\lambda(S)}{\lambda(A_0)}\leq \frac{\lambda(
\widetilde J)}{\lambda(J_0)}\leq C
\,\frac{\lambda(S)}{\lambda(A_0)}
\end{equation}
with $C$ an absolute constant. Hence,  from (\ref{SM}),
(\ref{cociente1}),  part (iii) of Theorem E and by taking $N_1$ large enough we get that for all
$\widetilde{J}_1\in \widetilde{\cal J}_1$
$$
 e^{-N_1(h_{\mu}+2\ep)}\le \frac
{\lambda(\widetilde{J}_1)}{\lambda(J_0)} \le
e^{-N_1(h_{\mu}-2\ep)}\,,
$$
From (\ref{cociente1}) we also get an estimate on the size of the family  ${\widetilde{\cal J}}_1$
$$
\lambda(\widetilde{\cal J}_1\cap J_0)=\lambda(\widetilde{\cal J}_1):=\sum_{\widetilde{J}_1\in\widetilde {\cal J}_1} \lambda(\widetilde{J}_1)  \geq \frac 1C \frac{\lambda(J_0)}{\lan(A_0)}\sum_{S\in{\cal S}_{N_1,M_1}}\lambda(S)\,.
$$
Therefore, from Proposition \ref{SMsinmalos}, and the part (iii)
of Theorem E  we get that
$$
\lambda(\widetilde{\cal J}_1 \cap J_0) \geq C\,
\lambda(J_0)\,\lambda (P(0,x_0))
$$
with $C>0$ depending on $P_1=A_0$ and $P_2=P(0.x_0)$.

\

Now, since $d_1=N_0+N_1\in\cal D$, by (\ref{permitidos2}) and  (\ref{permitidos1}) there exists an integer $k_1\in {\cal I}(x_0)$ such that
\begin{equation} \label{defdek1}
P(k_1,x_0)\subset  B(x_0, \,r_{d_1}).
\end{equation}
and
\begin{equation}\label{diametro}
\hbox{diam}(P(k_1,x_0))\geq r_{d_1+1}\,.
\end{equation}

Moreover, from (\ref{d1})  we have that $k_1\leq \frac{(1+\ep)^2\ellsup}{\log\beta}N_1$.
Since  we can take $N_1$ as
large as we want so that $k_1$ is large enough we can get that
\begin{equation}\label{clauint1}
\hbox{closure\,}(P(k_1,x_0)) \subset P(0,x_0)
\end{equation}
and also by taking $N_1$ large we
have  from the definition of $\dimsup$ and (\ref{diametro}) and (\ref{d1}) that
\begin{equation}\label{Ik1}
\lambda(P(k_1,x_0))\ge \hbox{diam} (P(k_1,x_0))^{\dimsup+\ep}
\,\geq   r_{d_1+1}^{\dimsup + \ep} \ge e^{-N_1(1+\ep)(\dimsup+\ep)(\ellsup+\ep)}.
 \end{equation}

\

For all $S\in{\cal S}_{N_1,M_1}$ let $G^{-1}_{S}$ denote the composition of the $N_1$ branchs of $T^{-1}$ such that $S=G^{-1}_{S}(P(0,x_0))$. In each set  $S$ in ${\cal S}_{N_1,M_1}$ we  take the subset $L_1:=G^{-1}_{S}(P(k_1,x_0))$ and we denote by ${\cal L}_1$ this family of sets.
To define the family ${\cal J}_1$  we just ``draw'' the sets $L_1$
in $J_0$. More precisely, ${\cal J}_1$ is the family
$G^{-1}_{J_0}(L_1)$ with $L_1\in{\cal L}_1$, see Figure 2.

\begin{center}
\includegraphics[scale=0.45]{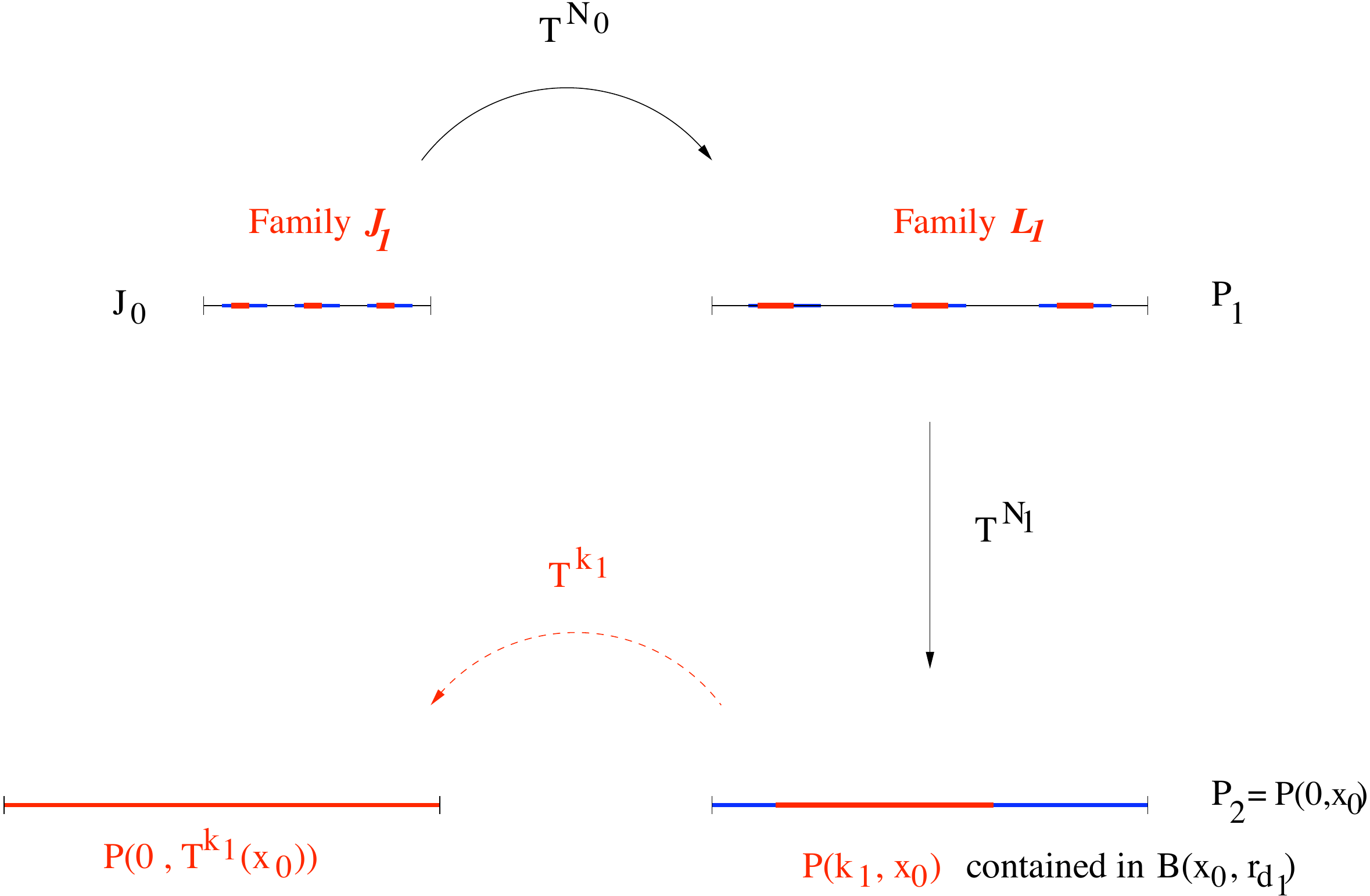}
\end{center}
\centerline{\it Figure $2$.}

Notice that by construction if $ J\in{ \cal J}_1$, then
\begin{equation}\label{d1dej1}
T^{d_1}(J)=P(k_1,x_0) \qquad \text{ and }\qquad
T^{d_1+k_1}(J)=P(0,T^{k_1}(x_0)).
\end{equation}

Hence if $x\in J \in {\cal J}_1$ then $T^{d_1}(x)\in
P(k_1,x_0)\subset B(x_0, \,r_{d_1})$, and   it follows that
$$
d(T^{d_1}(x),x_0)\leq  r_{d_1} \,.
$$

\

By construction we have that for all $ J\in{ \cal J}_1$
there exists an unique $\widetilde J\in\widetilde{ \cal J}_1$ such that $J\subset \widetilde J$,  and by using the condition (A.6) and (\ref{clauint1}) we have that
$$
\hbox{closure}(J_1)=G_{J_0}^{-1}(G_{S}^{-1}(\hbox{closure}(P(k_1,x_0))))\subset
G_{J_0}^{-1}(G_S^{-1}((P(0,x_0))) = \widetilde J_1\,.
$$

Also by (\ref{d1tilde}),  (\ref{d1dej1}) and Proposition \ref{elchachi} we get that
$$
 \frac 1C \, \lan(P(k_1,x_0))\leq \frac{\lan(J_1)}{\lan(\widetilde J_1)} \leq C \lan(P(k_1,x_0))
$$
with $C>0$ a constant depending on $\lan (P(0,x_0))$. And from
(\ref{Ik1}) by taking $N_1$ large, we have that
$$
 \frac{\lan(J_1)}{\lan(\widetilde J_1)} \geq e^{-N_1[(1+\ep)(\dimsup+\ep)(\ellsup+\ep)+\ep]}.
$$

\

Now, let us assume that we have already constructed the families
$\widetilde{\cal J}_j, \, {\cal J}_j$ and the numbers $N_j$ and
$k_j\in {\cal I}(x_0)$  for $j=1,\ldots,m$ with the following properties:

\medskip

Let $d_1=N_0+N_1$ and
$$d_j:=N_0+N_1+\cdots +N_j+k_1+\cdots
+k_{j-1}\qquad \text{ for } \qquad  j\geq 2$$.
\begin{itemize}
\item[(a)]  For all point $x$ in $J_j\in{\cal J}_j$
$$
d(T^{d_j}(x),x_0) \leq  r_{d_j}.
$$
\item[(b)] For all $\widetilde{J}_j\in\widetilde{\cal J}_j$ we have

\begin{itemize}
\item[(b1)] $ T^{d_j}(\widetilde{J}_j)=P(0,x_0)$  and
$$
T^{d_j-N_j}(\widetilde{J}_j)\bigcap(X_0\setminus E_{M_j}^{\ep})\neq\emptyset \qquad \text{ with }\quad M_j=[\ep N_{j-1}].
$$
\item[(b2)]
 There exists a unique $J_{j-1}\in{\cal J}_{j-1}$ so that $\widetilde{J}_j\subset J_{j-1}$ and
$$
e^{-N_j(h_{\mu}+2\ep)}\le \frac
{\lan(\widetilde{J}_j)}{\lan (J_{j-1})} \le
e^{-N_j(h_{\mu}-2\ep)}\,.
$$
\end{itemize}
\item[(c)] For
all $J_j\in{\cal J}_j$ we have
\begin{itemize}
\item[(c1)]  $T^{d_j}(J_j)=P(k_j,x_0)$.
\item[(c2)] There exists a unique $\widetilde J_{j}\in\widetilde {\cal J}_{j}$ so
that $\hbox{closure}(J_j) \subset \widetilde J_{j}$,
$$
 \frac {\lan(J_j)}{\lan(\widetilde J_{j})}\asymp \lan(P(k_j,x_0)) \qquad \hbox{and} \qquad
 \frac {\lan(J_j)}{\lan(\widetilde J_{j})}\geq
 e^{-N_j[(1+\ep)(\dimsup+\ep)(\ellsup+\ep)+\ep]}\,.
$$
Besides, for each  $\widetilde J_{j}\in\widetilde {\cal J}_{j}$ there exists a unique $J_{j}\in  {\cal J}_{j}$ so that $\hbox{closure}(J_j) \subset \widetilde J_{j}$
\item[(c3)] $k_j\leq\frac{(1+\ep)^2\ellsup}{\log\beta}N_j$.
\end{itemize}
\item[(d)] There exists an absolute constant $\widetilde{c}>1$ such that for
all $J_{j-1}\in {\cal J}_{j-1}$,
$$
\lan(\widetilde{\cal J}_j\cap
J_{j-1}):=\sum_{\widetilde{J}_j\in\widetilde{\cal J}_j, \,
\widetilde {J}_j\subset J_{j-1}} \lan(\widetilde{ J}_j)\geq  \frac
1{\widetilde{c}}\,  \lan(J_{j-1}) \,.
$$
\item[(e)] $N_j$ is big enough so that
$$
\frac{ j}{N_1+\cdots+N_j} < \frac 1{j}\,.
$$
\end{itemize}

We want to mention that the hypothesis on $x_0$ of being approximable is only required to obtain  an absolute constant $\widetilde{c}$ in the property (d).

Recall that we want to apply Corollary  \ref{azulesyrojos2}. In our case
\begin{equation} \label{enemigomio}
a=h_\mu+2\ep\,, \qquad
b=h_\mu-2\ep\,, \qquad
c=(1+\ep)(\dimsup+\ep)(\ellsup+\ep)+\ep \qquad \hbox{and} \qquad
\delta=1/\widetilde{c}\,.
\end{equation}

\

Now we start with the construction of the family ${\cal J}_{m+1}$.
We choose a natural number $N_{m+1}\ge N_m$ large enough so that
$$
d_{m+1}:=N_0+N_1+\cdots +N_{m+1}+k_1+\cdots + k_{m}\in {\cal D},
$$
property (e) holds for $j=m+1$, (\ref{Sn}) holds with $P_1=P(0,T^{k_m}(x_0))$, $P_2=P(0,x_0)$, $N=N_{m+1}$,  and
$M:=M_{m+1}=[\ep N_{m}]$, (\ref{kvsN}) holds for $d_{m+1}$, and also
\begin{equation} \label{historiainterminable}
d_{m+1}\leq (1+\ep) N_{m+1}\,, \qquad r_{d_{m+1}} \ge e^{-d_{m+1}(\ellsup+\ep)}\,.
\end{equation}

Let ${\cal S}_{N_{m+1},M_{m+1}}$ denote the collection of elements
in ${\cal P}_{N_{m+1}}$ given by Proposition \ref{SMsinmalos}.
Notice that the sets in this family verify (\ref{SM}) with
$N=N_{m+1}$. For each $J\in{\cal J}_m$ let $G_J^{-1}$ denote the
composition of the $d_m+k_m$ branchs of $T^{-1}$ such that
$J=G_J^{-1}(P(0,T^{k_m}(x_0)))$.  We define now $\widetilde{\cal
J}_{m+1}$ as
$$
\widetilde{\cal J}_{m+1} = \bigcup_{J\in {\cal J}_m} G_J^{-1} ({\cal S}_{N_{m+1},M_{m+1}}) \,.
$$
Notice that, by construction, if $\widetilde J\in\widetilde{ \cal J}_{m+1}$, then
\begin{equation}\label{dm1tilde}
T^{d_{m+1}}(\widetilde J)=P(0,x_0),
\end{equation}
and since $T^{d_{m+1}-N_{m+1}}(\widetilde J)=S$ for some $S\in{\cal S}_{N_{m+1},M_{m+1}}$.
$$
T^{d_{m+1}-N_{m+1}}(\widetilde J)\cup (X_0\setminus E_{M_{m+1}}^\ep)\neq\emptyset
\qquad \text{ with } \quad M_{m+1}=[\ep N_m].
$$
Since $J_m\in {\cal P}_{d_m+k_m}$, by Proposition \ref{elchachi} we have
that if $\widetilde{J}= G_{J_m}^{-1}(S)$ with $S\in{\cal
S}_{N_{m+1},M_{m+1}}$, then
\begin{equation}\label{propmt}
\frac 1C \, \frac{\lan(S)}{\lan(P(0,T^{k_m}(x_0)))}\leq
\frac{\lan(\widetilde J)}{\lan(J_m)}\leq C \,
\frac{\lan(S)}{\lan(P(0,T^{k_m}(x_0)))} \,,
\end{equation}
with $C$ an absolute constant.

If $\widetilde{J}_{m+1}\in\widetilde {\cal J}_{m+1}$, then there
are $J_m \in {\cal J}_m$ and $S\in{\cal S}_{N_{m+1},M_{m+1}}$ such
that $\widetilde{J}_{m+1}= G_{J_m}^{-1}(S)$. Hence from
(\ref{propmt}), (\ref{SM}), the definition of approximable
points and by taking $N_{m+1}$ large, we get the property (b) for $j=m+1$. We remark that
$\lan(P(0,T^{k_m}(x_0))$ does not depend on $N_{m+1}$.

Now from (\ref{propmt}), Proposition \ref{SMsinmalos} and again
the definition of approximable points we get
$$
\lan(\widetilde{{\cal J}}_{m+1}\cap J_m) = \sum_{\substack{\widetilde{J}_{m+1}\in\widetilde {\cal J}_{m+1} \\ \widetilde{J}_{m+1}\subset J_m}} \!\!\!\!\lan(\widetilde{J}_{m+1}) \ge
\frac 1C \frac {\lan(J_m)}{\lan(P(0,T^{k_m}(x_0)))} \sum_{S\in
{\cal S}_{N_{m+1},M_{m+1}}}\!\!\!\!\!\!\! \lan(S) \ge \frac 1{c'}
\,\lan(P(0,x_0))\,\lan(J_m) \,,
$$
and this gives property (d) for $j=m+1$. Observe that the constant
$c'$ depends on the comparability constant between $\lan$ and $\mu$
in $P_1=P(0,T^{k_m}(x_0))$ but from the definition of
approximable points we know that this constant is absolute.

Since $d_{m+1}\in\cal D$, by (\ref{permitidos2}) and (\ref{permitidos1}), there exists an integer
$k_{m+1}\in {\cal I}(x_0)$ such that
\begin{equation} \label{defdekm}
P(k_{m+1},x_0)\subset  B(x_0, r_{d_{m+1}}).
\end{equation}
and
\begin{equation}\label{diametrom}
\hbox{diam}(P(k_{m+1},x_0))\geq r_{d_{m+1}+1}\,.
\end{equation}

From (\ref{kvsN}) and since $d_{m+1}\leq (1+\ep)N_{m+1}$ we get the property (c3) for $j=m+1$.

\

As in the initial step from the definition of $\dimsup$, (\ref{diametrom}) and (\ref{historiainterminable}),
we have that
\begin{equation}\label{Ikm+1}
\lambda(P(k_{m+1},x_0))\ge  r_{d_{m+1}+1}^{\dimsup + \ep} \ge e^{-N_{m+1}(1+\ep)(\dimsup+\ep)(\ellsup+\ep)}.
 \end{equation}

In each set $S\in{\cal S}_{N_{m+1},M_{m+1}}$ we take the subset
$L_{m+1}:=G^{-1}_{S}(P(k_{m+1},x_0))$ and we call ${\cal L}_{m+1}$
to this family of sets. We recall that by $G_S^{-1}$ we denote the
composition of the $N_{m+1}$ branchs of $T^{-1}$ such that
$S=G_S^{-1}(P(0,x_0))$.

To define the family ${\cal J}_{m+1}$ we ``draw'' the family
${\cal L}_{m+1}$ in each one of the sets $J\in{\cal J}_m$. More
precisely, for each  $J\in{\cal J}_m$ let $G_J$ denote the
composition of the $d_m+k_m$ branchs of $T$  such that
$G_J(J)=P(0,T^{k_m}(x_0))$.  We define now ${\cal J}_{m+1}$ as
$$
{\cal J}_{m+1} = \bigcup_{J\in {\cal J}_m} G_J^{-1} ({\cal L}_{m+1}) \,.
$$

Notice that by construction if $ J\in{ \cal J}_{m+1}$, then
\begin{equation}\label{d1dejm}
T^{d_{m+1}}(J)=P(k_{m+1},x_0) \qquad \text{ and }\qquad
T^{d_{m+1}+k_{m+1}}(J)=P(0,T^{k_{m+1}}(x_0)).
\end{equation}
Therefore the condition (c1) holds for $j=m+1$. Besides, by (\ref{defdekm}), if $x\in J
\in {\cal J}_{m+1}$ then $T^{d_{m+1}}(x)\in P(k_{m+1},x_0)\subset
B(x_0,r_{d_{m+1}})$, and   therefore the condition (a) holds for
$j=m+1$.

By construction we have that for all $ J_{m+1}\in{ \cal J}_{m+1}$
there exists an unique $\widetilde J_{m+1}\in\widetilde{ \cal J}_{m+1}$ such that $J_{m+1}\subset \widetilde J_{m+1}$,  and by using the condition (A.6) and
(\ref{clauint1}) as in the initial step we have that
$$
\hbox{closure}(J_{m+1}) \subset \widetilde J_{m+1}\,.
$$
The estimates of $\lan(J_{m+1})/\lan(\widetilde J_{m+1})$ of the condition (c2) follows by applying
Proposition \ref{elchachi} and by  using (\ref{dm1tilde}),  (\ref{Ikm+1}) and  (\ref{d1dejm}).

\

We have already obtained the properties (a)-(e) for $j=m+1$, and
therefore we have concluded the construction of the Cantor-like
set ${\cal C}_{\ep}$. The property that for all $J_{m+1}\in{\cal J}_{m+1}$ there exists a unique  $J_m\in{\cal J}_m$ such that
$$
\hbox{closure}(J_{m+1}) \subset J_{m}
$$
implies that  ${\cal C}_{\ep}$ is not empty. And moreover  ${\cal
C}_{\ep}\subset X_0$ since by construction $P(m,x)$ is defined for
all $x\in  {\cal C}_{\ep}$. Hence the condition (a) implies that
${\cal C}_{\ep}$ is contained in the set $W$.

\

The estimate (\ref{DiminfCantor}) follows now directly from (\ref{enemigomio}), property (e) and Corollary \ref{azulesyrojos2}.


\

Next we will prove the estimate (\ref{dimporabajo}) for the
$\lan$-grid Hausdorff dimension of ${\cal C}_\ep$. We will use  the
subcollections $\{{\cal Q}_m\}$ of $\{{\cal P}_m\}$ given by
$$
{\cal Q}_m=\{ P(m,x) \, : \, x\in {\cal C}_{\ep}\}
$$
in order to apply Proposition \ref{haches}.
Since  $\Pibf$ is a $\lan$-regular grid, see Definition \ref{gridregular},
we only need to deal with the computation of the parameters
$\{a_m\}$ and $\{b_m\}$ of the subcolletions ${\cal Q}_m$. We recall that $a_m$
and $b_m$ are,  respectively,   a lower and an upper  bound for
$\lan(P(m,x))$ with $x\in{\cal C}_{\ep}$.

The easiest cases correspond to $m=d_n$ and $m=d_n+k_n$, i.e.  to the families $\widetilde{\cal J}_n$ and ${\cal J}_n$.
Since $P(d_n,x)$ belongs to $\widetilde{\cal J}_n$, from property (b2) of  ${\cal C}_{\ep}$ and by taking $N_n$ large enough  we have that
\begin{equation}\label{abdn}
a_{d_n}=e^{-N_n(h_{\mu}+3\ep)} \qquad \text{ and }\qquad
b_{d_n}=e^{-N_n(h_{\mu}-3\ep)}\,.
\end{equation}
Also, from property (c2) for $j=n$,
$$
\lan(P(d_n+k_n,x))\asymp \lan(P(d_n, x))\, \lan(P(k_n,x_0))\,.
$$
for all $x\in{\cal C}_{\ep}$, and therefore,
\begin{equation} \label{abdnmaskn}
a_{d_n+k_n}\asymp  a_{d_n}\lan(P(k_n,x_0)) \qquad \text{ and }\qquad
b_{d_n+k_n}\asymp b_{d_n}\lan(P(k_n,x_0))\,.
\end{equation}

To estimate $a_m$ and $b_m$ in the other cases  we need first  some estimate on the Jacobian. Specifically we need to estimate $\J_{d_n}(x)$ and $\J_{d_n+k_n}(x)$ for $x\in{\cal C}_{\ep}$. From (\ref{njacobiano}), Proposition \ref{cotajacobiano},  and properties (b1) and (c1) of ${\cal C}_{\ep}$ (for $ j=n$) we have that
$$
\lan(P(0,x_0))=\lan(T^{d_n}(P(d_n,x)))=\int_{P(d_n,x)}\J_{d_n}
d\lan \asymp \J_{d_n}(x)\, \lan(P(d_n,x))
$$
and
$$
\lan(P(0,T^{k_n}(x_0)))=\lan(T^{d_n+k_n}(P(d_n+k_n,x)))=\int_{P(d_n+k_n,x)}\!\!\!\!\J_{d_n+k_n}
d\lan \asymp \J_{d_n+k_n}(x)\, \lan(P(d_n+k_n,x))\,.
$$
Hence for all $x\in{\cal C}_{\ep}$
\begin{equation}\label{otrojacobiano1}
\frac{1}{\J_{d_n}(x)}\asymp\lan(P(d_n,x))
\end{equation}
and, by using property (c2) for $j=n$,
\begin{equation} \label{otrojacobiano2}
\frac{1}{\J_{d_n+k_n}(x)}\asymp \frac{\lan(P(d_n,x)) \,
\lan(P(k_n,x_0))}{\lan(P(0,T^{k_n}(x_0)))}
\end{equation}
with constants depending on $P(0,x_0)$.

\

For $d_n<m < d_n+k_n$ by  (\ref{njacobiano}) and Proposition \ref{cotajacobiano} we have that
$$
\lan(P(m-d_n,T^{d_n}(x)))=\lan(T^{d_n}(P(m,x)))=\int_{P(m,x)}\J_{d_n}
d\lan \asymp \J_{d_n}(x)\, \lan(P(m,x))\,.
$$
Then from (\ref{otrojacobiano1}) we get
$$
\lan(P(m,x))\asymp \lan(P(d_n,x))\, \lan(P(m-d_n, T^{d_n}(x)))\,.
$$
But, since $x\in{\cal C}_{\ep}$, by (c1)
$$
T^{d_n}(x)\in P(k_n,x_0)\subseteq P(m-d_n, x_0), \qquad \text{ for
} m\le d_n+k_n,
$$
and therefore $P(m-d_n, T^{d_n}(x))=P(m-d_n,x_0)$.
Hence we have that for $d_n \le m \le d_n+k_n$
\begin{equation}\label{abdn<dnkn}
a_{m}\asymp  a_{d_n}\lan(P(m-d_n,x_0)) \qquad \text{ and }\qquad
b_{m}\asymp b_{d_n}\lan(P(m-d_n,x_0))\,.
\end{equation}

Now for $d_n+k_n<m < d_{n+1}=d_n+k_n+N_{n+1}$ by  (\ref{njacobiano}) and  Proposition \ref{cotajacobiano} we have that
$$
\lan(T^{d_n+k_n}(P(m,x)))=\int_{P(m,x)}\J_{d_n+k_n} d\lan \asymp
\J_{d_n+k_n}(x)\, \lan(P(m,x))
$$
and from (\ref{otrojacobiano2}) we get
\begin{equation}\label{m<dn1}
\lan(P(m,x))\asymp \frac {\lan(P(d_n,x))\,
\lan(P(k_n,x_0)\,\lan(T^{d_n+k_n}(P(m,x)))}{\lan(P(0,T^{k_n}(x_0)))}\,.
\end{equation}
Therefore, we need to obtain upper and lower bounds of $\lan(T^{d_n+k_n}(P(m,x)))$ independent of $x\in {\cal C}_\ep$.

Notice that if $m \le d_{n+1}$, then
$$
T^{d_{n+1}-N_{n+1}}(P(d_{n+1},x))=T^{d_{n}+k_{n}}(P(d_{n+1},x))\subset
T^{d_n+k_n}(P(m,x))
$$
and, since $T^{d_{n+1}-N_{n+1}}(P(d_{n+1},x))=P(N_{n+1},T^{d_n+k_n}(x))$ is an element of the family ${\cal S}_{N_{n+1},M_{n+1}}$, from the property (b1) of ${\cal C}_{\ep}$ we can conclude
that there exists $z\in T^{d_n+k_n}(P(m,x))$ such that $z\not\in
E_{M_{n+1}}^{\ep}$ with $M_{n+1}=[\ep N_n]$. Hence, for $m\le d_{n+1}$,
\begin{equation}\label{bien}
T^{d_n+k_n}(P(m,x))=P(m-d_n-k_n,z) \, \qquad \hbox {with }\qquad
z\not\in E_{[\ep N_n]}^{\ep}
\end{equation}
and, for $m-d_n-k_n\ge M_{n+1}=[\ep N_n]$,
$$
\frac 1C e^{-(m-d_n-k_n)(h_{\mu}+\ep)}\leq
\lan(T^{d_n+k_n}(P(m,x)))\leq C e^{-(m-d_n-k_n)(h_{\mu}-\ep)}\,.
$$

Therefore, for $d_n+k_n+[\ep N_n]\leq m < d_{n+1}$,
\begin{equation}\label{abdnknnn<dn1}
a_{m}\asymp  \frac {a_{d_n}\lan(P(k_n,x_0))
e^{-(m-d_n-k_n)(h_{\mu}+\ep)}}{\lan(P(0,T^{k_n}(x_0)))} \quad
\text{ and }\quad b_{m}\asymp
\frac{b_{d_n}\lan(P(k_n,x_0))e^{-(m-d_n-k_n)(h_{\mu}-\ep)}}{\lan(P(0,T^{k_n}(x_0)))}\,.
\end{equation}

For $d_n+k_n<m<d_n+k_n+[\ep N_n]$ we have that
$$
P([\ep N_n],z)\subset T^{d_n+k_n}(P(m,x))\subset P(0,T^{k_n}(x_0))
$$
and therefore, from Lemma \ref{egorof} and the definition of
approximable point (recall that $k_n\in{\cal I}(x_0)$),
\begin{equation} \label{sinliminf}
\frac 1C \, e^{-[\ep N_n](h_\mu+\ep)} \le \lan(T^{d_n+k_n}(P(m,x)))\le \lan(P(0,T^{k_n}(x_0))).
\end{equation}
Hence,  from (\ref{m<dn1}) we get that
\begin{equation}\label{abdnkn<dnknnn}
a_m\asymp \frac
{a_{d_n}\lan(P(k_n,x_0))e^{-[\ep N_n](h_{\mu}+\ep)}}{\lan(P(0,T^{k_n}(x_0)))}
\qquad \text{ and }\qquad b_m \asymp b_{d_n}\lan(P(k_n,x_0))\,.
\end{equation}

\

In order to apply Proposition \ref{haches} we will show that
\begin{equation}\label{cocientegrid}
\limsup_{m\to\infty}\frac{\log \left (1/a_m\right)}{\log\left(
1/b_{m-1}\right)}\leq 1+C\ep
+\frac{\overtau\ellsup}{(h_{\mu}-3\ep)\log\beta}
\end{equation}
with $C$ an absolute constant. Then, from Proposition \ref{haches},
(\ref{cocientegrid}), and by taking $\ep\to 0$ we get the desired bound for the
$\lan$-Hausdorff dimension of the set ${\cal W}(U,x_0,\{r_n\})$.

Let  us define
$$
q_m:=\frac{\log \left (1/a_m\right)}{\log\left( 1/b_{m-1}\right)}\, .
$$

For $m=d_n$ we have  by (\ref{abdn}) and (\ref{abdnknnn<dn1})  that
\begin{equation}\label{phm1}
q_m\asymp \frac{N_n (h_\mu+3\ep)}{N_n(h_\mu-\ep)+C_{n-1}}
\end{equation}
with
$$
C_{n-1}=(h_{\mu}-3\ep)N_{n-1}+\log\frac{\lan(P(0,T^{k_{n-1}}(x_0)))}{\lan(P(k_{n-1},x_0))}\,
.
$$

Hence from  (\ref{phm1}) we get that
\begin{equation}\label{hm1}
q_m\leq 1+C\ep \qquad \text{ for } \quad  m=d_n
\end{equation}
with  $C>0$ an absolute constant. In order to obtain
(\ref{cocientegrid}) from (\ref{phm1}) we must also say that we
have taken $N_n$ large enough so that $C_{n-1}\ge -\ep N_n$.

\

For $m=d_n+1$ we get from (\ref{abdn}) and (\ref{abdn<dnkn}) that,
for $N_n$ large enough,
\begin{equation}\label{phm1punto5}
q_m\asymp\frac{N_n(h_\mu+3\ep) +\log\frac
1{\lan(P(1,x_0))}}{N_n(h_\mu-3\ep)}  \le 1+C\ep\,.
\end{equation}

For $d_n+1<m\leq d_n+ k_n$ we have from (\ref{abdn<dnkn})
\begin{equation}\label{phm2}
q_m\asymp\frac{ N_n(h_\mu+3\ep)+\log\frac 1{\lan(P(m-d_n,x_0))}}
{N_n(h_\mu-3\ep)+\log\frac{1}{\lan(P(m-d_n-1,x_0)}} \, .
\end{equation}
But since
$$
\overtau=\limsup_{n\to\infty}\frac{1}{n}\log\frac{\lan(P(n,x_0))}{\lan(P(n+1,x_0))},
$$
then given $\ep>0$ there exists $C'>0$ such that for all $n$
$$
\log\frac{\lan(P(n,x_0))}{\lan(P(n+1,x_0))} \leq C' +
n(\overtau+\ep)\, .
$$
Hence
$$
\log\frac{\lan(P(m-d_n-1,x_0))}{\lan(P(m-d_n,x_0))}\leq
C'+k_n(\overtau+\ep)
$$
and by property (c3) of ${\cal C}_{\ep}$ for $N_n$ large enough we get that
\begin{equation}\label{pphm2}
\log\frac{\lan(P(m-d_n-1,x_0))}{\lan(P(m-d_n,x_0))}\leq
C'+\frac{(\overtau+\ep) (1+\ep)^2\ellsup}{\log\beta} N_n.
\end{equation}
From (\ref{phm1punto5}), (\ref{phm2}) and (\ref{pphm2}) it follows
that for $N_n$ large
\begin{equation}\label{hm2}
q_m\leq 1+C\ep+\frac{\overtau\ellsup}{(h_{\mu}-3\ep)\log\beta}
\qquad \text{ for } \quad d_n< m\leq d_n+k_n
\end{equation}
with  $C>0$ an absolute constant.

\

For $d_n+k_n<m <d_n+k_n+[\ep N_n]$ we get from (\ref{abdnmaskn}),
(\ref{abdnknnn<dn1}) and (\ref{abdnkn<dnknnn}) that
$$
q_m\asymp\frac{N_n(h_\mu+3\ep) +[\ep N_n](h_{\mu}+\ep)+ \log\frac
1{\lan(P(k_n,x_0))}-\log\frac{1}{\lan(P(0,T^{k_n}(x_0)))}}
{N_n(h_{\mu}-3\ep)+\log\frac{1}{\lan(P(k_n,x_0))}}
$$
and therefore, for $N_n$ large enough, we have with $C$ an
absolute constant that
\begin{equation}\label{hm3}
q_m \leq 1+\frac{6\ep N_n+[\ep N_n]
(h_{\mu}+\ep)}{N_n(h_{\mu}-3\ep)}\le 1+C\ep \qquad \text{ for }
\quad d_n+k_n< m < d_n+k_n+ [\ep N_n] \,.
\end{equation}

\medskip

For $m=d_n+k_m+[\ep N_n]$ we get from  (\ref{abdnkn<dnknnn}) and
(\ref{abdnknnn<dn1}) that
$$
q_m\asymp\frac{N_n(h_\mu +3\ep)+ [\ep N_n]( h_{\mu}+\ep)
+\log\frac 1{P(k_n,x_0)}-\log\frac 1{\lan(P(0,T^{k_n}(x_0)))}}
{N_n(h_{\mu}-3\ep)+\log\frac{1}{\lan(P(k_n,x_0))}}
$$
and therefore, for $N_n$ large enough, we have, with $C$ an
absolute constant, that
\begin{equation} \label{hm3ymedio}
q_m \le \frac {6\ep N_n+[\ep N_n](h_\mu
+\ep)}{N_n(h_{\mu}-3\ep)}\le 1+C\ep \qquad \text{ for } \quad
 m = d_n+k_n+ [\ep N_n]\,.
\end{equation}

\medskip

For $d_n+k_n+[\ep N_n]<m <d_{n+1}$ we get from (\ref{abdnknnn<dn1}) that
\begin{equation}\label{phm4}
q_m\asymp\frac{N_n(h_\mu+3\ep)+ (m-d_n-k_n)(h_\mu+\ep)+\log\frac
1{P(k_n,x_0)}-\log\frac 1{\lan(P(0,T^{k_n}(x_0)))}}
{N_n(h_{\mu}-3\ep)+(m-d_n-k_n-1)(h_{\mu}+\ep)+\log\frac
1{P(k_n,x_0)}-\log\frac 1{\lan(P(0,T^{k_n}(x_0)))}}\, .
\end{equation}
Hence from (\ref{sinliminf}) and (\ref{phm4}) we have that
$$
q_m\leq 1+\frac{6\ep N_n+h_\mu+\ep}{N_n(h_{\mu}-3\ep)-\log\frac
1{\lan(P(0,T^{k_n}(x_0)))}}\leq 1+ \frac{6\ep
N_n+h_\mu+\ep}{N_n(h_{\mu}-3\ep)+\log\frac 1C-[\ep
N_n](h_\mu+\ep)} \, \,
$$
and so, for $N_n$ large enough, we have, with $C$ an absolute
constant, that
\begin{equation} \label{hm4}
q_m \le 1+C\ep \qquad \text{ for } \quad d_n+k_n+ [\ep N_n]<m
<d_{n+1}\,.
\end{equation}

From (\ref{hm1}),  (\ref{hm2}),  (\ref{hm3}), (\ref{hm3ymedio})
and (\ref{hm4}) we get (\ref{cocientegrid}). Using now Proposition
\ref{haches} it follows that
\begin{equation} \label{cocientedimensiones}
\frac {1-\Dim_\lan ({\cal C}_\ep)}{1-\Dim_{\Pi,\lan} ({\cal
C}_\ep)}\le 1+C\ep
+\frac{\overtau\ellsup}{(h_{\mu}-3\ep)\log\beta}\,.
\end{equation}
As ${\cal C}_\ep \subset {\cal W}(U,x_0,\{r_n\})$ for all $\ep>0$,
(\ref{dimporabajo}) follows now from (\ref{DiminfCantor}) and
(\ref{cocientedimensiones}) by letting $\ep\to 0$.

\

Finally, to prove (\ref{dimdiametros}) it is enough to show that, for all $x\in{\cal C}_\ep$,
\begin{equation} \label{frostman}
\nu(B(x,r)) \le C\, (\lan(B(x,r))^\eta \,, \qquad \hbox{with } \eta=1-(1+\ep)\frac{(1+\ep)(\dimsup+\ep)(\ellsup+\ep)+\ep}{s\log\beta}\,.
\end{equation}
First, notice that from the definition of the measure $\nu$ and  the properties (c2) and (d) of the definition of the Cantor set ${\cal C}_\ep$  it follows that
for all $x\in{\cal C}_\ep$:

\

(1) If $J_{n+1} \subset P(m,x) \subseteq \widetilde{J}_{n+1}$, then
$$
\nu(P(m,x))=\nu(\widetilde{J}_{n+1})=\nu(J_{n+1})\le C \lan( \widetilde{J}_{n+1})\, \frac{\nu(J_n)}{\lan(J_n)}\le C\, \frac 1{\lan(P(k_{n+1},x_0))}\, \frac{\nu(J_n)}{\lan(J_n)}\,\lan (P(m,x)) \,.
$$

(2) If $\widetilde{J}_{n+1} \subset P(m,x) \subseteq J_n$, then
\begin{align}
\nu(P(m,x)) & =
\sum_{\substack{{\tilde{J}}_{n+1} \in {\widetilde{\cal J}}_{n+1} \\
{\widetilde{J}}_{n+1} \subseteq P(m,x)}}  \nu (J_{n+1}) \,=
\sum_{\substack{{\tilde{J}}_{n+1} \in {\widetilde{\cal J}}_{n+1} \\
{\widetilde{J}}_{n+1} \subseteq P(m,x)}} \frac{\lan({\widetilde{J}_{n+1}})}{\lan({\widetilde{\cal J}}_{n+1}
\cap J_n)}\,\nu (J_{n}) \notag \\
&\le C\,  \frac{\nu (J_{n}) }{\lan( J_n)}\,\sum_{\substack{{\tilde{J}}_{n+1}
\in {\widetilde{\cal J}}_{n+1} \\
{\widetilde{J}}_{n+1} \subseteq
P(m,x)}}\lan({\widetilde{J}_{n+1}}) \leq C\,  \frac{\nu (J_{n})
}{\lan({ J}_{n})}\ \lan(P(m,x)) \notag \\
&\le  C\,  \frac{\lan (\widetilde{J}_{n})\, \nu(J_{n-1})}{\lan({ J}_{n})\, \lan({\widetilde{\cal J}}_{n}
\cap J_{n-1}) }\ \lan(P(m,x))\le C
\,  \frac{\lan (\widetilde{J}_{n})} {\lan({ J}_{n})  } \frac{ \nu(J_{n-1})} {\lan({ J}_{n-1}) }\ \lan(P(m,x)) \notag \\
 &\le C\, \frac 1{\lan(P(k_{n},x_0))}\, \frac{\nu(J_{n-1})}{\lan(J_{n-1})}\,\lan (P(m,x)) \,. \notag
\end{align}

In any case, by taking $N_{n+1}$ large enough (and therefore $k_{n+1}$ also large enough) or $N_n$ large enough (and therefore $k_{n}$ also large enough) we have that for $J_{n+1} \subset P(m,x) \subseteq J_n$,
\begin{equation} \label{doscasos}
\nu(P(m,x)) \le C\, \frac 1{\lan(P(k_{j},x_0))^{1+\ep}} \,\lan(P(m,x))
\end{equation}
with $j=n+1$ in the case (1) and $j=n$ in the case (2).

Recall also that by the property (C) of expanding maps we have
$$
\sup_{P\in {\cal P}_n} \hbox{diam}(P) \le C \, \frac 1{\beta^n}
$$
Now given a  ball $B=B(x,r)$ with center $x\in {\cal C}_\ep$ we define the natural number $m=m(B)$ given by
\begin{equation} \label{defdem}
\frac {2C}{\beta^m} \le \hbox{diam}(B) < \frac {2C}{\beta^{m-1}}
\end{equation}
and the family ${\cal P}(B)$ as the collection of blocks in $P\in{\cal P}_m$ such that
$P\cap B\neq \emptyset$. Let us also denote by $n=n(B)$ the natural number such that
$d_n+k_n\le m <d_{n+1}+k_{n+1}$.  It is clear that
$$
\nu(B) \le \sum_{P\in {\cal P}(B)} \nu(P)
$$
and using (\ref{doscasos}) we obtain that
$$
\nu(B) \le C  \frac 1{\lan(P(k_{j},x_0))^{1+\ep}} \sum_{P\in {\cal P}(B)} \,\lan(P)
$$
where $j=n$ if  $d_n+k_n\le m <d_{n+1}$  and $j=n+1$ if $d_{n+1}\le m < d_{n+1}+k_{n+1}$. Notice that, by (\ref{defdem}), it is clear that $\cup_{P\in {\cal P}(B)} P \subset 2B:=B(x,2r)$ and since $\lan$ is a doubling measure we have that
$$
\sum_{P\in {\cal P}(B)} \,\lan(P) \le C\lan (B)\,.
$$
Therefore, in each of the above cases, we have
\begin{equation} \label{aquel}
\nu(B) \le C\,  \frac 1{\lan(P(k_{j},x_0))^{1+\ep}} \,\lan(B)\,.
\end{equation}
But, by the property (c2) of the Cantor set ${\cal C}_\ep$
$$
\lan(P(k_{j},x_0) \ge  e^{-N_j[(1+\ep)(\dimsup+\ep)(\ellsup+\ep)+\ep]}\,
$$
and by (\ref{defdem}) we obtain that
$$
\lan (B) \le C\, \hbox{diam}(B)^s \le C\, \frac 1{\beta^{(m-1)s}} \le C\, e^{-N_j s \log\beta}
$$
where we have used that $d_{n+1}\asymp N_{n+1}$ and $d_n+k_n\asymp d_n \asymp N_n$. Hence,
\begin{equation} \label{este}
\lan(P(k_{j},x_0) \ge \lan(B)^{[(1+\ep)(\dimsup+\ep)(\ellsup+\ep)+\ep]/s\log\beta}\,.
\end{equation}
From (\ref{aquel}) and  (\ref{este}) we get (\ref{frostman}).

\end{proof}

\begin{remark} \label{Markovsi}
\rm If $X\subset \R$ and $\lan$ is Lebesgue measure, then Theorem \ref{dimension} holds also for all approximable point $x_0\in X_0^+$. In fact, in this case
we can not assure that $\hbox{closure}(P(k_1,x_0))\subset P(0,x_0)$. However it is true that
$\hbox{closure}(P(k_1,x_0))\subset P(0,x_0) \cup \{x_0\}$
and from this fact we can conclude easily that
$$
{\cal C}_\ep := \bigcap_{n=0}^\infty \bigcup_{J\in {\cal J}_n} J
\subset \bigcap_{n=0}^\infty \bigcup_{J\in {\cal J}_n} \hbox{closure\,}(J) \subset {\cal C}_\ep \cup
S\,,
$$
where $S$ is a countable set. Hence,  ${\cal C}_\ep$ and $ \bigcap_{n=0}^\infty \bigcup_{J\in {\cal J}_n} \hbox{closure\,}(J)$ have the same Hausdorff dimensions. Also, since $\lan(J)=\lan(\hbox{closure\,}(J))$ the proof of Theorem \ref{dimension} allows to estimate the Hausdorff dimensions of $ \bigcap_{n=0}^\infty \bigcup_{J\in {\cal J}_n} \hbox{closure\,}(J)$.

\end{remark}

The next result  follows  from the proof of Theorem \ref{dimension}. In this case the sequence of radii is constant and therefore we are estimating the set of points returning periodically to a neighbourhood of the given point $x_0$. The proof is much more simple because the constructed Cantor-like sets have a more regular pattern.

\medskip

\begin{corollary} \label{recurrentes}
Let $(X,d,{\cal A},\lambda,T)$ be an expanding system with finite
entropy $H_\mu({\cal P}_0)$ with respect to the partition ${\cal
P}_0$ where $\mu$ is the ACIPM associated to the system. Let us
consider the grid $\Pibf=\{{\cal P}_n\}$. Let $r>0$ and let $P$ be
a block of ${\cal P}_{N_0}$. Then, given $\ep>0$, for all point $x_0\in X_0$,
there exist $k$ depending on $x_0$ and $r$, and $\widetilde N$
depending on $P$, $x_0$, $r$ and $\ep$ such that  for all
$N\geq\widetilde N$ the grid Hausdorff dimensions of the set
${\cal R}(P,x_0,r,N)$ of points $x\in P\cap X_0$ such that
$$
d(T^{d_j} (x),x_0) < r \ \mbox{ for } \ d_j=N_0 + k+(j-1)(N+k)  \ \mbox{ for } j=1,2,\ldots
$$
verify
$$
\Dim_{\Pibf,\lan} ({\cal R}(P,x_0,r,N)) = \Dim_{\Pibf,\mu} ({\cal
R}(P,x_0,r,N)) \geq 1-\ep-C_1/N \,.
$$
where $C_1$ is an absolute constant. Moreover for all $x_0\in X_0$ we have
\begin{enumerate}
\item[\rm 1.]  If the grid $\Pibf$ is $\lan$-regular then,
$$
\Dim_\lan  ({\cal R}(U,x_0,r,N))=\Dim_\mu  ({\cal R}(U,x_0,r,N )) \geq 1-\ep-C_2/N .
$$
\item[\rm 2.] If $\lan$ is a doubling measure verifying that $\lan (B(x,r))\le C\,r^s$ for all ball $B(x,r)$, then
$$
\Dim_{\lan} ({\cal R}(P,x_0,r,N)) = \Dim_{\mu} ({\cal R}(P,x_0,r,N)) \ge 1-\frac {\log C_3}{(N+k)s\log\beta}\,,
$$
with $C_3\asymp  1/\lan(P(k,x_0))$.
\end{enumerate}
\end{corollary}

\begin{proof}
We have now $\ellsup=0$ and we can do the same construction that in Theorem \ref{dimension} with $N_j=N$ and $k_j=k$ for all $j\ge 1$. The result for $\Dim_{\Pibf,\lan}$ follows from Corollary
\ref{azulesyrojos2} by taking $\alpha_j=e^{-Na}$, $\beta_j=e^{-Nb}$  with $a=h_\mu+2\ep$, $b=h_\mu-2\ep$  and $\gamma_j$  a constant. Part 1 is a consequence of Proposition
\ref{haches}. The proof of Part 2 is similar to the corresponding one in the proof of Theorem \ref{dimension}. Now instead of (\ref{doscasos}) we have that $\nu(P(m,x))\le C\,C_3^n \lan(P(m,x))$.
\end{proof}

\medskip

\begin{lemma}
Let $\{A_k\}$ be a decreasing sequence of Borel sets in $X$ such that $\Dim_{\Pibf,\lan} (A_k)\ge \beta>0$. Then,
$Dim_{\Pibf,\lan} (\cap_k A_k)\ge \beta$.
\end{lemma}

\begin{proof}
If $0<\alpha<\beta$, then ${\cal H}_{\Pibf,\lan}^\alpha (A_k)=\infty$ for all $k$ and therefore
${\cal H}_{\Pibf,\lan}^\alpha (\cap_k A_k)=\lim_{k\to\infty}{\cal H}_{\Pibf,\lan}^\alpha (A_k)=\infty $. It follows that
$\Dim_\lan (\cap_k A_k)\ge\alpha$. The result follows by letting $\alpha\to\beta$.
\end{proof}

\begin{remark}
\rm The lemma also holds (with the same proof) for the $\lan$-Hausdorff dimension.
\end{remark}

\begin{corollary} \label{cotasparadimconradios2}
Under the hypotheses of Theorem {\rm \ref{dimension}} we have that
$$
 \Dim_{\Pibf,\lan} \left\{ x\in U: \ \liminf_{n\to\infty} \frac{d(T^n(x),x_0)}{r_n} = 0 \right\} \ge
 \frac {h_\mu}{h_\mu + \dimsup \,\ellsup}\,.
$$
Moreover, if the grid $\Pibf$ is $\lan$-regular, then
$$
 \Dim_{\lan} \left\{ x\in U: \ \liminf_{n\to\infty} \frac{d(T^n(x),x_0)}{r_n} = 0 \right\} \ge
\frac {h_{\mu}}{h_\mu + \dimsup \, \ellsup} \, \Big(
1-\frac{\overline{\tau}_\lan(x_0)\dimsup\ell^2}{h_{\mu}^2\log\beta}\Big)\,.
$$
\end{corollary}

\begin{proof}
Notice that if $x\in U$ verifies
$$
d(T^n(x),x_0) \le r_n \,, \ \hbox{for
infinitely many $n$} \quad \implies \quad \liminf_{n\to\infty} \frac {d(T^n(x),x_0)}{r_n }\le 1
$$
and from Theorem \ref{dimension} we obtain that
$$
\Dim_{\Pibf,\lan} \left\{ x\in U: \ \liminf_{n\to\infty} \frac{d(T^n(x),x_0)}{r_n} \le 1 \right\} \ge
\frac {h_\mu}{h_\mu + \dimsup \,\ellsup}\,.
$$
By applying this last result to the sequence $\{r_n/m\}_{n=1}^\infty$ for any $m\in\N$, we get that
$$
\Dim_{\Pibf,\lan} \left\{ x\in U: \ \liminf_{n\to\infty} \frac{d(T^n(x),x_0)}{r_n} \le \frac 1m \right\} \ge
\frac {h_\mu}{h_\mu + \dimsup \,\ellsup}\,,
$$
and since
$$
\left\{ x\in X : \  \liminf_{n\to\infty} \frac {d(T^n(x),x_0)}{r_n} = 0 \right\}
= \bigcap_{m=1}^\infty \left\{ x\in X : \  \liminf_{n\to\infty} \frac {d(T^n(x),x_0)}{r_n} \le \frac 1m  \right\}\,.
$$
the lower bound in the statement follows from the above lemma. The proof of the second statement is similar.
\end{proof}

\

As in the measure section we are also interested in the size of
the set
$$
\widetilde{\cal W}(x_0,\{t_n\}) = \{x\in X_0: \ T^k(x) \in
P(t_k,x_0) \ \hbox{for infinitely many $k$}\}
$$
with $\{t_k\}$ an increasing sequence of positive integers and
$x_0\in X_0^+$. We recall that if
$x_0=[\;i_0\;i_1\ldots\;]$, then $\widetilde{\cal W}(x_0,\{t_n\})$
is the set of points $x=[\;m_0\;m_1\ldots\;]\in X_0$ such that
$$
m_k=i_0, m_{k+1}=i_i,\;\ldots \;,m_{k+t_k}=i_{t_k}
$$
for infinitely many $k$.

\

\begin{theorem} \label{dimensioncode} Let $(X,d,{\cal A},\lambda,T)$ be an expanding system
with finite entropy $H_\mu({\cal P}_0)$ with respect to the
partition ${\cal P}_0$ where $\mu$ is the ACIPM associated to the system.
Let $\{t_n\}$ be a non decreasing
sequence of positive integers and let $U$ be an open set in $X$
with $\mu(U)>0$. Let us consider the grid $\Pibf=\{{\cal P}_n\}$.
Then, for all approximable point $x_0\in X_0$, the grid Hausdorff
dimensions of the set
$$
\widetilde{\cal W}(U,x_0,\{t_n\})=\{x\in U\cap X_0: \ T^k(x) \in
P(t_k,x_0) \ \hbox{for infinitely many $k$}\}\,,
$$
verify
$$
\Dim_{\Pibf,\lan} ({\cal \widetilde W}(U,x_0,\{t_n\})) =
\Dim_{\Pibf,\mu} ({\cal \widetilde W}(U,x_0,\{t_n\})) \ge \frac
{h_\mu}{h_\mu +\overline{L}(x_0)}  \,,
$$
where $\overline{L}(x_0)= \limsup_{n\to\infty} \frac 1n
\log\frac{1}{\lan(P( t_n,x_0))}$ and $h_\mu$ is the entropy of $T$
with respect to $\mu$.

Moreover, for all approximable point $x_0\in X_0$, the Hausdorff dimension of the set
${\cal \widetilde W}(U,x_0,\{t_n\})$ verify:

\begin{enumerate}
\item[\rm 1.]  If the grid $\Pibf$ is
$\lan$-regular, then
$$
\Dim_{\lan} ({\cal \widetilde W}(U,x_0,\{t_n\})) = \Dim_{\mu}
({\cal \widetilde W}(U,x_0,\{t_n\})) \ge \frac {h_\mu}{h_\mu
+\overline{L}(x_0)} \, \Big( 1- \frac{\overtau \,\overline{w}\,\overline{L}(x_0) }{h_{\mu}^2} \Big)\,,
$$
where $\overline{w}= \limsup_{n\to\infty} \frac {t_n}n$\,.

\item[\rm 2.] If $\lan$ is a doubling measure verifying that $\lan(B(x,r))\le C\,r^s$ for all ball $B(x,r)$, then
$$
\Dim_\lan  (\widetilde{\cal W}(U,x_0,\{t_n\}))=\Dim_\mu  ({\cal W}(U,x_0,\{t_n\}))\ge 1 - \frac {\overline{L}(x_0) \overline{w}}{s\log\beta}\,.
$$
\end{enumerate}

\end{theorem}

\begin{remark} \rm  We recall that from Remark
\ref{approx} and Lemma \ref{justificacion} we know that the
set of approximable points such that $\overtau=0$ has full
$\lan$-measure.
\end{remark}

As in the case of radii, we have the following consequence of the
proof of Theorem \ref{dimensioncode} when we take the sequence $\{
t_n:=t\}$ constant. We are estimating the set of points in whose
code appear periodically the first $t$ digits of the code of the
point $x_0$. The proof is similar.

\begin{corollary}
Let $(X,d,{\cal A},\lambda,T)$ be an expanding system with finite
entropy $H_\mu({\cal P}_0)$ with respect to the partition ${\cal
P}_0$ where $\mu$ is the ACIPM associated to the system. Let us
consider the grid $\Pibf=\{{\cal P}_n\}$. Let $t\in\N$ and let $P$
be an block of ${\cal P}_{N_0}$. Then, given $\ep>0$ for all point $x_0\in
X_0$ there exist $k$ depending on $x_0$ and $t$, and $\widetilde
N$ depending on $P$, $x_0$ $\ep$ and $t$, such that  for all
$N\geq\widetilde N$ the grid Hausdorff dimensions of the set
$\widetilde{{\cal R}}(P,x_0,r,N)$ of points
$x=[\;m_0\;m_1\ldots\;]\in P\cap X_0$ such that for $j=1,2,\ldots$
$$
m_{d_j}=i_0,\,  m_{d_j+1}=i_1,\;\ldots \;,m_{d_j+t}=i_{t}
\ \mbox{ with } \ d_j=N_0+k +(j-1)(N+k)
$$
verify
$$
\Dim_{\Pibf,\lan} ({\cal R}(P,x_0,r,N)) = \Dim_{\Pibf,\mu} ({\cal R}(P,x_0,r,N)) \geq 1-\ep-C_1/N \,,
$$
where $C_1$ is an absolute constant. Moreover for all $x_0\in X_0$ we have
\begin{enumerate}
\item[\rm 1.]  If the grid $\Pibf$ is $\lan$-regular then,
$$
\Dim_\lan  ({\cal R}(U,x_0,r,N))=\Dim_\mu  ({\cal R}(U,x_0,r,N )) \geq 1-\ep-C_2/N .
$$
\item[\rm 2.] If $\lan$ is a doubling measure verifying that $\lan (B(x,r))\le C\,r^s$ for all ball $B(x,r)$, then
$$
\Dim_{\lan} ({\cal R}(P,x_0,r,N)) = \Dim_{\mu} ({\cal R}(P,x_0,r,N)) \ge 1-\frac {\log C_3}{(N+k)s\log\beta}\,,
$$
with $C_3\asymp  1/\lan(P(k,x_0))$.
\end{enumerate}
\end{corollary}

\begin{proof}[{\it Proof of Theorem {\rm \ref{dimensioncode}}.}]
We may assume that $\overline{L}(x_0)$, $\overtau$ and
$\overline{w}$ are all finite, since otherwise our Hausdorff
dimension estimates are obvious. Now, the proof is similar to the
proof of Theorem \ref{dimension}. For each $\ep>0$ we construct a
Cantor-like  set ${\cal C}_{\ep}\subset \widetilde{\cal
W}(U,x_0,\{t_n\})$. Recall that in the proof of Theorem
\ref{dimension} we defined an increasing sequence $\cal D$ of
allowed indexes. Here, we define $\cal D$ in the following way:
Let ${\cal I}(x_0)=\{p_i\}$ denote the sequence associated to $x_0$ given
by Definition \ref{aproximable}. For each $p_k \in {\cal I}(x_0)$ let
$n(k)$ denote the greatest natural number such that $t_{n(k)}\le
p_k$. We denote by $\cal D$ the set of these allowed indexes. We
will write ${\cal D}=\{d_i\}$ with $d_i<d_{i+1}$\,.

With this new definition of ${\cal D}$ we have that if $d\in\cal
D$, then there exists $k(d)\in {\cal I}(x_0)$ such that
$$
t_d \le k(d) <t_{d+1} \qquad \text{and} \qquad P(k(d),x_0)\subset
P(t_d,x_0)\,.
$$
These two properties substitute to (\ref{permitidos2}) and
(\ref{permitidos1}). For all $d$ large enough we have that
$$
\frac{k(d)}{d}\leq (1+\ep) \overline{w}\,.
$$
This inequality substitute to (\ref{kvsN}).
With the above considerations and proceeding as in the proof of Theorem \ref{dimension}
we construct  the families $\widetilde{\cal J}_j, \, {\cal J}_j$ and the numbers $N_j$ and
$k_j\in {\cal I}(x_0)$  with the properties (b), (c1), (c3),  (d) and (e).
The corresponding properties (a) and (c2) are now the following ones:
\begin{itemize}
\item[(a)]  For all point $x$ in $J_j\in{\cal J}_j$
$$
T^{d_j}(x)\in P(t_{d_j},x_0)\,.
$$
\item[(c2)] For all $J_j\in{\cal J}_j$ there exist a unique $\widetilde{J}_{j}\in\widetilde{\cal J}_{j}$ so
that $\text{closure\,}(J_j)\subset \widetilde{J}_{j}$ and
$$
\frac {\lan(J_j)}{\lan(\widetilde{J}_{j})} \asymp \lan(P(k_j,x_0)) \qquad \hbox{and} \qquad
\frac {\lan(J_j)}{\lan(\widetilde{J}_{j})}\ge
e^{-N_j[(1+\ep)(\overline{L}(x_0)+\ep)+\ep]}
\,.
$$
\end{itemize}
The rest of the proof is similar.
\end{proof}

\begin{remark} \label{Markovsicode} \rm Using the same argument that in Remark \ref{Markovsi} we get that
if $X\subset \R$ and $\lan$ is Lebesgue measure, then Theorem \ref{dimensioncode}
holds also for all approximable point $x_0\in X_0^+$.
\end{remark}

For points $x_0\in X_1$ we have the following version of the above theorem.

\begin{theorem} \label{dimensioncode2} Let $(X,d,{\cal A},\lambda,T)$
be an expanding system with finite entropy $H_\mu({\cal P}_0)$
with respect to the partition ${\cal P}_0$, where $\mu$ is the
ACIPM associated to the system. Let
$\{t_n\}$ be a non decreasing sequence of positive integers and
let $U$ be an open set in $X$ with $\mu(U)>0$. Let us consider the
grid $\Pibf=\{{\cal P}_n\}$. Then, for all approximable point
$x_0\in X_1$, and therefore for $\lan$-almost point $x_0\in X$,
the grid Hausdorff dimensions of the set
$$
\widetilde{\cal W}(U,x_0,\{t_n\})=\{x\in U\cap X_0: \ T^k(x) \in
P(t_k,x_0) \ \hbox{for infinitely many $k$}\}\,,
$$
verify
$$
\Dim_{\Pibf,\lan} ({\cal \widetilde W}(U,x_0,\{t_n\})) =
\Dim_{\Pibf,\mu} ({\cal \widetilde W}(U,x_0,\{t_n\})) \ge \frac
{1}{1 +\overline{w}} \,,
$$
where $\overline{w}= \limsup_{n\to\infty} \frac {t_n}n$\,.
Moreover, if the grid $\Pibf$ is $\lan$-regular, then also
$$
\Dim_{\lan} ({\cal \widetilde W}(U,x_0,\{t_n\})) = \Dim_{\mu}
({\cal \widetilde W}(U,x_0,\{t_n\})) \ge \frac {1}{1
+\overline{w}}\,.
$$
\end{theorem}

\begin{proof} We may assume that $\overline{w}<\infty$ since
otherwise our estimations are trivial. The proof is similar to the
proof of Theorem \ref{dimensioncode} but using that  for $d$ large
$$
\lan(P(t_d,x_0)) \ge e^{-t_d(h_\mu+\ep)} \ge
e^{-d(\overline{w}+\ep)(h_\mu+\ep)}\,.
$$
\end{proof}

\subsection{Upper bounds of the dimension}\label{upperbounds}

We will prove now some upper bounds for the $\lan$-grid Hausdorff
dimension of ${\cal W}(U,x_0,\{r_n\})$ and $\widetilde{\cal W}(U,x_0,\{t_n\})$ in the case that the partition ${\cal P}_0$ is finite.

\begin{proposition} \label{dimensionporarribacode} Let $(X,d,{\cal A},\lambda,T)$ be an expanding system such that the partition ${\cal P}_0$ is finite. Let $\mu$ be the ACIPM associated to the system. Let $\{t_n\}$ be a non decreasing
sequence of positive integers and $U$ be an open set in $X$ with
$\mu(U)>0$.

Then, if $x_0 \in X_0^+$, we have that
the grid Hausdorff dimensions of the set
$$
\widetilde{\cal W}(U,x_0,\{t_n\})=\{x\in U\cap X_0: \ T^k(x) \in
P(t_k,x_0) \ \hbox{for infinitely many $k$}\}\,,
$$
verify
$$
\Dim_{\Pibf,\lan} ({\widetilde{\cal W}}(U,x_0,\{t_n\})) =
\Dim_{\Pibf,\mu} ({\widetilde{\cal W}}(U,x_0,\{t_n\}))  \le \min \left\{1,  \frac {\log D}{h_\mu
+ \underline{L}(x_0)} \right\}\,.
$$
where $\underline{L}(x_0)= \liminf_{n\to\infty} \frac
1n\log\frac{1}{\lan(P( t_n,x_0))}$,  $h_\mu=h_{\mu}(T)$  is the entropy of $T$ with respect to the measure  $\mu$ and $D$ is the cardinality of ${\cal P}_0$. Moreover, if $x_0\in X_1$, then
$$
\Dim_{\Pibf,\lan} ({\widetilde{\cal W}}(U,x_0,\{t_n\})) =
\Dim_{\Pibf,\mu} ({\widetilde{\cal W}}(U,x_0,\{t_n\}))  \le \min \left\{1, \frac {\log D}{(1+\underline w)h_\mu} \right\}\,,
$$
where $\underline{w}=\liminf_{n\to\infty} \frac 1n \, t_n \,.$
\end{proposition}

\begin{proof} We define the collections
$$
{\cal F}_n=\{T^{-n}(P(t_n,x_0))\cap Q:  \ Q\in {\cal P}_n\} \,.
$$
Then the collection
$$
{\cal G}_N = \bigcup_{n=N}^\infty {\cal F}_n \,.
$$
covers the set ${\widetilde{\cal W}}(U,x_0,\{t_n\})$ for all $N\in\N$.
Using Proposition \ref{yaesta} we get that
\begin{equation} \label{todos}
\sum_{k=N}^\infty \sum_{F\in {\cal F}_k} \mu(F)^\tau  \le C \sum_{k=N}^\infty \mu(P(t_k,x_0))^\tau \sum_{Q\in {\cal P}_k} \mu(Q)^\tau \,.
\end{equation}
Let us consider now the following two subcollections of ${\cal P}_k$:
$$
{\cal P}_{k, {\rm small}} = \{Q \in {\cal P}_k : \ \mu(Q)\le e^{-kh_\mu} \} \,, \qquad
{\cal P}_{k,{\rm big}} = \{Q \in {\cal P}_k : \ \mu(Q) > e^{-kh_\mu} \}\,.
$$
Then,
$$
\sum_{Q\in {\cal P}_{k,{\rm small}}} \mu(Q)^\tau \le D^k e^{-k\tau h_\mu}
$$
and
$$
\sum_{Q\in {\cal P}_{k,{\rm big}}} \mu(Q)^\tau = \sum_{Q\in {\cal P}_{k,{\rm big}}}
\frac 1{\mu(Q)^{1-\tau}} \,\mu(Q) \le e^{kh_\mu(1-\tau)}\,.
$$
Since $h_\mu \le H_\mu({\cal P}_0)\le \log D$ we have that
$$
e^{kh_\mu(1-\tau)} \le D^k e^{-k\tau h_\mu}
$$
and therefore using (\ref{todos}) we obtain that
$$
\sum_{k=N}^\infty \sum_{F\in {\cal F}_k} \mu(F)^\tau  \le 2C \sum_{k=N}^\infty D^k e^{-k\tau h_\mu} \mu(P(t_k,x_0))^\tau
$$
By part (iii) of Theorem E we know that $\mu(P(t_k,x_0))\asymp
\lan(P(t_k,x_0))$. Then for $t_k$ large enough
$$
\mu(P(t_k,x_0))\asymp \lambda(P(t_k,x_0))\leq e^{-k(\underline{L}-\ep)} \,.
$$
Hence
$$
\sum_{G\in {\cal G}_N}  \mu(G)^\tau \to 0 \qquad \text{ when} \quad
N\to\infty
$$
for
$$
\tau > \frac {\log D}{h_\mu + \underline{L}(x_0)-\ep}\,.
$$
For these $\tau$'s the $\tau$-dimensional $\mu$-grid Hausdorff measure
of  $\widetilde{\cal W}(U,x_0,\{t_n\})$ is zero and therefore
$$
\Dim_{\Pibf,\mu} ({\cal \widetilde W}(U,x_0,\{t_n\})) \le \frac
{\log D}{h_\mu + \underline{L}(x_0)-\ep}\,.
$$
The result follows by taking $\ep$ tending to zero and using Lemma \ref{dimensioncoincide}. Finally, if $x_0\in X_1$, we have that $\underline{L}(x_0)=h_\mu \, \underline{w}$.
\end{proof}

%

\begin{remark}
\rm Even in the case that the partition ${\cal P}_0$ is not finite a slight modification of the above proof shows that
\begin{equation} \label{intX1}
\Dim_{\Pibf,\lan} (\widetilde{\cal W}(U, x_0,\{t_n\})\cap X_1) \le \frac {h_\mu}{h_\mu+\underline{L}(x_0)}\,.
\end{equation}
To see this, notice that if we define for any $\ep>0$ the subcollections:
$$
{\cal P}_{n, \ep, {\rm big}} = \{Q \in {\cal P}_n : \ \mu(Q) > e^{-n(h_\mu+\ep)} \}\,.
$$
then the set $\widetilde{\cal W}(U, x_0,\{t_n\})\cap X_1$ can be covered by the collections
$$
{\cal G}_N^\ep = \bigcup_{n=N}^\infty {\cal F}_{n,{\rm big}}^\ep \,.
$$
where
$$
{\cal F}_{n.{\rm big}}^\ep =\{T^{-n}(P(t_n,x_0))\cap Q:  \ Q\in {\cal P}_{n,\ep,{\rm big}}\} \,.
$$
The proof of (\ref{intX1}) follows now easily.
\end{remark}

\begin{proposition} \label{dimensionporarriba} Let $(X,d,{\cal A},\lambda,T)$
be an expanding system with $X\subset \R$ and such that the partition
${\cal P}_0$ is finite. Let $\mu$ be the ACIPM associated to the system. Let $\{r_n\}$ be a non increasing sequence of positive
numbers and $U$ be an open set in $X$ with $\mu(U)>0$.

Then, for $x_0\in X_0$, the Hausdorff
dimensions of the set
$$
{\cal W}(U,x_0,\{r_n\})=\{ x\in U\cap X_0: \  d(T^n (x),x_0)\le r_n \ \mbox{for
infinitely many} \ n \}
$$
verify
\begin{equation} \label{dimporarriba}
\Dim_{\lan} ({\cal W}(U,x_0,\{r_n\}))=\Dim_{\mu} ({\cal W}(U,x_0,\{r_n\})) \le \min \left\{ 1, \frac {\log D}{h_\mu +
\diminf \ellinf} \right\}\,,
\end{equation}
where $D$ is the cardinality of ${\cal P}_0$, $h_\mu=h_{\mu}(T)$ is the entropy of $T$ with respect to  $\mu$ and
$$
\ellinf= \liminf_{n\to\infty} \frac 1n\log\frac{1}{r_n}
$$
\end{proposition}

\begin{proof}  We define the collections
$$
{\cal F}_n=\{T^{-n}(B(x_0,r_n))\cap Q:  \ Q\in {\cal P}_n\} \,.
$$
Notice that for $n$ large enough $B(x_0,r_n)\subseteq P(0,x_0)$ and then
$T^{-n}(B(x_0,r_n))\cap Q$ is an interval for any  $Q\in {\cal P}_n$.
Therefore, the collection of intervals
$$
{\cal G}_N = \bigcup_{n=N}^\infty {\cal F}_n
$$
covers the set ${\cal W}(U,x_0,\{r_n\})$ for all $N\in\N$ large enough. Using Proposition \ref{yaesta} we get that
$$
\sum_{k=N}^\infty \sum_{F\in {\cal F}_k} \mu(F)^\tau  \le C \sum_{k=N}^\infty \mu(B(x_0,r_k))^\tau \sum_{Q\in {\cal P}_k} \mu(Q)^\tau \,.
$$
By estimating $\sum_{Q\in {\cal P}_k} \mu(Q)^\tau$ as in
the proof of Proposition \ref{dimensionporarribacode} we get
$$
\sum_{k=N}^\infty \sum_{F\in {\cal F}_k} \mu(F)^\tau  \le 2C \sum_{k=N}^\infty D^k e^{-k\tau h_\mu} \mu(B(x_0,r_k))^\tau
$$
For $r_k$ small enough we have that $B(x_0,r_k)\subset P(0,x_0)$ and then
$\mu(B(x_0,r_k))\asymp\lan(B(x_0,r_k))$ by part (iii) of Theorem E. Hence, from
the definition of $\diminf$ and Lemma \ref{egorof}, we conclude
that given $\ep>0$,
$$
\sum_{k=N}^\infty \sum_{F\in {\cal F}_k} \mu(F)^\tau  \le C \sum_{k=N}^\infty r_k^{(\diminf-\ep)\tau}  D^k e^{-k\tau h_\mu}
\longrightarrow 0
$$
as $N\to\infty$, if
$$
\tau > \frac {\log D}{h_\mu + (\diminf-\ep)
(\ellinf-\ep)}\,.
$$
For these $\tau$'s the $\tau$-dimensional $\mu$-Hausdorff measure
of ${\cal W}(U,x_0,\{r_n\})$ is zero  and therefore
$$
\Dim_{\mu} ({\cal W}(U,x_0,\{r_n\})) \le  \frac {\log D}{h_\mu +
(\diminf-\ep) (\ellinf-\ep)} \,.
$$
The result follows by taking $\ep$ tending to zero and using Lemma \ref{dimensioncoincide}.
\end{proof}

\section{Applications}\label{aplicaciones}

\subsection{Markov transformations}

Let $\lan$ be Lebesgue measure in $[0,1]$. A map
$f:[0,1]\longrightarrow [0,1]$ is a {\it Markov transformation} if
there exists a  family ${\cal P}_0=\{P_j\}$ of disjoint
open intervals in $[0,1]$ such that

\begin{itemize}
\item[(a)] $\lan([0,1]\setminus \cup_j P_j)=0$.

\item[(b)] For each $j$, there exists a set $K$ of indices such
that $f(P_j)=\cup_{k\in K} P_k$ (mod $0$).

\item[(c)] $f$ is derivable in $\cup_j P_j$ and there exists
$\sigma>0$ such that $|f'(x)|\ge \sigma$ for all $x\in\cup_j P_j$.

\item[(d)] There exists $\gamma>1$ and a non zero natural number
$n_0$ such that if $f^m(x)\in \cup_j P_j$ for all $0\le m \le
n_0-1$, then $|(f^{n_0})'(x)|\ge\gamma$.

\item[(e)] There exists a non zero natural number $m$ such that
$\lan(f^{-m}(P_i)\cap P_j)>0$ for all $i,j$.

\item[(f)] There exist constants $C>0$ and $0<\alpha\le 1$ such that, for all $x,y \in P_j$,
$$
\Big| \frac{f'(x)}{f'(y)}-1 \Big| \le C|f(x)-f(y)|^\alpha \,.
$$
\end{itemize}

Markov transformations are expanding maps with parameters $\alpha$
and $\beta=\gamma^{1/n_0}$, see \cite{M}, p.171, and therefore, by
Theorem E, there exists a unique $f$-invariant probability measure
$\mu$ in $[0,1]$ which is absolutely continuous with respect to
Lebesgue measure and satisfies properties (i)-(v) in Theorem E. As
a consequence of our results we obtain

\begin{theorem} \label{markovradios}
Let $f:[0,1]\longrightarrow [0,1]$ be a {\it Markov
transformation} and $\{r_n\}$ be a non increasing sequence of
positive numbers. Then,

\begin{itemize}

\item[\rm (1)] If $\sum_n r_n^{1+\ep}=\infty$ for some $\ep>0$,
then for almost all $x_0\in [0,1]$ we have that
$$
\liminf_{n\to\infty} \frac {|f^n(x)-x_0|}{r_n} = 0\,, \qquad
\text{for almost all $x\in [0,1]$.}
$$

\item[\rm (2)] If $\sum_n r_n <\infty$, then for all $x_0\in
\cup_j P_j$ we have that
$$
\liminf_{n\to\infty} \frac {|f^n(x)-x_0|}{r_n} = \infty \,, \qquad
\text{for almost all $x\in [0,1]$.}
$$

\item[\rm (3)] If $H_\mu ({\cal P}_0)<\infty$, then, for almost all $x_0\in [0,1]$, we have that
$$
 \hbox{\rm Dim}\left\{x\in [0,1]: \
\liminf_{n\to\infty} \frac {|f^n(x)-x_0|}{r_n} = 0 \right\}\ge \frac{h_\mu}{h_\mu +\ellsup}\,,
$$
where  $\ellsup=\limsup
\frac 1n \log \frac 1{r_n}$, $h_\mu= \int_0^1 \log |f'(x)|\, d\mu(x)$
and $\Dim$ denotes Hausdorff dimension\,.
\end{itemize}
If the partition ${\cal P}_0$ is finite, then all the statements hold for all $x_0\in [0,1]$.

\end{theorem}

Let us observe that if $\sum_n r_n = \infty$ the theorem does not tell us what is the
measure of the set where
$$
\liminf_{n\to\infty} \frac {|f^n(x)-x_0|}{r_n} = 0
$$
but by part (3) we know that this set is big since has positive
Hausdorff dimension.

\begin{remark} \label{maspuntos} \rm
For sake of simplicity we have stated the above theorem for almost all point $x_0$.
However, the results in the previous sections give more information if we choose an specific $x_0$, see Remarks \ref{tau=0} and \ref{approx}-\ref{opeque–a}.
\end{remark}

\begin{proof}[Proof of Theorem {\rm \ref{markovradios}}]
Part (1) follows from Lemma \ref{justificacion} and
Corollary \ref{hartos2}. Part (2) is a consequence of Corollary
\ref{hartosinf}. Finally, part (3) follows from Remark \ref{approx}, Remark
\ref{regularreal}, Lemma \ref{justificacion} and Corollary
\ref{cotasparadimconradios2}. Finally, to get the result when the partition
${\cal P}_0$ is finite, we use additionally  that,  in this case, $X_0^+=[0,1]$, $\overtau=0$ for all $x_0\in [0,1]$ and Remarks \ref{ensoporte}, \ref{approxfinito} and \ref{Markovsi}.
\end{proof}

Recall that, as we saw in Section \ref{codes}, given an expanding map we have a code for almost all point $x_0$, and more precisely for all $x_0\in X_0^+$.
The following result summarizes our results about coding for Markov transformations.

\begin{theorem} \label{markovcode}
Let $f:[0,1]\longrightarrow [0,1]$ be a Markov transformation and
$\{t_n\}$ be a non decreasing sequence of natural numbers. Given a
point $x_0=[\,i_0,i_1,\dots \,]\in X_0^+$, let $\widetilde{\cal
W}(x_0,\{t_n\})$ be the set of points $x=[\,m_0,m_1,\dots\,]\in X_0$ such
that
$$
m_n=i_0\,, \; m_{n+1}=i_1\,,\; \dots\;, \; m_{n+t_n} = i_{t_n} \;, \qquad \text{for infinitely many $n$}.
$$
Then,

\begin{itemize}
\item[\rm (1)] If $\sum_n \lan(P(t_n,x_0))=\infty$, then $\lan(\widetilde{W}(x_0,\{t_n\})))=1$.
Moreover, if the partition ${\cal P}_0$ is finite or if $f(P)=[0,1]$ $($mod $0)$ for all $P\in {\cal P}_0$, then we have the following quantitative version:
$$
\lim_{n\to\infty} \frac {\#\{i\le n: \ f^i(x)\in
P(t_i,x_0)\}}{\sum_{j=1}^n \mu(P(t_j,x_0))} =1\,,  \qquad
\hbox{for $\lan$-almost every $x$.}
$$
\item[\rm (2)] If $\sum_n \lan(P(t_n,x_0))<\infty$, then $\lan(\widetilde{W}(x_0,\{t_n\})))=0$.
\item[\rm (3)] If $H_\mu ({\cal P}_0)<\infty$, then, for almost all $x_0\in X$, we have that
$$
\Dim (\widetilde{\cal W}(x_0,\{t_n\}) \ge \frac 1{1+\overline{w}}  \,.
$$
where $\overline{w}=\limsup_{n\to\infty} \frac {t_n}n$
and $\Dim$ denotes Hausdorff dimension.
\end{itemize}
\end{theorem}

\begin {remark} \rm
Even though part (3) is stated for almost every $x_0$ a more precise result
for an specific $x_0$ follows from Theorem \ref{dimensioncode} and Remark \ref{Markovsicode}. Recall also that any grid
contained in $\R$ is regular.
\end{remark}

\begin{proof}[Proof of Theorem {\rm \ref{markovcode}}]
Part (1) and (2) follow from Theorem \ref{elbuenocode} and
Proposition \ref{serieconverge}, respectively. Part (3) is
a consequence of Lemma \ref{justificacion}, Remark
\ref{approx}, Remark \ref{regularreal} and Theorem
\ref{dimensioncode2}.
\end{proof}

\subsubsection{Bernoulli shifts and subshifts of finite type} \label{markovchain}

Given a natural number $D$ let $\Sigma$ denote the space of all
infinite sequences $\{(i_0\,,\,i_1\,,\,\dots)\}$ with
$i_n\in\{0,1,\dots,D-1\}$ endowed with the product topology. The
left shift $\sigma: \Sigma \longrightarrow \Sigma$ is the
continuous map defined by
$$
\sigma (i_0\,,\,i_1\,,\,\dots) =(i_1\,,\,i_2\,,\,\dots)\,.
$$
For every positive numbers $p_0,p_1,\dots,p_{D-1}$ verifying
$\sum_{i=0}^{D-1} p_i=1$ we define the function
$$
\nu (C^{j_0,j_1,\dots,j_t}_{i_0,i_1,\dots,i_t})=p_{i_0} p_{i_1} \cdots p_{i_t} \,,
$$
where $C^{j_0,j_1,\dots,j_t}_{i_0,i_1,\dots,i_t}$ is the cylinder
$$
C^{j_0,j_1,\dots,j_t}_{i_0,i_1,\dots,i_t} = \{(k_0\,,k_1\,,\,\dots)\in \Sigma:
\ k_{j_s}=i_s \ \hbox{for all $s=0,1,\dots,t$}\}\,.
$$
It is well known that we can extend the set function $\nu$ to a
probability measure defined on the $\sigma$-algebra of the Borel
sets of $\Sigma$\,. The space $(\Sigma, \sigma, \nu)$ is called a
(one-sided) {\it Bernoulli shift}.

We can generalize the full shift space $(\Sigma, \sigma,
\nu)$ by considering the set $\Sigma_A$ defined by
$$
\Sigma_A=\{(i_0\,,\,i_1\,,\,\dots)\in\Sigma: \
a_{i_k,i_{k+1}}=1 \ \hbox{for all $k=0,1,\ldots$}\}\,,
$$
where $A=(a_{i,j})$ is a $D\times D$ matrix with entries
$a_{i,j}=0$ or $1$. The matrix $A$ is known as a {\it transition
matrix}. Let us consider now a new $D\times D$ matrix $M=(p_{i,j})$ such
that $p_{i,j}=0$ if $a_{i,j}=0$, and
\begin{eqnarray}
&(1)& \ \sum_{j=0}^{D-1} p_{i,j}=1\,, \quad \hbox{for every $i=0,1,\dots,D-1$}. \notag \\
&(2)& \ \sum_{i=0}^{D-1} p_i p_{i,j}=p_j\,, \quad \hbox{for every
$j=0,1,\dots,D-1$}. \notag
\end{eqnarray}
The numbers $p_{i,j}$ are called the {\it transition
probabilities} associated to the transition matrix $A$ and the
matrix $M$ is called a {\it stochastic matrix}. Observe that the
probability vector $(p_0,\dots,p_{D-1})$ is an eigenvector of the
matrix $M$.

We introduce now a probability measure $\nu$ on all Borel
subsets of $\Sigma_A$ by extending the set function defined by
$$
\nu(C^{n,\dots,n+t}_{i_0,i_1,\dots,i_t}) = p_{i_0} p_{i_0,i_1} \cdots
p_{i_{t-1},i_t}\,.
$$
The space $(\Sigma_A,\sigma,\nu)$ is called a (one-sided)
{\it subshift of finite type} or a (one-sided) {\it Markov chain}.

We will explain now how to associate to (one-sided) Bernoulli
shifts or (one-sided) subshift of finite type a Markov
transformation:

Let $(\Sigma, \sigma, \nu)$ be a (one-sided) Bernoulli shift and
let $\lan$ denote the Lebesgue measure in $[0,1]$. Consider a
partition $\{P_0, \dots, P_{D-1}\}$ of $[0,1]$ in $D$ consecutive
open intervals such that $\lan(P_j)=p_j$ for $j=0,1,\dots,D-1$. We
define now a function $f:[0,1]\longrightarrow [0,1]$ by letting
$f$ to be linear and bijective from each $I_j$ onto $(0,1)$, i.e.
$$
f(x)= \frac 1{p_j}\, \left( x- \sum_{k=0}^{j-1}
p_k\right) \,, \qquad \hbox{if $x\in P_j$}\,,
$$
and $f$ equal to zero on the boundaries of the intervals $P_j$. It is easy to check that $f$ is a Markov transformation and therefore an expanding map.

Define now a mapping $\pi: \Sigma \longrightarrow [0,1]$ by
$$
\pi((i_0,i_1,\dots))=\bigcap_{n=0}^\infty \hbox{closure\,} (f^{-n} (
{P}_{i_n}))\,.
$$
Then it is not difficult to see that $\pi$ is continuous and
$f\circ \pi= \pi\circ \sigma$.
$$
\begin{CD}
\Sigma @>\sigma>> \Sigma \\
@V\pi VV @VV \pi V \\
[0,1] @>f>> [0,1]
\end{CD}
$$
Notice that the space of codes associated through $f$ to the
points in $X_0$ (as we explained in Section \ref{codes}) is
precisely the set $\Sigma_0$ of all sequences $(i_0,i_1,\dots)$
such that there is not $k_0\in\N$ such that $i_j=0$ for all $j\ge
k_0$ or $i_j=D-1$ for all $j\ge k_0$, and that $\pi$ is bijective
from this set onto $X_0=[0,1]\setminus \cup_{n=0}^\infty
f^{-n}(\cup_i \partial P_i)$.

It is also easy to check that the image measure of the product
measure $\nu$ under $\pi$ is precisely the Lebesgue measure in
$[0,1]$ and that $f$ preserves Lebesgue measure. Therefore the
dynamical systems $(\Sigma,\sigma,\nu)$ and $([0,1],f,\lan)$
are isomorphic. We have also that the Hausdorff dimension of a
Borel subset $B\subset [0,1]$ coincides que the $\nu$-Hausdorff
dimension of $\pi^{-1}(B)$.

The subshifts of finite type whose stochastic matrix $M$ is transitive
(i.e. there exists $n_0>0$ such that all entries of $M^{n_0}$ are positive)
can be also thought as Markov transformations:
Let  $(\Sigma_A, \sigma, \nu)$ be a  subshift
of finite type with respect to the stochastic matrix $M=(p_{i,j})$
and the probability vector $(p_0,p_1,\dots,p_{D-1})$. Consider as
before a partition $\{P_0, P_1, \dots, P_{D-1}\}$ of $[0,1]$ in
$D$ consecutive open intervals such that $\lan(P_i)=p_i$ for
$i=0,1,\dots,D-1$. We divide now each interval $P_i$ into $D$
consecutive open intervals $P_{i,j}$ such that $\lan(P_{i,j})=p_i
p_{i,j}$ for $j=0,1,\dots, D-1$. If $p_{i,j}=0$ we take
$P_{i,j}=\emptyset$. Notice that by property (1) of stochastic
matrices we have that $\lan(P_i)=\sum_j \lan (P_{i,j})$.

We define now a function $f:[0,1]\longrightarrow [0,1]$ in the following way:
for each $P_{i,j}\ne \emptyset$,
$$
f(x) = \frac 1{p_{j,i}} \left( x - \sum_{k=0}^{j-1} p_k p_{k,i} - \sum_{\ell=0}^{i-1} p_\ell \right) +
\sum_{k=0}^{j-1} p_k \,, \qquad \hbox{if } x\in P_{i,j} \,.
$$
We define also $f$ on the points beloging to $P_i \cap \partial
P_{i,j}$ in such a way that $f$ is continuous in that points.
Finally we define $f$ to be zero on the boundaries of the
intervals $P_i$.

When the stochastic matrix $M$ verifies that there
exists $n_0$ such that all entries of $M^{n_0}$ are positive,
it is easy to check  that $f$ is a Markov transformation and
therefore an expanding map with respect to the partition ${\cal
P}_0 = \{P_{i}: \ i=0,1,\dots,D-1\}$. The condition on $M$ it is necessary
only to assure property (e) of Markov transformations, see \cite{M}, Lemma 12.2.
Notice also that the
condition (2) in the definition of stochastic matrices means that
$f$ preserves Lebesgue measure.

As in the case of Bernoulli shifts, the dynamical systems
$(\Sigma_A,\sigma,\nu)$ and $([0,1],f,\lan)$ are isomorphic and
also the Hausdorff dimension of a Borel subset $B\subset [0,1]$
coincides que the $\nu$-Hausdorff dimension of $\pi^{-1}(B)$.

We obtain the following result:

\begin{corollary} \label{onesidemarkovchaincode}
Let $(\Sigma_A,\sigma,\nu)$ be a subshift of
finite type whose stochastic matrix verifies that there exists $n_0$ such that all entries of $M^{n_0}$ are positive. Let $\{t_n\}$ be a non decreasing sequence of
natural numbers. Given a sequence $s_0=(i_0,i_1,\ldots)\in\Sigma_A$,
let $\widetilde{\cal W}(s_0,\{t_n\})$ be the set of sequences
$s=(m_0,m_1,\ldots)\in\Sigma_A$ such that
$$
m_n=i_0\,, \; m_{n+1}=i_1\,,\; \dots\;, \; m_{n+t_n} = i_{t_n} \;,
\qquad \text{for infinitely many $n$}.
$$
Then,

$(1)$ \ If $\sum_n p_{i_0} p_{i_0,i_1} \cdots
p_{i_{t_n-1},i_{t_n}}=\infty$, then $\nu(\widetilde{\cal W}(s_0,\{t_n\}))=1$. Besides, we have that
$$
\lim_{n\to\infty} \frac {\#\{j\le N: \ \sigma^j(s)\in C^{0,1,\dots,t_j}_{i_0,i_1,\dots,i_{t_j}}  \}}{\sum_{n=1}^N  p_{i_0} p_{i_0,i_1} \cdots p_{i_{t_n-1},i_{t_n}}} =1\,, \qquad
\hbox{for $\nu$-almost every $s\in \Sigma_A$.}
$$

$(2)$ \ If $\sum_n p_{i_0} p_{i_0,i_1} \cdots
p_{i_{t_n-1},i_{t_n}}<\infty$, then $\nu(\widetilde{\cal W}(s_0,\{t_n\}))=0$.

$(3)$ \ In any case we have that
$$
 \frac {h}{h+\overline L} \le \Dim_{\Pibf,\nu} (\widetilde{\cal W}(s_0,\{t_n\})) \le
 \min \left\{ 1, \frac{\log D}{h+\underline{L}}\right\}     \,,
$$
where $h=\sum_{i,j} p_i p_{i,j}\log  (1/p_{i,j})$,
$$
 \underline L=\liminf_{n\to\infty} \frac 1n \log \frac 1{p_{i_0}p_{i_0,i_1} \cdots p_{i_{t_{n-1},i_{t_n}}}}     \qquad \hbox{and} \qquad
\overline L=\limsup_{n\to\infty} \frac 1n \log \frac 1{p_{i_0}p_{i_0,i_1} \cdots p_{i_{t_{n-1},i_{t_n}}}} \,.
$$
\end{corollary}

\begin{proof} First, let us observe that for subshifts of finite type  $X_0^+=[0,1]$ because  ${\cal P}_0$ is a
finite partition of intervals. Also, since $f$ preserves the
Lebesgue measure $\lan$ we have that the ACIPM $\mu$, whose
existence is assured by Theorem E, coincides with $\lan$. Then,
parts (1) and (2) follow from Theorem \ref{markovcode}.
Finally any point in $[0,1]$ is an approximable point. Besides,
from Lemma \ref{justificacion} we have that $\overtau=0$ for
$x_0\in X_0^+=[0,1]$. Therefore, part (3)  follows from Remark
\ref{regularreal}, Theorem \ref{dimensioncode}, Remark \ref{Markovsi} and Proposition \ref{dimensionporarribacode}.
\end{proof}

The special properties of Bernoulli shifts allow us to get a
better upper bound for the Hausdorff dimension of the set
$\widetilde{\cal W}(s_0,\{t_n\})$ in that case. To prove it we
will use the following concentration inequality (see, for example
\cite{H}).

\

\noindent {\bf Lemma (Hoeffding's tail inequality)}. {\it Let
$(\Omega,{\cal A}, \mu)$ be a probability space and let
$X_1,\dots, X_n$ be independent copies of a bounded random variable $X$ taking values in the interval
$(a,b)$ almost surely. Then, for any $t>0$,
$$
\mu \Big[\sum_{i=1}^n  X_i - n\, E(X) \ge t \Big] \le e^{-2t^2/(n(b-a)^2)}\,.
$$}

\

\begin{theorem}
Let $(\Sigma,\sigma,\nu)$ be a Bernoulli shift, $\{t_n\}$ be a non
decreasing sequence of natural numbers and
$s_0=(i_0,i_1,\ldots)\in\Sigma$. Then
$$
 \frac {h}{h+\overline L} \le\Dim_{\Pibf,\nu} (\widetilde{\cal W}(s_0,\{t_n\})) \le
\frac {\sqrt{(h+\uL)^2+2\uL(\log \frac {\max_j p_j}{ \min_j p_j})^2}+h-\uL}
{\sqrt{(h+\uL)^2+2\uL(\log \frac {\max_j p_j}{ \min_j p_j})^2}+h+\uL} \,.
$$
where $h=\sum_{i} p_i \log  (1/p_{i})$,
$$
 \underline L=\liminf_{n\to\infty} \frac 1n \log \frac 1{p_{i_0}p_{i_1} \cdots p_{i_{t_n}}}\qquad \hbox{and} \qquad
\overline L=\limsup_{n\to\infty} \frac 1n \log \frac 1{p_{i_0}p_{i_1} \cdots p_{i_{t_n}}} \,.
$$
In particular, if $p_0=p_1=\cdots=p_{D-1}=1/D$ and $\overline L=\uL=L$, then
$$
\Dim_{\Pibf,\nu} (\widetilde{\cal W}(s_0,\{t_n\}))=  \frac {h}{h+ L}\,.
$$
\end{theorem}

\begin{proof} The lower inequality follows from Corollary \ref{onesidemarkovchaincode}. So we only need to deal with the upper one.
We define the collections
$$
{\cal F}_n = \{ C_{j_0,\dots,j_{n-1},i_0,\dots,i_{t_n}}^{0,\dots,n+t_n} : \ j_0.j_1,\dots,j_{n-1} \in
\{0,\dots,D-1\} \}\,.
$$
Then, the collection
$$
{\cal G}_N = \cup_{n=N}^\infty {\cal F}_n
$$
covers the set $\widetilde{\cal W}(s_0,\{t_n\})$ and
\begin{equation} \label{cover}
\sum_{n=N}^\infty \sum_{F\in {\cal F}_n} \nu(F)^\tau = \sum_{n=N}^\infty
\nu(C_{i_0,\dots,i_{t_n}}^{0,\dots,t_n})^\tau \sum_{j_0,\dots,j_{n-1}=1}^{D-1}
\nu(C_{j_0,\dots,j_{n-1}}^{0,\dots,n-1})^\tau \,.
\end{equation}
For each $n$, we will divide  the partition of $\Sigma $
$$
{\cal P}_{n-1}=
\{C_{j_0,\dots,j_{n-1}}^{0,\dots,n-1} : \ j_0.j_1,\dots,j_{n-1} \in
\{0,\dots,D-1\}\}
$$
in the following three subcollections:
\begin{align}
{\cal P}_{n-1,{\rm big}} & =\left\{C\in {\cal P}_{n-1}: \ \nu(C)\ge e^{-nh}  \right\} \,, \notag \\
{\cal P}_{n-1,{\rm middle}} & =\left\{C\in {\cal P}_{n-1}: \ e^{-\beta n} < \frac {\nu(C)}{e^{-nh}} < 1  \right\} \,, \notag\\
{\cal P}_{n-1,{\rm small}} & =\left\{C\in {\cal P}_{n-1}: \  \frac {\nu(C)}{e^{-nh}} \le e^{-\beta n}  \right\} \,,\notag
\end{align}
where $\beta=(1+\ep)(1-\tau)/(\alpha K)$, $0<\alpha<1$ and $\ep>0$. Then
\begin{equation} \label{big}
\sum_{C\in {\cal P}_{n-1,{\rm big}}} \nu(C)^\tau = \sum_{C\in {\cal P}_{n-1,{\rm big}}}
\frac {1}{\nu(C)^{1-\tau}} \,\nu(C) \le e^{nh(1-\tau)}
\end{equation}
and
$$
\sum_{C\in {\cal P}_{n-1,{\rm middle}}} \nu(C)^\tau  \le (\#  {\cal P}_{n-1,{\rm middle}})^{1-\tau}
\left( \sum_{C\in {\cal P}_{n-1,{\rm middle}}} \nu(C)\right)^\tau \,.
$$
Since $(\#  {\cal P}_{n-1,{\rm middle}}) \, e^{-(h+\beta)n} \le  \sum_{C\in {\cal P}_{n-1,{\rm middle}}} \nu(C)$, we deduce that
\begin{equation} \label{middle}
\sum_{C\in {\cal P}_{n-1,{\rm middle}}} \nu(C)^\tau  \le e^{n(h+\beta)(1-\tau)}  \sum_{C\in {\cal P}_{n-1,{\rm middle}}} \nu(C) \le  e^{n(h+\beta)(1-\tau)}  \,.
\end{equation}
For each $n\in\N$ we will choose an increasing sequence $\{a_{n,k}\}$, $a_{n,k} \to\infty$ as $k\to\infty$, with $a_{n,1}=\beta n$. Using this sequence we divide the collection ${\cal P}_{n-1,{\rm small}}$ in the following way:
$$
{\cal P}_{n-1,{\rm small}} = \bigcup_{k=1}^\infty {\cal P}^k_{n-1,{\rm small}}\,, \qquad
{\cal P}^k_{n-1,{\rm small}}=\{C\in {\cal P}_{n-1}: \  e^{-a_{n,k+1}} \le \frac {\nu(C)}{e^{-nh}} \le
e^{-a_{n,k}} \} \,.
$$
Then,
$$
\sum_{C\in {\cal P}_{n-1,{\rm small}}} \nu(C)^\tau = \sum_{k=1}^\infty \sum_{C\in {\cal P}^k_{n-1,{\rm small}}} \nu(C)^\tau  \le \sum_{k=1}^\infty (\#  {\cal P}^k_{n-1,{\rm small}})^{1-\tau}
\left( \sum_{C\in {\cal P}^k_{n-1,{\rm small}}} \nu(C)\right)^\tau \,.
$$
Since $(\#  {\cal P}^k_{n-1,{\rm small}}) \, e^{-nh} e^{-a_{n,k+1}} \le  \sum_{C\in {\cal P}^k_{n-1,{\rm small}}} \nu(C)$, we deduce that
\begin{equation} \label{small}
\sum_{C\in {\cal P}_{n-1,{\rm small}}} \nu(C)^\tau  \le  e^{nh(1-\tau)} \sum_{k=1}^\infty e^{(1-\tau)a_{n,k+1}} \sum_{C\in {\cal P}^k_{n-1,{\rm small}}} \nu(C)   \,.
\end{equation}
Now, observe that
$$
\bigcup_{C\in  {\cal P}^k_{n-1,{\rm small}}} C = \left\{ s\in\Sigma: \ a_{n,k}\le \log \frac 1{\nu( P(n-1,s))} -nh <
a_{n,k+1} \right\} \,.
$$
For each $j=0,1,\ldots$,
let $Z_j:\Sigma\longrightarrow \R$ be the ramdom variable defined by
$$
Z_j(i_0,i_1,\ldots)=\log \frac 1{p_{i_j}}\,.
$$
These ramdom variables are independent and identically distributed with expectated value
$$
E(Z_j)=\sum_{i=0}^{D-1} p_i \log \frac 1{p_i} = h \,.
$$
Moreover,
$$
S_n(s):=\sum_{j=0}^{n-1} Z_j(s) =  \log \frac 1{\nu( P(n-1,s))}
$$
and therefore
$$
\bigcup_{C\in  {\cal P}^k_{n-1,{\rm small}}} C = \left\{ s\in\Sigma: \ a_{n,k}\le S_n(s) - E(S_n) <
a_{n,k+1} \right\} \,.
$$
Hence, from Hoeffding's tail inequality we have that, for all $\ep>0$,
$$
\nu \Big( \bigcup_{C\in  {\cal P}^k_{n-1,{\rm small}}} C \Big) \leq e^{ - K a_{n,k}^2/n}\,, \qquad
\hbox{with} \quad K= \frac 2{\big(\log \frac {\max_j p_j}{ \min_j p_j}+\ep\big)^2}\,.
$$
Using now (\ref{small}) we get that
$$
\sum_{C\in {\cal P}_{n-1,{\rm small}}} \nu(C)^\tau  \le  e^{nh(1-\tau)} \sum_{k=1}^\infty e^{(1-\tau)a_{n,k+1}}
e^{ - K a_{n,k}^2/n} \,.
$$
Notice that for $0<\alpha<1$ the sequence defined by
$$
a_{n,k+1} = \frac {K\alpha}{1-\tau} \frac {a_{n,k}^2}n \,, \qquad a_{n,1}=\beta n \,,
$$
verifies that
$$
a_{n,k}= \frac {(1-\tau)n}{K\alpha} \Big( \frac {K\alpha}{1-\tau}\,\beta \Big)^{2^k} = \frac {(1-\tau)n}{K\alpha} (1+\ep)^{2^k} \longrightarrow \infty \,, \qquad \hbox{as } k\to\infty
$$
and therefore
$$
\sum_{C\in {\cal P}_{n-1,{\rm small}}} \nu(C)^\tau  \le  e^{nh(1-\tau)} \sum_{k=1}^\infty
e^{ - K(1-\alpha) a_{n,k}^2/n} \,.
$$
But
$$
 \sum_{k=1}^\infty e^{ - K(1-\alpha) a_{n,k}^2/n} \le  \int_{\beta n}^\infty
e^{ - K(1-\alpha) x^2/n} dx \le \frac {\Gamma(1/2)}{2\sqrt{(1-\alpha)K}}\, \sqrt{n} \,.
$$
Hence, for any $\eta>0$ and $n$ large enough we have
\begin{equation} \label{small3}
\sum_{C\in {\cal P}_{n-1,{\rm small}}} \nu(C)^\tau  \le  e^{nh(1-\tau)(1+\eta)} \,.
\end{equation}
Using now (\ref{cover}), (\ref{big}), (\ref{middle}) and (\ref{small3}) we deduce, for $N$ large enough,
$$
\sum_{n=N}^\infty \sum_{F\in {\cal F}_n} \nu(F)^\tau \le 3 \sum_{n=N}^\infty
\nu(C_{i_0,\dots,i_{t_n}}^{0,\dots,t_n})^\tau e^{n(h+\beta)(1-\tau)}\,.
$$
Since, given $\ep>0$, for $N$ large enough we have that $\nu(C_{i_0,\dots,i_{t_n}}^{0,\dots,t_n})\le e^{-n(\underline{L}(s_0)-\ep)}$ we conclude that
$$
\sum_{n=N}^\infty \sum_{F\in {\cal F}_n} \nu(F)^\tau \le 3 \sum_{n=N}^\infty
e^{-n\tau(\underline{L}(s_0)-\ep)}  e^{n(h+\beta)(1-\tau)}
\to 0 \qquad \hbox{as } N\to\infty
$$
if
$$
\tau> \frac {\sqrt{(h+\uL-\ep)^2+4(1+\ep)(\uL-\ep)/(K\alpha)}+h-\uL+\ep}{\sqrt{(h+\uL-\ep)^2+4(1+\ep)(\uL-\ep)/(K\alpha)}+h+\uL-\ep}\,.
$$
Therefore,
$$
\Dim_{\Pibf,\nu} (\widetilde{\cal W}(s_0,\{t_n\})) \le  \frac {\sqrt{(h+\uL-\ep)^2+2(1+\ep)(\uL-\ep)(\log \frac {\max_j p_j}{ \min_j p_j}+\ep)^2/\alpha}+h-\uL+\ep}
{\sqrt{(h+\uL-\ep)^2+2(1+\ep)(\uL-\ep)(\log \frac {\max_j p_j}{ \min_j p_j}+\ep)^2/\alpha}+h+\uL-\ep} \,.
$$
The result follows by taking $\ep\to 0$ and $\alpha\to 1$\,.

\end{proof}

The above results allow us to get, for example, the
following one:

\begin{corollary} \label{bernouillicode}
Let $(\Sigma,\sigma,\nu)$ be a Bernoulli shift.
\begin{itemize}
\item[\rm (1)] Let $t_n=[\log n]$. Then, for every sequence
$(\,i_0,i_1,\dots)\in \Sigma$ we have that, for $\nu$-almost
all sequence $(m_0,m_1,\dots)\in \Sigma$,
$$
m_n=i_0\,, \; m_{n+1}=i_1\,,\; \dots\;, \; m_{n+t_n} = i_{t_n} \;,
\qquad \text{for infinitely many $n$}.
$$

\item[\rm (2)] Let $t_n=[n^\kappa]$ with $\kappa>0$. Then, for
every sequence $(i_0,i_1,\dots)\in \Sigma$, the set $\widetilde{\cal W}$ of sequences
$(m_0,m_1,\dots)\in \Sigma$ such that
$$
m_n=i_0\,, \; m_{n+1}=i_1\,,\; \dots\;, \; m_{n+t_n} = i_{t_n} \;,
\qquad \text{for infinitely many $n$},
$$
has zero $\nu$-measure. Moreover, the $\nu$-grid Hausdorff dimension of $\widetilde{\cal W}$ is $1$ if $0<\kappa<1$ and zero if $\kappa>1$.
\end{itemize}
\end{corollary}

\begin{proof}
Notice that if $t_n=[n^\kappa]$, we have
$$
\sum_n p_{i_0} p_{i_1}\cdots p_{i_{t_n}} \le \sum_n (\max_j p_j )^{t_n}\sum_n
(\max_j p_j )^{n^\kappa-1} <\infty
$$
and if  $t_n=[\log n]$ we have that
$$
\sum_n p_{i_0} p_{i_1}\cdots p_{i_{t_n}} \ge
 \sum_n (\min_j p_j )^{t_n} \ge  \sum_n (\min_j p_j )^{1+\log n} = \infty \,.
$$
Also if $t_n=[n^\kappa]$ we have that $L=0$ if $0<\kappa<1$ and $L=\infty$ if $\kappa>1$.
\end{proof}

When $p_0=\cdots=p_{D-1}=1/D$ we can identify the Bernoulli shift
$(\Sigma, \sigma, \nu)$ with the set of $D$-base representations
of numbers in the interval $[0,1]$. The associated expanding map
$f$ is then the map $f(x)=Dx$ (mod $1$), and the measure results
contained in Corollary \ref{bernouillicode} in this particular
case, are well known (see \cite{Ph}).

\

\subsubsection{Gauss transformation}

Let us consider now the map $\phi:[0,1]\longrightarrow [0,1]$
given by
$$
\phi(x)=
\begin{cases}
\dfrac 1x - \left[ \dfrac 1x \right] \,, \quad &\text{if } x\ne 0
\,, \vspace{3pt}
\\
0 \,, \quad &\text{if } x=0 \,.
\end{cases}
$$
Here $[x]$ denotes the integer part of $x$. The map $\phi$ is
called the {\it Gauss transformation} and it is very close related
with the theory of continued fractions. Recall that given $0<x<1$
we can write it as
$$
x=\frac 1{n_0 + \phi(x)} \,, \qquad \text{with } n_0:= \left[\frac
1x \right] \,.
$$
If $\phi(x)\ne 0$, i.e. if $x\notin \{1/n: \ n\in\N\}\cup\{0\}$,
we can repeat the process with $\phi(x)$ to obtain
$$
x= \frac 1{n_0+\dfrac 1{n_1+\phi^2(x)}}\,, \qquad \text{with }
n_1:= \left[\frac 1{\phi(x)} \right] \,.
$$
If $\phi^n(x)\ne 0$ for all $n$, or equivalently if $x$ is
irrational, we can repeat the process for all $n$ and associate in
this way to $x$ the infinite sequence $\{n_j\}$, with
$n_j=[1/\phi^{j}(x)]$ and we write
$$
x:= [n_0 \; n_1 \; n_2\; \dots ] = \lim_{j\to\infty} \frac
1{n_0+\dfrac 1{n_1+\dfrac 1{n_2+\dfrac 1{\ddots +n_j}}}}\;.
$$
Observe that if we denote by $I_n$ the interval
$I_n=(1/(n+1),1/n)$, then the sequence $n_j$ is determined by the
property $\phi^{j}(x) \in I_{n_j}$.

If $x$ is rational the above expansion is finite (ending with $n$
such that $\phi^n(x)=0$. We call to the code $[n_0 \; n_1 \;
n_2\; \dots ]$ the continued fraction expansion of $x$.
It is clear that the Gauss transformation acts on the continued
fraction expansions as the left shift
$$
x=[n_0 \; n_1 \; n_2\; \dots ] \qquad \implies \qquad \phi(x)= [
n_1 \; n_2\;  \dots ]\,.
$$

It is not difficult to check that the Gauss transformation $\phi$
is a Markov transformation with respect to the partition ${\cal P}_0=\{I_n\}$ and that
the continued fraction expansion of $x$ coincide with the code associated to
an expanding map given in Section \ref{codes}.
It is also easy to check that $\phi$ preserves the so called
Gauss measure which is defined by
$$
\mu(A)= \frac 1{\log 2} \int_A \frac 1{1+x} \, d\lambda(x)
$$
where $\lan$ denotes Lebesgue measure. Since this measure is
obviously absolutely continuous with respect to $\lan$, we
conclude that the Gauss measure is the unique $\phi$-invariant
absolutely continuous probability whose existence is assured by
Theorem E.

\

The next theorem is an example of the kind of statements  that  we
can obtain when we apply our results to  the Gauss transformation.

\

\begin{corollary} \label{gaussradios}
\begin{itemize}
\item[]
\item[\rm (1)] If $\alpha>1$ then, for almost all $x_0\in [0,1]$, and more precisely,
if $x_0=[\,i_0\,,\;i_1\,, \; \dots ]$ is an irrational number such that
$\log i_n = o(n)$ as $n\to\infty$, we have that
$$
\liminf_{n\to\infty} n^{1/\alpha}  |\phi^n(x)-x_0| = 0 \,, \qquad
\text{for almost all $x\in [0,1]$.}
$$

\item[\rm (2)] If $\alpha <1$, then for all $x_0\in [0,1]$ we have that
$$
\liminf_{n\to\infty} n^{1/\alpha} |\phi^n(x)-x_0| = \infty \,,
\qquad \text{for almost all $x\in [0,1]$.}
$$

\item[\rm (3)] If $x_0$ verifies the same hypothesis than in part $(1)$, then
$$
\Dim \left\{x\in [0,1]: \ \liminf_{n\to\infty} n^{1/\alpha}
|\phi^n(x)-x_0| =0 \right\} = 1\,, \qquad \text{for any
$\alpha>0$}.
$$
and
$$
\Dim \left\{x\in [0,1]: \ \liminf_{n\to\infty} e^{n\kappa}
|\phi^n(x)-x_0| =0 \right\} \ge \frac {\pi^2}{\pi^2+6\kappa\log2}
\,, \qquad \text{for any $\kappa>0$}.
$$
\end{itemize}
\end{corollary}

\begin{proof}
Let us observe first that now $\diminf=\dimsup=1$ for all $x_0\in
[0,1]$ and that obviously $\lan$ and $\mu$ are comparable in
$[0,1]$. With this facts in mind, part (2) is a consequence of part (2) of Theorem \ref{markovradios}
if $x_0\in \cup_j I_j$. Part (2) is also true if $x_0=1/m$ for
some $m\in\N$, since $\lan$ and $\mu$ are comparable in $[0,1]$
and then we do not need that $B(x_0,r_k)\subset P(0,x_0)$ in the
proof of Proposition \ref{serieconverge}.

Since $T(P)=(0,1)$ for all $P\in {\cal P}_0$ we can use Proposition
\ref{elchachi} for the case $j=n+1$ to get that
$$
\frac {\lan(P(n,x_0))}{\lan(P(n+1,x_0))} \asymp \frac
{1}{\lan(P(0,T^{n+1}(x_0)))}\,.
$$
But $T^{n+1}(x_0)=[\,i_{n+1}\,,\;i_{n+2}\,, \; \dots ]$ and
therefore $P(0,T^{n+1}(x_0))=(1/(i_{n+1}+1),1/i_{n+1})$. Hence
$$
\log \frac {\lan(P(n,x_0))}{\lan(P(n+1,x_0))} \asymp \log i_{n+1}
$$
and   we conclude that $\tau(x_0)=0$ if $\log i_n=o(n)$ as
$n\to\infty$. Part (1) follows now from
Corollary \ref{hartos2}, since in this case the set $X_0$ is
precisely the set of irrational numbers in $[0,1]$.

Since $\lan$ and $\mu$ are comparable in $[0,1]$ we have that all
irrational number is approximable (see Definition
\ref{aproximable}) and as we have just seen $\tau(x_0)=0$ if $\log
i_n=o(n)$, we can use  Remark
\ref{regularreal} and Corollary \ref{cotasparadimconradios2} to
obtain that
$$
\Dim \left\{x\in [0,1]: \ \liminf_{n\to\infty}
\frac{|\phi^n(x)-x_0|} {r_n} =0 \right\} \ge \frac {h}{h+\ell}
$$
for any non increasing sequence $\{r_n\}$ of positive numbers such
that there exists $\ell:=\lim_{n\to\infty} \frac 1n \log \frac
1{r_n}$. Here $h$ denotes de entropy of the Gauss transformation
which is known to be
$$
h=\frac 2{\log 2} \int_0^1 \frac {\log (1/x)}{1+x}\, dx = \frac
{\pi^2}{6\log 2}\;.
$$
Part (3) follows now from the fact that if $r_n=n^{-1/\alpha}$
with $\alpha>0$ then $\ell=0$, and if $r_n=e^{-n\kappa}$ with
$\kappa>0$ we have $\ell=\kappa$.
\end{proof}

For continued fractions expansions there is an analogous to Corollary \ref{bernouillicode}. However
we have preferred to state the following result involving the digits appearing in the continued fraction expansion of $x_0$.

%

%



\begin{corollary}
Let $x_0\in [0,1]$ be an irrational number with continued fraction expansion $x_0=[\,i_0,i_1,\dots \,]$ and let $t_n$ be a non decreasing sequence of natural numbers. Let $\widetilde W$ be the set of
points $x=[\,m_0,m_1,\dots\,]\in [0,1]$ such that
$$
m_n=i_0\,, \; m_{n+1}=i_1\,,\; \dots\;, \; m_{n+t_n} = i_{t_n} \;,
\qquad \text{for infinitely many $n$}.
$$

\begin{itemize}
\item[\rm (1)] \ $\lan(\widetilde W)=1$, if
$$
 \sum_n \frac 1{(i_0+1)^2\cdots (i_{t_n}+1)^2} =\infty \,.
$$

\item[\rm (2)] \ $\lan(\widetilde W)=0$, if
$$
\sum_n \frac 1{i_0^2\cdots i_{t_n}^2} <\infty \,.
$$

\item[\rm (3)] In any case, if  $\log i_n = o(n)$ as $n\to\infty $, then
$$
\Dim (\widetilde W) \ge \frac h{h+\limsup_{n\to\infty}
\frac 1n \log  (i_0+1)^2\cdots (i_{t_n}+1)^2} \,.
$$
\end{itemize}
\end{corollary}

\begin{proof} It is easy to check that, for all $n\in\N$,
$$
\frac 1{(i_0+1)^2\cdots (i_{n}+1)^2} \le \lan (P(n,x_0)) \le  \frac 1{i_0^2\cdots i_{n}^2}\,.
$$
Then, parts (1) and (2) follow from Theorem \ref{markovcode}  and part (3) is a consequence of Theorem \ref{dimensioncode} and Remark \ref{regularreal}.
\end{proof}

\subsection{Inner functions}

A Blaschke product is a complex function of the type
$$
B(z)= \prod_{k=1}^\infty \frac {|a_k|}{a_k} \frac{z-a_k}{1-\overline{a_k}z}\,, \qquad \hbox{$|a_k|<1$}\,,
$$
verifying the Blaschke condition $\sum_{k=1}^\infty (1-|a_k|)<\infty$. The function $B(z)$ is holomorphic in the unit disk $\D=\{z: |z|<1\}$ of the complex plane and it is an example of an inner function, i.e a holomorphic function $f$
defined on $\D$ and with values in $\D$ whose radial limits
$$
f^*(\xi):=\lim_{r\to 1^-}  f(r\xi)
$$
(which exists for almost every $\xi$ by Fatou's Theorem) have
modulus $1$ for almost every $\xi\in\p\D$.
Therefore an inner function $f(z)$ induces a mapping $f^*:\p\D\longrightarrow \p\D$.
It is well known that any inner function can be written as
$$
f(z) = e^{i\phi} B(z) \, \exp
\left( -\int_{\p\D} \frac{\xi+z}{\xi-z} \,d\nu(\xi) \right)
$$
where $B(z)$ is a Blaschke product and $\nu$ is a finite positive singular measure on $\p\D$.

\smallskip

For inner functions it is well known the following result, see e.g. \cite{R}:

\smallskip

\noindent {\bf Theorem G (L\"owner's lemma).} {\it If $f:\D
\longrightarrow \D$ is an inner function then $f^*:\p\D
\longrightarrow \p\D$  preserves Lebesgue measure if and only if
$f(0)=0$.}

\smallskip

We recall that, by the Denjoy-Wolff theorem \cite{D}, for any
holomorphic function $f:\D\longrightarrow \D$ which is not
conjugated to a rotation, there exists a point $p\in\overline\D$,
the so called Denjoy-Wolff point of $f$, such the iterates $f^n$
converge to $p$ uniformly on compact subsets of $\D$. Also, if
$p\in\D$ then $f(p)=p$ and if $p\in\p\D$ then $f^*(p)=p$. Hence,
if $f$ is an inner function which is not conjugated to a rotation
and does not have a fixed point $p\in\D$ then its Denjoy-Wolff
point $p$ belongs to $\p\D$ and $f^n$ converges to $p$ uniformly
on compact subsets of $\D$.
Bourdon, Matache and Shapiro \cite{BMS} and Poggi-Corradini \cite{PC} have
proved independently that if $f$ is inner with a fixed point in
$p\in\p\D$, then $(f^*)^n$ can converge to $p$ for almost every
point in $\p\D$. In fact, see Theorem 4.2  in \cite{BMS}, $(f^*)^n \to p$
almost everywhere in $\p\D$ if and only if $\sum_n
(1-|f^n(0)|)<\infty$.

If $f$ is inner with a fixed point in $\D$, $f$ preserves the
harmonic measure $\omega_p$. We recall that $\omega_p$ can be
defined as the unique probability measure such that, for all
continuous function $\phi:\p\D\longrightarrow \R$,
$$
\int_{\p\D} \phi \, d\omega_p = \widetilde \phi (p) \,,
$$
where $\widetilde\phi$ is the unique extension of $\phi$ which is
continuous in $\overline\D$ and harmonic in $\D$. It follows that
if $A$ is an arc in $\p\D$, then $\omega_p(A)$ is the value at the
point $p$ of the harmonic function whose radial limits take the
value $1$ on $A$ and the value $0$ on the exterior of $A$.

If $f$ is inner with a fixed point in $\D$, but it is not
conjugated to a rotation, J. Aaronson \cite{A} and J.H. Neuwirth \cite{N}
proved, independently, that $f^*$ is exact with respect to
harmonic measure and therefore mixing and ergodic. In fact, inner
functions are also ergodic with respect to $\alpha$-capacity
\cite{FPR}.  An interesting study of some dynamical properties of inner
functions is contained in the works of M. Craizer. In \cite{C1} he
proves that if $f'$ belongs to the Nevanlinna class, then the
entropy of $f^*$ is finite and it can be
calculated by the formula
$$
h(f^*) = \frac 1{2\pi} \int_0^{2\pi} \log |(f^*)'(x)|\, dx\,,
$$
where $(f^*)'$ denotes the angular derivative of $f$. He also
proves that the Rohlin invertible extension of an inner function with a fixed
point in $\D$ is equivalent to a generalized Bernoulli shift, see
\cite{C2}.

The mixing properties of inner functions are even stronger. In
this sense Ch. Pommerenke \cite{P} has shown the following

\

\noindent {\bf Theorem H (Ch. Pommerenke).} {\it Let
$f:\D\longrightarrow\D$ be an inner function with $f(0)=0$,  but
not a rotation. Then, there exists a positive absolute constant
$K$ such that
$$
\left| \frac{\lan[B\cap (f^*)^{-n}(A)]}{\lan(A)} - \lan(B) \right| \le
K\,e^{-\alpha n} \,,
$$
for all $n\in \N$, for all arcs $A,B\subset\p\D$, where
$\alpha=\max\{1/2, |f'(0)|\}$ and $\lan$ denotes normalized Lebesgue measure.}

\

In the terminology of \cite{FMP} this imply that inner functions with
$f(p)=p$  $(p\in\D)$ are uniformly mixing at any point of $\p\D$ with respect to the harmonic measure $\omega_p$. In
particular, we have that the correlation coefficients of
characteristic functions {\it of balls} have exponential decay. As a consequence
of Theorem 3 in \cite{FMP}, and the arguments of the proofs of Corollaries \ref{hartosinf}
and \ref{hartos2} we have  that if $\xi_0$ is any point in $\p\D$ and
$\{r_n\}$ is a non increasing sequence of positive numbers, then we have that

\begin{itemize}
\item[(A)] If $\sum_{n=1}^\infty r_n < \infty$, then
$$
\liminf_{n\to\infty} \frac {d((f^*)^n(\xi), \xi_0)}{r_n} = \infty\,, \qquad \hbox{for almost every
$\xi\in\p\D$} \,.
$$
\item[(B)] If $\sum_{n=1}^\infty r_n = \infty$, then
$$
\lim_{N\to\infty} \dfrac{\#\{n\le N: \  d((f^*)^n(\xi), \xi_0)<r_n
\}}{\sum_{n=1}^{N}r_n} = 1\,, \qquad \hbox{for almost every
$\xi\in\p\D$} \,.
$$
and
$$
\liminf_{n\to\infty} \frac {d((f^*)^n(\xi), \xi_0)}{r_n} =0 \,,
\qquad \hbox{for almost every $\xi\in\p\D$}\,.
$$
\end{itemize}

A finite Blaschke product $B$ (with, say, $N$ factors) is a rational function of degree
$N$  and therefore it is a covering of order $N$ of $\p\D$. As a consequence $B$ has a fixed point in $\p\D$ if $N\ge 3$ or if $N=2$ and $B(0)=0$. Hence, we can choose a branch of the argument of $B(e^{i\theta})$ mapping $0$ on $0$ and $[0,2\pi]$ onto $[0,2N\pi]$. Also $B(z)$ is $C^\infty$ at the boundary $\p\D$ of the unit disk and its derivative verifies
$$
|B'(z)| = \sum_{k=1}^N \frac {1-|a_k|^2}{|z-a_k|^2} \,,
\qquad \hbox{if $|z|=1$}\,.
$$
Therefore, if $B(0)=0$, we have that $|B'(z)|>C>1$ for all
$z\in\p\D$, and the dynamic of $B^*$ on $\p\D$ is isomorphic to
the dynamic of a Markov transformation with a finite partition
${\cal P}_0$ (it has $N$ elements) and having the Bernoulli
property. Besides, since the Lebesgue measure is exact we have
that the ACIPM measure of the system is precisely Lebesgue measure
$\lan$. Hence, we obtain the following improvement of statement
(A):

\begin{theorem} \label{innerp}
Let $B:\D \longrightarrow \D$ be a finite Blaschke product with a
fixed point $p\in \D$, but not an automorphism which is conjugated
to a rotation. Let also $\xi_0$ be any point in $\p\D$ and let
$\{r_n\}$ be a non increasing sequence of positive numbers. Then
$$
\Dim \left\{ \xi\in\p\D: \liminf_{n\to\infty} \frac
{d((B^*)^n(\xi), \xi_0)}{r_n} =0\right\} \ge \frac {h}{h+\ellsup}
$$
where  $\ellsup=\limsup_{n\to\infty} \frac 1n \log \frac 1{r_n}$, $h=
\int_{\p\D} \log |B'(z)|\, d\lan(z)$ and $\Dim$ denotes Hausdorff
dimension\,. The result is sharp in the sense that we get equality when $B(z)=z^N$ and
$\ellsup=\ell= \lim_{n\to\infty} \frac 1n \log \frac 1{r_n}$.
\end{theorem}

\begin{proof}
In the case that $p=0$ the result follows from the above comments and Theorem \ref{markovradios}. In the general case, let $T:D\longrightarrow \D$ be a M\"obius transformation such that $T(p)=0$. Then, $g=T\circ B\circ T^{-1}$ is a finite Blasckhe product with $g(0)=0$. Besides, it is easy to see that
$$
\{\xi\in\p\D : \  d((g^*)^n(\xi),\xi_0)<r_n \ \hbox{i.o.}\} \subseteq T(\{\xi\in\p\D : \  d((B^*)^n(\xi),T^{-1}(\xi_0))<Cr_n \ \hbox{i.o.}\} )
$$
where $C$ is a constant depending on $T$. Therefore the lower bound follows from the case $p=0$. The equality for $B(z)=z^N$ follows from Proposition \ref{dimensionporarriba}.
\end{proof}

Theorem \ref{innerp} is also true for the following infinite Blaschke product:
$$
B(z)= \prod_{k=0}^\infty \frac{z-a_k}{1-a_kz}\,, \qquad \hbox{$a_k=1-2^{-k}$}\,.
$$
since as we will see, the dynamic of $B^*$ on $\p\D$ is isomorphic
to the dynamic of a Markov transformation with a countable
partition ${\cal P}_0$ and with the Bernoulli property. Notice
also that $B^*$ is exact with respect to Lebesgue measure  and
therefore we have that the ACIPM measure is Lebesgue measure.

For this Blaschke product $B$ is defined in $\p\D\setminus\{1\}$ and in fact it is
$C^\infty$ there and
$$
|B'(z)| = \sum_{k=0}^\infty \frac {1-a_k^2}{|z-a_k|^2} \,,
\qquad \hbox{if $|z|=1, \ z\ne1$}\,.
$$
If we denote $B(e^{2\pi it})=e^{2\pi iS(t)}$ then $S'(t)=|B'(e^{2\pi it})|>C>1$. Moreover,
it follows from Phragm\'en-Lindel\"of Theorem that the image of $S(t)$ is $(-\infty,\infty)$
and so we can define the intervals $P_j= \{t\in (0,1): \ j<S(t)<j+1\}$. The transformation
$T:[0,1]\longrightarrow [0,1]$ given by $T(t)=S(t)$ (mod $1$), $T(0)=T(1)=0$, is a Markov transformation with partition ${\cal P}_0=\{P_j\}$. To see this we only left to prove property (f). We define the following collection of subarcs of $\p\D$: $I^+_k=\{e^{i\alpha}:
\theta_{k+1}<\alpha<\theta_{k}\}$ $(k\ge 0)$, where $\theta_{0}=\pi$ and for each $k\ge 1$ we denote by $e^{i\theta_k}$ $(\theta_k\in (0,\pi))$ the point whose distance to $1$ is $1-a_{k-1}=2^{-(k-1)}$. We define also $I_k^-=\{z\in\p\D: \bar z \in I_k^+\}$. It is geometrically clear that  if $z\in I^\pm_j$, then $|\sin 2\pi t|\le C\,2^{-j}$ and also that
$$
|z-a_k| \ge
\begin{cases}
C\, 2^{-j} \,, \qquad \hbox{for $k\ge j$} \\
C\, 2^{-k} \,, \qquad \hbox{for $k < j$}\,,
\end{cases}
$$
Now, if $z=e^{2\pi it}\in I^\pm_j$, we have that
\begin{gather}
S'(t)=|B'(e^{2\pi it})|= \sum_{k=0}^\infty \frac {1-a_k^2}{|e^{2\pi it}-a_k|^2}\ge C\, \frac {1-a_j}{2^{-2j}}=C\, 2^j \,, \notag \\
S'(t) =  \sum_{k=0}^\infty \frac {1-a_k^2}{|e^{2\pi it}-a_k|^2} \le C \sum_{k<j} \frac {2^{-k}}{2^{-2k}} + C \sum_{k\ge j} \frac {2^{-k}}{2^{-2j}} \le C\, 2^j
\notag
\end{gather}
and
$$
|S''(t)| \le C \sum_{k=0}^\infty \left| \frac {a_k(1-a_k^2)\sin 2\pi t}{(e^{2\pi it}-a_k)^4} \right|
\le C \sum_{k<j} \frac {2^{-k}2^{-j}}{2^{-4k}} + C \sum_{k\ge j} \frac {2^{-k}2^{-j}}{2^{-4j}}
\le C\, 2^{2j}\,.
$$
Therefore, since $\lan(I_j^\pm)\asymp 2^{-j}$ we have that $\int_{I_j^\pm} S'(t) \, dt \asymp C$ and so each $P_j\in {\cal P}_0$ contains at most a fixed constant number of consecutive intervals $I_k^\pm$. Hence, there exists an absolute constant $C$ such that, if $t_1,t_2,t_3 \in P_j$, then
$$
\frac {|T''(t_1)|}{T'(t_2)\,T'(t_3)} \le C
$$
and this implies that $T$ verify property (f) of Markov transformations.

Finally, the entropy $h$ of $B^*$ (or $T(t)$) is finite, because
\begin{align}
h=\int_{\p\D} \log |B'(z)|\, d\lan(z) &=2\sum_{j=0}^\infty  \int_{I_j}
\log |B'(z)|\, d\lan(z) \notag \\
&\le 2 \sum_{j=0}^\infty \log\Big(C\sum_{k=0}^\infty \frac {2^{-k}}{2^{-2j}}\Big) \frac 1{2^{j+1}} = 2 \sum_{j=0}^\infty \log\Big(C 2^{2j}\Big) \frac 1{2^{j+1}}
<\infty\,. \notag
\end{align}

The singular inner functions
$$
f(z)=e^{c\frac {1+z}{1-z}} \,, \qquad \hbox{for $c<-2$}.
$$
also verify Theorem \ref{innerp}. These inner functions have only
one singularity at $z=1$ and its Denjoy-Wolff point $p$ is real
and it verifies $0<p<1$. It is easy to see that if $f(e^{2\pi
it})= e^{2\pi i S(t)}$ for $t\in [0,1]$, then $S(t) = \frac
{c}{2\pi} \, \cot \pi t$. and the dynamic of $f^*$ on $\p\D$ is
isomorphic to the dynamic of the Markov transformation $T(t)=S(t)$
(mod 1). We have that the partition ${\cal P}_0$ for $T$ is
countable, ${\cal P}_0=\{P_j: j\in\Z\}$ where $P_j = \{t\in (0,1):
\ j<S(t)<j+1\}$, and $T$ has the Bernoulli property, i.e.
$T(P_j)=(0,1)$. Notice also that  $T'(t)=\frac {|c|}2 \,\csc^2 \pi
t>1$ and that, for $x,y\in P_j$,
$$
\left| \frac{T'(x)}{T'(y)}-1 \right| = \frac {|T(x)+T(y)|}{T(y)^2+ (c/2\pi)^2}\, |T(x)-T(y)|\le
\frac {2j+2}{j^2+ (c/2\pi)^2}\, |T(x)-T(y)|\le C\, |T(x)-T(y)|\,.
$$
It is known that the entropy of $f$ is finite (see \cite{Ma})
$$
h(f)= \log \left( \frac {1}{1-p^2} \log \frac 1{p^2} \right)<\infty\,.
$$

More generally, Theorem \ref{innerp} holds for inner functions $f$ with a fixed point $p\in \D$ and finite entropy such that the transformation $T$ defined as in these examples is Markov. This happens, for example, if the set of singularities of $f$ in $\p\D$ is finite,
the lateral limits of $f^*$ at the singular points are $\pm\infty$
and $T$ verifies properties (d) and (f) of Markov transformations. Notice that the condition on the lateral limits holds, for example, for Blasckhe products whose singular set is finite and each singular point $\xi$ is an accumulation point of zeroes inside of a Stolz cone with vertex $\xi$.
However we think that Theorem \ref{innerp} is true for any inner function  with a fixed point $p\in \D$ and finite entropy.

\subsection{Expanding endomorphisms} \label{expend}

Let $M$ be a compact Riemannian manifold. A $C^1$ map
$f:M\longrightarrow M$ is an {\it expanding endomorphism} if there
exists a natural number $n\ge 1$ and constants $C>0$ and $\beta>1$ such that
$$
\|(D_x f^n)u\| > C\,\beta^n \|u\|\,, \qquad \text{for all } x\in M, u\in T_x
M\,.
$$
A $C^1$ expanding endomorphism of a compact connected Riemannian
manifold $M$ whose derivative $D_x f$ is a H\"older continuous
function of $x$ is an expanding map with respect to Lebesgue
measure $\lan$ and a finite Markov partition ${\cal P}_0$, see \cite{M},
p.171. Therefore, the unique $f$-invariant probability measure
whose existence is assured by Theorem E is comparable to $\lan$ in
the whole $M$.  Our results also apply for this dynamical system.





%

\end{document}